\documentclass[a4paper,11pt,leqno]{article}
\pdfoutput=0
\usepackage[T1]{fontenc}	
\usepackage[utf8]{inputenc}
\usepackage[italian, english]{babel}
\usepackage{amssymb}
\usepackage{graphicx}
\usepackage{enumerate}
\usepackage{amsthm}
\usepackage[toc,page]{appendix}
\usepackage{amsmath}
\usepackage{makeidx}
\theoremstyle{plain}
\newtheorem{thm}{Theorem}

\newtheorem{lem}[thm]{Lemma}
\newtheorem{prop}[thm]{Proposition}
\theoremstyle{definition}
\newtheorem{defn}[thm]{Definition}

\newtheorem*{remark}{Remark}
\theoremstyle{remark}
{\left\lbrace\begin{array}{@{}l@{}}}%
{\end{array}\right.}
\usepackage{geometry} 
\usepackage{fullpage}
\usepackage{eucal}
\usepackage[babel]{csquotes}
\usepackage[numbers]{natbib}
\numberwithin{equation}{section}
\numberwithin{thm}{section}

\usepackage[hang,flushmargin]{footmisc} 

\usepackage{setspace}

\usepackage[textsize=small]{todonotes}
\newcommand{\tinytodo}[2][]
   {\todo[caption={#2}, size=\small, #1, disable]{\renewcommand{\baselinestretch}{0.5}\selectfont#2\par}}

\title{Global existence and asymptotics for quasi-linear one-dimensional Klein-Gordon equations with mildly decaying Cauchy data}

\author{A. Stingo\\
Universit{\'e} Paris 13,\\
Sorbonne Paris Cité, LAGA, CNRS (UMR 7539),\\
99, Avenue J.-B. Cl{\'e}ment,\\
F-93430 Villetaneuse}

\date{2015}

\begin{document}

\nocite{*}

\maketitle

\begin{abstract}
\singlespace
Let u be a solution to a quasi-linear Klein-Gordon equation in one-space dimension, $\Box u + u = P(u, \partial_t u, \partial_x u; \, \partial_t \partial_x u, \partial^2_x u)$, where $P$ is a homogeneous polynomial of degree three, and with smooth Cauchy data of size $\varepsilon \rightarrow 0$.
It is known that, under a suitable condition on the nonlinearity, the solution is global-in-time for compactly supported Cauchy data.
We prove in this paper that the result holds even when data are not compactly supported but just decaying as $\langle x \rangle^{-1}$ at infinity, combining the method of Klainerman vector fields with a semiclassical normal forms method introduced by Delort.
Moreover, we get a one term asymptotic expansion for $u$ when $t \rightarrow +\infty$.
\end{abstract}

\let\thefootnote\relax\footnotetext{ 
The author is supported by a PhD fellowship funded by the FSMP and the Labex MME-DII.}

\section*{Introduction}
\renewcommand{\theequation}{\arabic{equation}}
The goal of this paper is to prove the global existence and to study the asymptotic behaviour of the solution $u$ of the one-dimensional nonlinear Klein-Gordon equation, when initial data are small, smooth and slightly decaying at infinity.
We will consider the case of a quasi-linear cubic nonlinearity, namely a homogeneous polynomial $P$ of degree 3 in $(u, \partial_tu, \partial_xu; \, \partial_t\partial_x u , \partial^2_x u)$, affine in $(\partial_t\partial_x u , \partial^2_x u)$, so the initial valued problem is written as
\begin{equation} \label{KG introduction}
\begin{cases}
& \Box u + u = P(u, \partial_t\partial_xu, \partial^2_xu\,; \partial_tu,\partial_xu) \\
& u(1,x)=\varepsilon u_0(x) \\
& \partial_tu(1,x)=\varepsilon u_1(x)
\end{cases}  \qquad t\ge 1 \, , x \in \mathbb{R} \, , \varepsilon \in ]0,1[ \,.
\end{equation}

\vspace{0.5cm}
\noindent Our main concern is to obtain results for data which have only mild decay at infinity (i.e. which are $O(|x|^{-1})$, $x\rightarrow +\infty$), while most known results for quasi-linear Klein-Gordon equations in dimension 1 are proved for compactly supported data.
In order to do so, we have to develop a new approach, that relies on semiclassical analysis, and that allows to obtain for Klein-Gordon equations results of global existence making  use of Klainerman vector fields and usual energy estimates, instead of $L^2$ estimates on the hyperbolic foliation of the interior of the light cone, as done for instance in an early work of Klainerman \cite{klainerman:global_existence} and more recently in the paper  of LeFloch, Ma \cite{LeFloch_Ma}.

\vspace{0.5cm}
\noindent We recall first the state of the art of the problem.
In general, the problem in dimension 1 is critical, contrary to the problem in higher dimension which is subcritical. 
In fact, in space dimension $d$, the best time decay one can expect for the solution is $\|u(t,\cdot)\|_{L^{\infty}}=O(t^{-\frac{d}{2}})$: therefore, in dimension 1 the decay rate is $t^{-\frac{1}{2}}$, and for a cubic nonlinearity, depending for example only on $u$, one has $\|P(u)\|_{L^2}\le Ct^{-1}\|u(t,\cdot)\|_{L^2}$, with a time factor $t^{-1}$ just at limit of integrability.
It is well known from works of Klainerman \cite{klainerman:global_existence} and Shatah\tinytodo{Je ne crois pas que le résultat
  de Shatah soit limité au cas à support compact. Vérifier et adapter si besoin cette phrase}~\cite{shatah:normalForms} that the analogous problem in space dimension $d\ge 3$ has global-in-time solutions if $\varepsilon$ is sufficiently small.
In \cite{klainerman:global_existence}, Klainerman proved it for smooth, compactly supported initial data, with nonlinearities at least quadratic, using the Lorentz invariant properties of $\Box + 1$ to derive uniform decay estimates and generalized energy estimates for solutions $u$ to linear inhomogeneous Klein-Gordon equations. 
Simultaneously, in \cite{shatah:normalForms} Shatah proved this result for smooth and integrable initial data, extending Poincaré's theory of normal forms for ordinary differential equations to the case of nonlinear Klein-Gordon equations.
An earlier work from Simon \cite{simon:wave_operator}, and from Simon, Taflin in \cite{simontaflin:wave_operators} for coupled Klein-Gordon equations with several masses, established the global existence for data given at $t=\infty$.
In \cite{hormander:non-linear}, H\"{o}rmander refined Klainerman's techniques to obtain new time decay estimates of solutions to linear inhomogeneous Klein-Gordon equations and he showed that, for quadratic nonlinearities, the solution exists over $[-T_{\varepsilon},T_{\varepsilon}]$ with an existence time $T_{\varepsilon}$ such that $\lim_{\varepsilon\rightarrow 0} \varepsilon \log T_{\varepsilon}=\infty$ when $d=2$, while $\lim_{\varepsilon \rightarrow 0}\varepsilon^2 T_{\varepsilon}= \infty$ for $d=1$.
In addition, he presented two conjectures: for quadratic nonlinearities, $T_{\varepsilon}=\infty$ in two space dimensions, while for space dimension one $\liminf_{\varepsilon\rightarrow 0}\varepsilon^2 \log T_{\varepsilon} >0$. 
The first conjecture has been proved by Ozawa, Tsutaya and Tsutsumi in \cite{ozawa:global-existence} in the semi-linear case, after partial results by Georgiev, Popivanov in \cite{georgiev:global_solutions}, Kosecki \cite{kosecki:unit_condition} (for nonlinearities verifying some "suitable null conditions"), and Simon, Taflin in \cite{simon:cauchy_problem}. 
Later, in \cite{ozawaTT:remarks} Ozawa, Tsutaya and Tsutsumi announced the extension of their proof to the quasi-linear case and studied scattering of solutions.
Still in dimension 2, Delort, Fang and Xue proved in \cite{delortFangXue:global_existence} the global existence of solutions for a quasi-linear system of two Klein-Gordon equations, with masses $m_1, m_2$, $m_1\ne 2m_2$ and $m_2 \ne 2 m_1$ for small, smooth, compactly supported Cauchy data, extending the result proved by Sunagawa in \cite{sunagawa:global_smallamplitude} in the semi-linear case. Moreover, they proved that the global existence holds true also when $m_1=2m_2$ and a convenient null condition is satisfied by nonlinearities. 
The same result in the resonant case is also proved by Katayama, Ozawa \cite{katayama:null_condition}, and by Kawahara, Sunagawa \cite{kawahara:KGsystems}, in which the structural condition imposed on nonlinearities includes the Yukawa type interaction, which was excluded from the \emph{null condition} in the sense of \cite{delortFangXue:global_existence}.
In this context, we cite also the work of Germain \cite{germain:global_existence}, and of Ionescu, Pausader \cite{Ionescu_Pausader}, for a system of coupled Klein-Gordon equations with different speeds in dimension 3, with a quadratic nonlinearity, respectively in the semilinear case for the former, and in the quasi-linear one for the latter. For data small, smooth, and localized, they prove that a global solution exists and scatters.

\vspace{0.5cm}
\noindent In dimension 1, Moriyama, Tonegawa and Tsutsumi \cite{moriyama:almost_global} have shown that the solution exists on a time interval of length longer or equal to $e^{c/\varepsilon^2}$, where $\varepsilon$ is the Cauchy data's size, with a nonlinearity vanishing at least at order three at zero, or semi-linear. 
They also proved that the corresponding solution asymptotically approaches  the free solution of the Cauchy problem for the linear Klein-Gordon equation.
The fact that in general the solution does not exist globally in time was proved by Yordanov in \cite{yordanov:blow-up}, and independently by Keel and Tao \cite{keeltao:small_data}. 
However, there exist examples of nonlinearities for which the corresponding solution is global-in-time: on one hand, if $P$ depends only on $u$ and not on its derivatives; on the other hand, for seven special nonlinearities considered by Moriyama in \cite{moriyama:normal_forms}.
A natural question is then posed by H\"{o}rmander, in \cite{hormander:lifespan, hormander:non-linear}: can we formulate a structure condition for the nonlinearity, analogous to the null condition introduced by Christodoulou \cite{christodoulou:global_solutions} and Klainerman \cite{klainerman:null_condition} for the wave equation, which implies global existence?
In \cite{delort:existence_global, delort:erratum} Delort proved that, when initial data are compactly supported, one can find a \emph{null condition}, under which global existence is ensured.
This condition is likely optimal, in the sense that when the structure hypothesis is violated, he constructed in \cite{delort:minoration} approximate solutions blowing up at $e^{A/\varepsilon^2}$, for an explicit constant $A$.
This suggests that also the exact solution of the problem blows up in time at $e^{A/\varepsilon^2}$, but this remains still unproven.

\vspace{0.5cm}
\noindent In most of above mentioned papers dealing with the one dimensional problem, two key tools are used: normal forms methods and/or Klainerman vector fields $Z$.
In particular, the latter are useful since they have good properties of commutation with the linear part of the equation, and
their action on the nonlinearity $ZP(u)$ may be expressed from $u, Zu$ using Leibniz rule.
This allows one to prove easily energy estimates for $Z^ku$ and then to deduce from them $L^{\infty}$ bounds for $u$, trough Klainerman-Sobolev type inequalities.
However, in these papers the global existence is proved assuming small, \emph{compactly supported} initial data. This
\tinytodo{J'ai rajouté cette phrase}  is
related to the fact that the aforementioned authors use in an essential way a change of variable in hyperbolic coordinates,
that does not allow for non compactly supported Cauchy data.
Our aim is to extend the result of global existence for cubic quasi-linear nonlinearities in the case of small compactly supported Cauchy data of \cite{delort:existence_global, delort:erratum}, to the more general framework of data with mild polynomial decay.
To do that, we will combine the Klainerman vector fields' method with the one introduced by Delort in \cite{delort:semiclassical}. 

\vspace{0.5cm}
\noindent In \cite{delort:semiclassical}, Delort  develops a semiclassical normal form method to study global existence for nonlinear hyperbolic equations with small, smooth, decaying Cauchy data, in the critical regime and when the problem does not admit Klainerman vector fields.
The strategy employed is to construct, through semiclassical analysis, some \emph{pseudo-differential} operators which commute with the linear part of the considered equation, and which can replace vector fields when combined with a microlocal normal form method.
Our aim here is to show that one may combine these ideas together with the use of Klainerman vector fields to obtain, in one dimension, and for nonlinearities satisfying the null condition, global existence and modified scattering.

\vspace{0.5cm}
\noindent In our paper, we prove the global existence of the solution $u$ by a \emph{boostrap} argument, namely by showing that we can propagate some suitable \emph{a priori} estimates made on $u$.
We propagate two types of estimates: some energy estimates on $u$, $Zu$, and some uniform bounds on $u$.
To prove the propagation of energy estimates is the simplest task.
We essentially write an energy inequality for a solution $u$ of the Klein-Gordon equation in the quasi-linear case (the main reference is the book of H\"{o}rmander \cite{hormander:non-linear}, chapter 7), and then we use the commutation property of the Klainerman vector fields $Z$ with the linear part of the equation to derive an inequality also for $Zu$. 
Moreover, $Z$ acts like a derivation on the nonlinearity, so the Leibniz rule holds and we can estimate $ZP$ in term of $u, Zu$.
Injecting \emph{a priori} estimates in energy inequalities and choosing properly all involved constants allow us to obtain the result.

\begin{figure}[ht]
\centering
\includegraphics[height=5cm , scale=0.65]{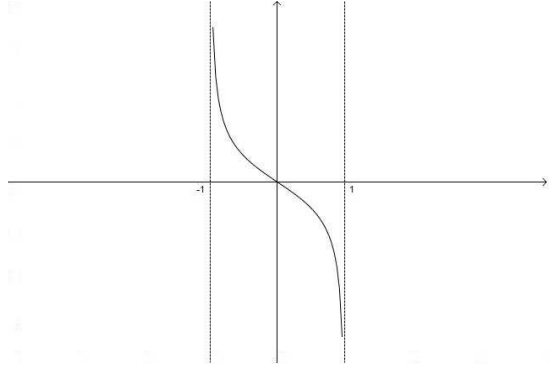}
\caption{$\Lambda$ for the Klein Gordon equation. } \label{fig: Lambda portrait}
\end{figure}

\vspace{0.5cm}
\noindent The
\tinytodo{Modifications effectuées dans ce paragraphe}
 main difficulty is to prove that the uniform estimates hold and can be propagated. Actually, as mentioned
above, the one dimensional Klein-Gordon equation is critical, in the sense that the expected decay for
$\|u(t,\cdot)\|^2_{L^{\infty}}$  is in $t^{-1}$, so is not integrable. A drawback of that is that one cannot prove energy
estimates that would be uniform as time tends to infinity. Consequently, a Klainerman-Sobolev inequality, that would
control $\|u(t,\cdot)\|_{L^{\infty}}$ by $t^{-1/2}$ times the $L^2$ norms of $u, Zu$, would not give the expected optimal
  $L^\infty$-decay of the solution, but only a bound in $t^{-\frac{1}{2}+\sigma}$ for some positive $\sigma$, which is
  useless to close the bootstrap argument. 
The idea to overcome this difficulty is, following   the approach of Delort in \cite{delort:semiclassical}, to rewrite
\eqref{KG introduction} in semiclassical coordinates, for some new unknown function $v$. The goal is then to deduce from the
PDE satisfied by $v$ an ODE from which one will be able to get a uniform $L^\infty$ bound for $v$ (which is equivalent to
the optimal $t^{-1/2}$ $L^\infty$-decay of $u$).  
Let us describe our approach for a simple model of Klein-Gordon equation:
\tinytodo{J'ai rajouté des paramètres $\alpha, \beta, \gamma, \delta$ pour faire apparaître explicitement   la conditon nulle sur cet exemple.}
\begin{equation} \label{model equation}
(D_t - \sqrt{1+D_x^2})u = \alpha u^3 + \beta|u|^2u + \gamma|u|^2\bar{u} + \delta\bar{u}^3,
\end{equation}
where $\alpha, \beta, \gamma, \delta$ are constants, $\beta$ being \emph{real} (this last assumption reflecting the null condition on that example).
Performing a semiclassical change of variables and unknowns $u(t,x)=\frac{1}{\sqrt{t}}v(t,\frac{x}{t})$, we rewrite this equation as
\begin{equation} \label{model equation on v}
[D_t - Op_h^w(\lambda_h(x,\xi))] v= h(\alpha v^3+ \beta|v|^2v + \gamma|v|^2\bar{v} + \delta\bar{v}^3) \,,
\end{equation}
where $\lambda_h(x,\xi)=x\xi+\sqrt{1+\xi^2}$, the semiclassical parameter $h$ is defined as $h:=1/t$, and the Weyl
quantization of a symbol $a$ is given by
\begin{equation*}
Op_h^w(a)v = \frac{1}{2\pi h}\int_{\mathbb{R}}\int_{\mathbb{R}} e^{\frac{i}{h}(x-y)\xi}a\big(\frac{x+y}{2},\xi\big)v(y) \, dyd\xi \, .
\end{equation*}
One introduces the manifold $\Lambda= \{(x,\xi) \, | \, x + \frac{\xi}{\sqrt{1+\xi^2}}=0\}$ as in figure \ref{fig: Lambda portrait}, which is the graph of the smooth function $d\varphi(x)$, where $\varphi:
]-1,1[\rightarrow \mathbb{R}$ is $\varphi(x)=\sqrt{1-x^2}$.
\noindent One can deduce an ODE from \eqref{model equation on v}, developing the symbol $\lambda_h(x,\xi)$ on $\Lambda$, i.e. on $\xi = d\varphi(x)$.
One obtains a first term $a(x)$ independent of $\xi$ and a remainder, which turns out to be integrable in time as may be shown using some ideas of Ifrim-Tataru \cite{ifrimtataru:global_bounds} and the $L^2$ estimates verified by $v$ and by the action of the Klainerman vector field on $v$.
Essentially, one uses that operators whose symbols are localised in a neighbourhood of $\Lambda$ of size $O(\sqrt{h})$ have $\mathcal{L}(L^2; L^{\infty})$ norms which are a $O(h^{-\frac{1}{4}-\sigma})$, for a small positive $\sigma$, instead of $O(h^{-\frac{1}{2}})$ as would follow from a direct application of the Sobolev injection.
In this way, one proves that $v$ is solution of the equation
\begin{equation}
D_t v = a(x) v + h \beta|v|^2v + \mbox{non characteristic terms} + \mbox{remainder of higher order in }h \, .
\end{equation}
Then the idea is to eliminate \emph{non characteristic} terms by a normal forms argument, introducing a new function $f$ which will be finally solution of an ordinary differential equation
\begin{equation}
D_t f = a(x) f + h \beta |f|^2f + \mbox{remainder of higher order in }h \,.
\end{equation}
From this equation, one easily derives an uniform control $L^{\infty}$ on $f$, and then on the starting solution $u$.

\vspace{0.5cm}
\noindent The
\tinytodo{J'ai effectué des modifications dans ce paragraphe.}
 analysis of the above  ODE provides as well a one term asymptotic expansion of the solution of
equation~\eqref{model equation} (or, more generally of the solution~\eqref{KG introduction}), as proved in the last section of this paper. This expansion shows that, in
general, scattering does not hold, and that one has only modified scattering. This is in contrast with higher dimensional
problems for the Klein-Gordon equation, where
 global solutions have at infinity the same behaviour as free solutions.
In space dimension one, only few results were known regarding asymptotics of solution, including for the simpler equation
$$\Box u + u = \alpha u^2  +\beta u^3 + \mbox{order} \, 4 \,.$$
For this equation, Georgiev and Yordanov \cite{georgievYordanov:asymptotic} proved that, when $\alpha=0$, the distance between the solution $u$ and linear solutions cannot tend to 0 when $t\rightarrow \infty$, but they do not obtain an asymptotic description of the solution (except for the particular case of sine-Gordon $\Box u + \sin u =0$, for which they use methods of "nonlinear scattering").
In \cite{Lindblad_Soffer:scattering}, Lindblad and Soffer studied the scattering problem for long range nonlinearities, proving that for all prescribed asymptotic solutions there is a solution of the equation with such behavior, for some choice of initial data, and finding the complete asymptotic expansion of the solutions.
Their method is based on the reduction of the long range phase effects to an ODE, via an appropriate ansatz. 
In \cite{hayashi:quadratic_kleingordon}, a sharp asymptotic behaviour of small solutions in the quadratic, semilinear case is proved by Hayashi and Naumkin, without the condition of compact support on initial data, using the method of normal forms of Shatah.
In \cite{delort:existence_global}, Delort studied asymptotics in the quasi-linear case, obtaining a one term asymptotic expansion for the solution, under the assumption of small, compactly supported Cauchy data, and showing that in general the solution does not behave as in the linear case.
The only other cases in dimension one for which the asymptotic behaviour is known concern nonlinearities studied by Moriyama
in \cite{moriyama:normal_forms}, where he showed that solutions have a free asymptotic behaviour, assuming the initial data to be sufficiently small and decaying at infinity.
\tinytodo[inline]{Est-ce dans le
  cas de données à support compact ou pas ? Rajouter aussi une référence à mon article, en précisant que le comportement
  asymptotique n'y est obtenu que pour des données à support compact.}

\section{Statement of the main results}
\renewcommand{\theequation}{\thesection.\arabic{equation}}

The Cauchy problem we are considering is
\begin{equation} \label{KG}
\begin{cases}
& \Box u + u = P(u, \partial_t\partial_xu, \partial^2_xu\,; \partial_tu,\partial_xu) \\
& u(1,x)=\varepsilon u_0(x) \\
& \partial_tu(1,x)=\varepsilon u_1(x)
\end{cases}  \qquad t\ge 1 \, , x \in \mathbb{R}
\end{equation}
where $\Box:= \partial^2_t - \partial^2_x$ is the D'Alembert operator, $\varepsilon \in ]0,1[$, $u_0, u_1$ are smooth enough functions.
$P$ denotes a homogeneous polynomial of degree three, with real constant coefficients, affine in $(\partial_t\partial_xu, \partial^2_xu)$. 
We can highlight this particular dependence on second derivatives following the approach of \cite{delort:existence_global} and decomposing $P$ as 
\begin{equation}
P(u, \partial_t\partial_xu, \partial^2_xu\,; \partial_tu,\partial_xu) = P'(u\, ; \partial_tu, \partial_x u) + P''(u, \partial_t\partial_xu, \partial^2_xu\,; \partial_tu,\partial_xu) \,,
\end{equation}
where $P', P''$ are homogeneous polynomials of degree three, $P''$ linear in $(\partial_t\partial_xu, \partial^2_xu)$.
Moreover
\begin{equation} \label{polynomial decomposition}
\begin{split}
P'(X_1 ; \,Y_1, Y_2)&= \sum_{k=0}^3 i^k P_k'(X_1 ; -i Y_1, -i Y_2) \\
 P''(X_1, X_2, X_3 ;\, Y_1, Y_2) & = \sum_{k=0}^2 i^k P_k''(X_1, -X_2, -X_3; -i Y_1, -i Y_2)
 \end{split}
\end{equation}
where $P_k'$ is homogeneous of degree $k$ in $(Y_1, Y_2)$ and of degree $3-k$ in $X_1$, while $P_k''$ is homogeneous of degree 1 in $(X_2, X_3)$ and of degree $k$ in $(Y_1, Y_2)$.
We denote $P_k=P'_k+P''_k$.
For $x \in ]-1,1[$, define 
\begin{equation}
\begin{split}
& \omega_0(x) := \frac{1}{\sqrt{1-x^2}} \,, \\
& \omega_1(x) := \frac{-x}{\sqrt{1-x^2}} \, ,
\end{split}
\end{equation}
and 
\begin{equation} \label{Phi null condition}
\Phi(x) := P'_1(1; \omega_0(x), \omega_1(x)) + P''_1(1, \omega_0(x)\omega_1(x), \omega_1^2(x); \omega_0(x), \omega_1(x)) + 3 P'_3(1; \omega_0(x), \omega_1(x)) \, .
\end{equation}
\begin{defn} \label{definition of null condition}
We say that the nonlinearity $P$ satisfies the \emph{null condition} if and only if $\Phi \equiv 0$.
\end{defn}

\vspace{0.5cm}
\noindent Our goal is to prove that there is a global solution of \eqref{KG} when $\varepsilon$ is sufficiently small, $u_0, u_1$ decay rapidly enough at infinity, and when the cubic nonlinearity satisfies the \emph{null condition}. 
We state the main theorem below.

\begin{thm}[Main Theorem] \label{main theorem}
Suppose that the nonlinearity $P$ satisfies the null condition.
Then there exists an integer $s$ sufficiently large, a positive small number $\sigma$, an $\varepsilon_0 \in ]0,1[$ such that, for any real valued $(u_0, u_1)\in H^{s+1}(\mathbb{R})\times H^s(\mathbb{R})$ satisfying
\begin{equation}
\|u_0\|_{H^{s+1}}+ \|u_1\|_{H^s}+ \|xu_0\|_{H^2}+\|xu_1\|_{H^1}\le 1 \,, \label{conditions on data}
\end{equation}
for any $0<\varepsilon<\varepsilon_0$, the problem \eqref{KG} has an unique solution $u \in C^0([1,+\infty[; H^{s+1})\cap C^1([1,+\infty[; H^s)$.
Moreover, there exists a $1$-parameter family of continuous function $a_{\varepsilon}: \mathbb{R}\rightarrow \mathbb{C}$, uniformly bounded and supported in $[-1,1]$, a function $(t,x)\rightarrow r(t,x)$ with values in $L^2(\mathbb{R})\cap L^{\infty}(\mathbb{R})$, bounded in $t\ge 1$, such that, for any $\varepsilon \in ]0,\varepsilon_0]$, the global solution $u$ of \eqref{KG} has the asymptotic expansion
\begin{equation} \label{asymptotics for u}
u(t,x) = \Re\left[\frac{\varepsilon}{\sqrt{t}}a_{\varepsilon}\left(\frac{x}{t}\right)\exp\left[it \varphi\left(\frac{x}{t}\right)+i\varepsilon^2\left|a_{\varepsilon}\left(\frac{x}{t}\right)\right|^2\Phi_1\left(\frac{x}{t}\right)\log t\right] \right] + \frac{\varepsilon}{t^{\frac{1}{2}+\sigma}}r(t,x) \,, 
\end{equation}
where $\varphi(x)= \sqrt{1-x^2}$, and
\begin{equation} \label{def of Phi1}
\begin{split}
\Phi_1(x)  = \frac{1}{8}\langle \omega_0(x)\rangle^{-4} & \left[3P_0(1, \omega_0(x)\omega_1(x),\omega_1(x)^2;\,\omega_0(x),\omega_1(x)) \right.\\
&\left. + P_2(1, \omega_0(x)\omega_1(x),\omega_1(x)^2;\,\omega_0(x),\omega_1(x))\right]\, ,
\end{split}
\end{equation}
with $\langle x \rangle = \sqrt{1+x^2}$.
\end{thm}

\vspace{0.5cm}
\noindent We denote by $Z$ the Klainerman vector field for the Klein-Gordon equation, that is
$Z:= t\partial_x + x\partial_t$, and by $\Gamma$ a generic vector field in the set $\mathcal{Z}=\{Z, \partial_t, \partial_x\}$.  
The most remarkable properties of these vector fields are the commutation with the linear part of the equation in \eqref{KG}, namely
\begin{equation} \label{commutation property Gamma e Box}
[\Box + 1, \Gamma]=0 \, ,
\end{equation}
and the fact that they act like a derivation on the cubic nonlinearity.
Using the notation $D=\frac{1}{i}\partial$, we also denote by $W^{t, \rho,\infty}$ a modified Sobolev space, made by functions $t\rightarrow\psi(t,\cdot)$ defined on an interval, such that $\langle D_x\rangle^{\rho-i}D^i_t u \in L^{\infty}$, for $i \le 2$, with the norm
\begin{equation}
\|\psi(t,\cdot)\|_{W^{t, \rho,\infty}(\mathbb{R})}:= \sum_{i=0}^2\| \langle D_x\rangle^{\rho-i}D^i_t \psi(t,\cdot)\|_{L^{\infty}(\mathbb{R})}\, .
\end{equation}

\vspace{0.5cm}
\noindent The proof of the main theorem is based on a \emph{bootstrap} argument.
In other words, we shall prove that we are able to propagate some \emph{a priori} estimates made on a solution $u$ of \eqref{KG} on some interval $[1,T]$, for some $T>1$ fixed, as stated in the following theorem.

\begin{thm}[Bootstrap Theorem] \label{bootstrap theorem}
There exist two integers $s, \rho$ large enough, $s \gg \rho$, an $\varepsilon_0 \in ]0,1[$ sufficiently small, and two constants $A, B >0$ sufficiently large such that, for any $0<\varepsilon<\varepsilon_0$, if $u$ is a solution of \eqref{KG} on some interval $[1,T]$, for $T>1$ fixed, and satisfies
\begin{subequations} \label{bootstrap hypothesis}
\begin{align}
& \|u(t,\cdot)\|_{W^{t, \rho, \infty}}  \le A \varepsilon t^{-\frac{1}{2}} \\
& \|Zu(t,\cdot)\|_{H^1}  \le B\varepsilon t^{\sigma} \,,\qquad \|\partial_tZu(t,\cdot)\|_{L^2}\le B\varepsilon t^{\sigma} \label{energy bound1}\\
& \|u(t,\cdot)\|_{H^s} \le B\varepsilon t^{\sigma} \,,\qquad \hspace{0.3cm} \|\partial_tu(t,\cdot)\|_{H^{s-1}} \le B\varepsilon t^{\sigma} \label{energy bound2}\,,
\end{align}
\end{subequations}
for every $t \in [1,T]$, for some $\sigma\ge 0$ small, then it verifies also
\begin{subequations} \label{bootstrap thesis}
\begin{align}
& \|u(t,\cdot)\|_{W^{t, \rho, \infty}}  \le \frac{A}{2} \varepsilon t^{-\frac{1}{2}}\label{A' bootstrap} \\
& \|Zu(t,\cdot)\|_{H^1}  \le \frac{B}{2}\varepsilon t^{\sigma} \,, \qquad \|\partial_tZu(t,\cdot)\|_{L^2}\le \frac{B}{2}\varepsilon t^{\sigma}\\
& \|u(t,\cdot)\|_{H^s} \le \frac{B}{2}\varepsilon t^{\sigma} \,,\qquad \hspace{0.3cm}\|\partial_tu(t,\cdot)\|_{H^{s-1}} \le \frac{B}{2}\varepsilon t^{\sigma} \,.
\end{align}
\end{subequations}
\end{thm}

\noindent In section \ref{section Generalised energy estimates} we show that energy bounds \eqref{energy bound1}, \eqref{energy bound2} can be propagated, simply recalling an energy inequality obtained by H\"{o}rmander in \cite{hormander:non-linear} for a solution $u$ of a quasi-linear Klein-Gordon equation, and applying it to $\partial^{s-1}_xu$ and $Zu$.
Sections from \ref{Semiclassical Pseudo-differential Operators.} to \ref{Study of the ODE and End of the Proof} concern instead the proof of the uniform estimate's propagation.
Furthermore, in section \ref{Study of the ODE and End of the Proof} we derive also the asymptotic behaviour of the solution $u$.

\vspace{0.5cm}
\noindent To conclude, we can mention that we will mainly focus on not very high frequencies, for it is easier to control what happens for very large frequencies which correspond to points on $\Lambda$ in figure \ref{fig: Lambda portrait} close to vertical asymptotic lines.
This is justified by the fact that contributions of frequencies of the solution larger than $h^{-\beta}$, for a small positive $\beta$, have $L^2$ norms of order $O(h^N)$ if $s\beta \gg N$, assuming small $H^s$ estimates on $v$.
In this way, most of the analysis is reduced to frequencies lower than $h^{-\beta}$.

\section{Generalised energy estimates} \label{section Generalised energy estimates}

\noindent With notations introduced in the previous section, we define 
\begin{equation}
E_0(t,u)=\big( \|\partial_tu(t,\cdot)\|_{L^2}^2+\|\partial_xu(t,\cdot)\|_{L^2}^2 + \|u(t,\cdot)\|_{L^2}^2 \big)^{1/2}
\end{equation}
as the square root of the energy associated to the solution $u$ of \eqref{KG} at time $t$, and $E^{\Gamma}_N(t,u)=\displaystyle\sum_{k=0}^N\big(E_0(t,\Gamma^ku)^2\big)^{1/2}$, for a fixed $\Gamma$.
The goal of this section is to obtain an energy inequality involving $E^{\Gamma}_N(t,u)$.
In particular, since the aim is to propagate \emph{a priori} energy bounds on $u$, i.e. $\|u(t,\cdot)\|_{H^s}$, $\|\partial_tu(t,\cdot)\|_{H^{s-1}}$, $\|Zu(t,\cdot)\|_{H^1}$ and $\|\partial_tZu(t,\cdot)\|_{L^2}$, we will consider on one hand $E^{\partial_x}_{s-1}(t,u)$ where all $\Gamma$ are equal to $\partial_x$, and on the other $E^Z_1(t,u)$ where $\Gamma=Z$. 
Often in what follows we will denote partial derivatives with respect to $t$ and $x$ respectively by $\partial_0$ and $\partial_1$.

\vspace{0.5cm}
\noindent We will use the following result, which concerns the specific energy inequality for the Klein-Gordon equation in the quasi-linear case, and which is presented here without proof (see lemma 7.4.1 in \cite{hormander:non-linear} for further details). 
\begin{lem} \label{Hormander energy lemma}
Let u be a solution of 
\begin{equation} 
\Box u + u + \gamma^{01}\partial_0\partial_1u + \gamma^{11} \partial_1^2 u + \gamma^0\partial_0 u + \gamma^1 \partial_1 u = f \, ,
\end{equation}
where functions $\gamma^{ij}=\gamma^{ij}(t,x)$, $\gamma^j=\gamma^{j}(t,x)$ are smooth, such that $\displaystyle\sum_{i,j=0}^1|\gamma^{ij}|+|\gamma^{j}|\le \frac{1}{2}$.
Then,
\begin{equation} \label{energy inequality Hormander}
E_0(t,u)  \le C \big[ E_0(1,u) + \int_1^t \|f(\tau,\cdot)\|_{L^2} d\tau \big] \exp \big( \int_1^t C(\tau)\,d\tau \big) \, ,
\end{equation}
where $C(\tau):= \displaystyle\sum_{i,j,h=0}^1 \displaystyle\sup_{x} \left( |\partial_h\gamma^{ij}(\tau,x)|+|\partial_h\gamma^{j}(\tau,x)|\right) $.
\end{lem}

\vspace{0.5cm}
\noindent We can rewrite the equation in \eqref{KG} in the same form as in lemma \ref{Hormander energy lemma}, especially highlighting the linear dependence on second derivatives,
\begin{equation} \label{eq: KG in quasi linear structure}
\Box u + u + \gamma^{01}\partial_0\partial_1u + \gamma^{11} \partial_1^2 u + \gamma^0\partial_0 u + \gamma^1 \partial_1 u = 0 \, ,
\end{equation}
where coefficients $\gamma^{ij}, \gamma^j$ are homogeneous polynomials of degree two in $(u, \partial_0 u, \partial_1 u)$.
Let us apply $\partial^{s'}_1$, $s':=s-1$, to this equation.
If $u$ is a solution of \eqref{eq: KG in quasi linear structure}, then $\partial^{s'}_1u$ satisfies
\begin{equation}
\Box \partial^{s'}_1u + \partial^{s'}_1u + \partial^{s'}_1 \big(\gamma^{01}\partial_0\partial_1u + \gamma^{11} \partial_1^2 u + \gamma^0\partial_0 u + \gamma^1 \partial_1 u \big)=0 \, ,
\end{equation}
and applying the Leibniz rule, we obtain that $\partial^{s'}_1 u$ is solution of the equation
\begin{equation}
\Box \partial^{s'}_1u + \partial^{s'}_1u +  \gamma^{01}\partial_0\partial_1(\partial^{s'}_1u) + \gamma^{11} \partial_1^2 (\partial^{s'}_1u) + \gamma^0\partial_0 (\partial^{s'}_1u) + \gamma^1 \partial_1 (\partial^{s'}_1u) = f^{s'} \, ,
\end{equation}
where $f^{s'}$ is a linear combination of terms of the form
\begin{equation}
\begin{split}
& (\partial_1^{s'_1}\partial_i^{\alpha_1} u)\, (\partial_1^{s'_2}\partial_j^{\alpha_2}u)\, (\partial_1^{s'_3}\partial^2_{ij}u )\ , \\
& (\partial_1^{s'_1}\partial_i^{\alpha_1} u)\, (\partial_1^{s'_2}\partial_j^{\alpha_2}u)\, (\partial_1^{s'_3}\partial_h u )\,,
\end{split}
\end{equation}
for $i,j,h, \alpha_1, \alpha_2=0,1$, $s'_1+s'_2 + s'_3=s'$, $s'_3<s'$.
So taking the $L^2$ norm and observing that at most one index $s'_j$ can be larger than $s'/2$, we have
\begin{equation}
\|f^{s'}(t,\cdot)\|_{L^2}\le \Big(\sum_{\substack{i+j=0 \\ j\le 2}}^{[\frac{s'}{2}]+2}\|\partial^i_x\partial^j_tu(t,\cdot)\|^2_{L^{\infty}}\Big) E^{\partial_1}_{s'}(t, u) \le \|u(t,\cdot)\|^2_{W^{t,\rho,\infty}} E^{\partial_1}_{s'}(t, u)\,,
\end{equation}
for a $\rho \ge [\frac{s'}{2}]+3$.
Rewriting inequality \eqref{energy inequality Hormander} for $\partial^{s'}_1 u$, where $s'=s-1$ and $C(\tau) \le \|u(\tau,\cdot)\|^2_{W^{t,2,\infty}}$, we obtain
\begin{equation} \label{energy estimate for Es partial1}
E^{\partial_1}_{s-1}(t, u) \le C \left[E^{\partial_1}_{s-1}(1,u) + \int_1^t\|u(\tau,\cdot)\|^2_{W^{t, \rho,\infty}}E^{\partial_1}_{s-1}(\tau,u) d\tau\right] \exp\left(\int_1^t \|u(\tau,\cdot)\|^2_{W^{t,2,\infty}} d\tau\right) .
\end{equation}

\vspace{0.5cm}
\noindent On the other hand, we want to obtain an analogous of \eqref{energy estimate for Es partial1} for $E_1^Z(t,u)$.
Applying $Z$ to \eqref{eq: KG in quasi linear structure}, Leibniz rule and commutations, we derive that $Zu$ is solution of the equation
\begin{equation}
\Box Zu + Zu + \gamma^{01}\partial_0\partial_1 Zu + \gamma^{11} \partial_1^2 Zu + \gamma^0\partial_0 Zu + \gamma^1 \partial_1 Zu = f^Z \, ,
\end{equation}
where $f^Z$ is linear combination of $[\gamma^{ij}\partial^2_{ij}, Z]u$ and $[\gamma^h\partial_h, Z]u$.
We calculate for instance the term $[\gamma^{01}\partial^2_{01}, Z]u$ and we find that it is equal to $-(Z\gamma^{01})\partial^2_{01}u -\gamma^{01}[\partial^2_{01},Z]u$, that is a linear combination of  
\begin{equation}
\begin{split}
& (\partial^{\alpha_1}_iu)(\partial^{\alpha_2}_jZu)(\partial^2_{01}u) \, , \\
& (\partial^{\alpha_1}_i u)(\partial^{\alpha_2}_j u)(\partial^2_{hk}u) \,,
\end{split}
\end{equation}
for $i,j,h,k,\alpha_1, \alpha_2 = 0,1$.
Therefore, the $L^2$ norm of $f^Z$ can be estimated as follows
\begin{equation}
\|f^Z(t,\cdot)\|_{L^2}\le \Big(\sum_{i+j=0}^2\|\partial^i_x\partial^j_t u(t,\cdot)\|^2_{L^{\infty}} \Big) E^Z_1(t,u)\le \|u(t,\cdot)\|^2_{W^{t,3,\infty}}E^Z_1(t,u) \,, 
\end{equation}
and applying lemma \ref{Hormander energy lemma} for $Zu$, we derive
\begin{equation} \label{energy estimate for EZ}
E^Z_1(t,u) \le C \left[E^Z_1(1,u) + \int_1^t \|u(\tau,\cdot)\|^2_{W^{t,3,\infty}}E_1^Z(\tau,\cdot)\, d\tau \right] \exp\left(\int_1^t \|u(s,\cdot)\|^2_{W^{t,2,\infty}} ds\right)\, .
\end{equation}
\begin{remark}
  To make the above proof fully correct, one should chek as well that the energy of $Zu$ is actually finite at every fixed
  positive time. One may do that either using that the vector field $Z$ is the infinitesimal generator of the action on the
  equation of a one parameter group, along the lines of appendix~A.2 in \cite{Alazard-Delort}. Alternatively, one may instead
  exploit finite propagation speed, remarking that if the data are cut off on a compact set, the solution remains compactly
  supported at every fixed time, so that the energy of $Zu$ is actually finite, and that the bounds we get are uniform in
  terms of the cut off.
\end{remark}
\begin{prop}[Propagation of Energy Estimates] \label{Propagation of Energy Estimates}
There exist an integer $s$ large enough, a $\rho\ge [\frac{s-1}{2}]+3$, $\rho \ll s$, an $\varepsilon_0 \in ]0,1[$ sufficiently small, a small $\sigma\ge 0$, and two constants $A, B >0$ sufficiently large such that, for any $0<\varepsilon<\varepsilon_0$, if $u$ is a solution of \eqref{KG} on some interval $[1,T]$, for $T>1$ fixed, and satisfies 
\begin{subequations} \label{a priori energy est}
\begin{align}
& \|u(t,\cdot)\|_{W^{t, \rho, \infty}}  \le A \varepsilon t^{-\frac{1}{2}} \,,\label{uniform bound in energy estimate}\\
& E_{s-1}^{\partial_1}(t,u)  \le B\varepsilon t^{\sigma}\,, \label{bound on Zu}\\
& E_1^Z(t,u) \le B\varepsilon t^{\sigma}\,, \label{bound on Hs norm of u}
\end{align}
\end{subequations}
for every $t \in [1,T]$, then it verifies also
\begin{subequations} 
\begin{align}
& E_{s-1}^{\partial_1}(t,u)  \le \frac{B}{2}\varepsilon t^{\sigma} \,,\label{bound on Zu propagated}\\
& E_1^Z(t,u) \le \frac{B}{2}\varepsilon t^{\sigma}\,. \label{bound on Hs norm of u propagated}
\end{align}
\end{subequations}
\proof
Both estimates \eqref{bound on Zu} and \eqref{bound on Hs norm of u} can be propagated injecting \emph{a priori} estimates \eqref{a priori energy est} in energy inequalities \eqref{energy estimate for Es partial1} and \eqref{energy estimate for EZ} derived before, obtaining
\begin{equation*}
\begin{split}
E^{\partial_1}_{s-1}(t,u) & \le C\Big[E^{\partial_1}_{s-1}(1,u) + A^2B\varepsilon^3 \int_1^t \tau^{-1+\sigma} d\tau \Big] \le CE^{\partial_1}_{s-1}(1,u) + \frac{A^2BC\varepsilon^3}{\sigma}t^{\sigma} \,, \\
E^Z_1(t,u) & \le C\Big[E^Z_1(1,u) + A^2B\varepsilon^3 \int_1^t \tau^{-1+\sigma}  d\tau \Big]\le C E^Z_1(1,u) + \frac{A^2BC\varepsilon^3}{\sigma}t^{\sigma}  \,.
\end{split}
\end{equation*}
Then we can choose $B>0$ sufficiently large such that $C E^{\partial_1}_{s-1}(1,u) + CE_1^Z(1,u)\le \frac{B}{4}\varepsilon$, and $\varepsilon_0>0$ sufficiently small such that $\frac{A^2 C\varepsilon^2}{\sigma}\le \frac{1}{4}$, to obtain \eqref{bound on Zu propagated}, \eqref{bound on Hs norm of u propagated}.
\end{prop}

\section{Semiclassical Pseudo-differential Operators.} \label{Semiclassical Pseudo-differential Operators.}

As told in the introduction, in order to prove an $L^{\infty}$ estimate on $u$ and on its derivatives we need to reformulate the starting problem \eqref{KG} in term of an ODE satisfied by a new function $v$ obtained from $u$, and this will strongly use the semiclassical pseudo-differential calculus.
In the following two subsections, we introduce this semiclassical environment, defining classes of symbols and operators we shall use and several useful properties, some of which are stated without proof.
More details can be found in \cite{dimassi:spectral} and \cite{zworski:semiclassical}.

\subsection{Definitions and Composition Formula}

\begin{defn}
An order function on $\mathbb{R}\times\mathbb{R}$ is a smooth map from $\mathbb{R}\times\mathbb{R}$ to $\mathbb{R}_+$ : $(x,\xi)\rightarrow M(x,\xi)$ such that there exist $N_0\in \mathbb{N}$, $C>0$ and for any $(x,\xi), (y,\eta)\in \mathbb{R}\times\mathbb{R}$
\begin{equation} \label{def ineq order function}
M(y,\eta) \le C \langle x-y\rangle^{N_0} \langle\xi-\eta\rangle^{N_0} M(x, \xi) \, ,
\end{equation}
where $\langle x\rangle=\sqrt{1+x^2}$.
\end{defn}
\noindent Examples of order functions are $\langle x\rangle$, $\langle\xi\rangle$, $\langle x\rangle \langle\xi\rangle$.

\begin{defn}
Let \emph{M} be an order function on $\mathbb{R}\times\mathbb{R}$, $\beta\ge 0$, $\delta \ge 0$. One denotes by $S_{\delta, \beta}(M)$ the space of smooth functions
\begin{align*}
(x,\xi, h)  & \rightarrow a(x,\xi, h) \\
\mathbb{R}\times\mathbb{R}\times ]0,1] & \rightarrow \mathbb{C}
\end{align*}
satisfying for any $\alpha_1, \alpha_2, k, N \in \mathbb{N}$ bounds
\begin{equation} \label{symbol in S delta beta M}
|\partial_x^{\alpha_1}\partial_{\xi}^{\alpha_2}(h\partial_h)^k a(x,\xi, h)| \le C M(x,\xi)\, h^{-\delta(\alpha_1+\alpha_2)}(1+\beta h^{\beta}|\xi|)^{-N}\, .
\end{equation}
\end{defn}

\noindent A key role in this paper will be played by symbols $a$ verifying \eqref{symbol in S delta beta M} with $M(x,\xi)= \langle \frac{x+f(\xi)}{\sqrt{h}} \rangle^{-N}$, for $N \in \mathbb{N}$ and a certain smooth function $f(\xi)$. This function $M$ is no longer an order function because of the term $h^{-\frac{1}{2}}$ but nevertheless we continue to keep the notation $a \in S_{\delta,\beta}(\langle \frac{x+f(\xi)}{\sqrt{h}} \rangle^{-N})$.

\begin{defn}
We will say that $a(x,\xi)$ is a symbol of order $r$ if $a \in S_{\delta,\beta}(\langle \xi \rangle^r)$, for some $\delta \ge 0$, $\beta \ge 0$.
\end{defn}

\noindent Let us observe that when $\beta>0$, the symbol decays rapidly in $h^{\beta}|\xi|$, which implies the following inclusion for $r\ge 0$
\begin{equation}
S_{\delta,\beta}(\langle\xi\rangle^r) \subset h^{-\beta r}S_{\delta, \beta}(1) \, ,
\end{equation}
which will be often use in all the paper. 
This means that, up to a small loss in $h$, this type of symbols can be always considered as symbols of order zero.
In the rest of the paper we will not indicate explicitly the dependence of symbols on $h$, referring to $a(x,\xi,h)$ simply as $a(x,\xi)$.

\begin{defn}
Let $a\in S_{\delta, \beta}(M)$ for some order function $M$, some $\delta\ge 0$, $\beta \ge 0$. 
\begin{enumerate}[(i)]
\item We can define the \emph{Weyl quantization} of $a$ to be the operator $Op_h^w(a)=a^{w}(x, hD)$ acting on $u \in \mathcal{S}(\mathbb{R})$ by the formula :
\begin{equation}
Op_h^w(a(x,\xi))u(x) = \frac{1}{2\pi h}\int_{\mathbb{R}}\int_{\mathbb{R}} e^{\frac{i}{h}(x-y)\xi} a(\frac{x+y}{2}, \xi)\, u(y)\; dy d\xi \, ;
\end{equation} 
\item We define also the \emph{standard quantization} :
\begin{equation}
Op_h(a(x,\xi))u(x) = \frac{1}{2 \pi h} \int_{\mathbb{R}}\int_{\mathbb{R}} e^{\frac{i}{h}(x-y)\xi} a(x, \xi)\, u(y)\; dy d\xi \, .
\end{equation}
\end{enumerate}
\end{defn}

\noindent It is clear from the definition that the two quantizations coincide when the symbol does not depend on $x$.

\vspace{0.5cm}
\noindent We introduce also a semiclassical version of Sobolev spaces, on which is more natural to consider the action of above operators.

\begin{defn} 
\begin{enumerate}[(i)] \label{def of Sobolev spaces}
\item Let $\rho \in \mathbb{N}$.
We define the semiclassical Sobolev space $W^{\rho, \infty}_h(\mathbb{R})$ as the space of  families $(v_h)_{h \in ]0,1]}$ of tempered distributions, such that $\langle hD \rangle^\rho v_h := Op_h(\langle\xi\rangle^\rho)v_h$ is a bounded family of $L^{\infty}$, i.e.
\begin{equation}
W_h^{\rho, \infty}(\mathbb{R}): = \left\{ v_h \in \mathcal{S}'(\mathbb{R})\, \Big| \, \sup_{h\in ]0,1]}\|\langle hD\rangle^{\rho} v_h \|_{L^{\infty}(\mathbb{R})}<+ \infty\right\}\, .
\end{equation}
\item \label{def of H^s_h} Let $s\in\mathbb{R}$. We define the semiclassical Sobolev space $H^s_h(\mathbb{R})$ as the space of families $(v_h)_{h \in ]0,1]}$ of tempered distributions such that $\langle hD \rangle^s v_h := Op_h(\langle\xi\rangle^s)v_h$ is  a bounded family of $L^2$, i.e.
\begin{equation}
H^s_h(\mathbb{R}) := \left\{v_h \in \mathcal{S}'(\mathbb{R})\, \Big| \, \sup_{h\in ]0,1]}\int_{\mathbb{R}} (1+|h\xi|^2)^s |\hat{v}_h(\xi)|^2 \, d\xi < + \infty  \right\} \, .
\end{equation}
\end{enumerate}
\end{defn}

\vspace{0.5cm}
\noindent For future references, we write down the semiclassical Sobolev injection,
\begin{equation} \label{semiclassical Sobolev injection}
\|v_h\|_{W_h^{\rho, \infty}} \le C_{\theta}h^{-\frac{1}{2}} \|v_h\|_{H^{\rho + \frac{1}{2}+\theta}_h}\, , \qquad \forall \theta > 0\, .
\end{equation}

\vspace{0.5cm}
\noindent The following two propositions are stated without proof. They concern the adjoint and the composition of pseudo-differential operators we are considering, and a full detailed treatment is provided in chapter 7 of \cite{dimassi:spectral}, or in chapter 4 of \cite{zworski:semiclassical}.
\vspace{0.5cm}
\begin{prop}[Self-Adjointness]
If $a$ is a real symbol, its Weyl quantization is self-adjoint, 
\begin{equation}
\big(Op_h^w(a)\big)^*=Op_h^w(a)\, .
\end{equation}
\end{prop}

\begin{prop}[Composition for Weyl quantization] \label{Composition for Weyl quantization}
Let $a, b \in \mathcal{S}(\mathbb{R})$. Then
\begin{equation}
Op_h^w(a)\circ Op_h^w(b) = Op_h^w (a\sharp b) \, ,
\end{equation}
where
\begin{equation} \label{a sharp b integral formula}
a \sharp b \,(x,\xi) := \frac{1}{(\pi h)^2}\int_{\mathbb{R}}\int_{\mathbb{R}}\int_{\mathbb{R}}\int_{\mathbb{R}} e^{\frac{2i}{h}\sigma (y, \eta; \, z, \zeta)} a(x+z, \xi + \zeta) b(x+y, \xi +\eta) \; dy d\eta dz d\zeta 
\end{equation}
and
\begin{equation*}
\sigma (y, \eta; \, z, \zeta) =  \eta z - y \zeta\, .
\end{equation*}
\end{prop}

\vspace{0.5cm}
\noindent It is often useful to derive an asymptotic expansion for $a \sharp b$, which allows easier computations than the integral formula \eqref{a sharp b integral formula}.
This expansion is usually obtained by applying the stationary phase argument when $a, b \in S_{\delta, \beta}(M)$, $\delta \in [0,\frac{1}{2}[$ (as shown in \cite{zworski:semiclassical}).
Here we provide an expansion at any order even when one of two symbols belongs to $S_{\frac{1}{2},\beta_1}(M)$ (still having the other in $S_{\delta,\beta_2}(M)$ for $\delta<\frac{1}{2}$, and $\beta_1,\beta_2$ either equal or, if not, one of them equal to zero), whose proof is based on the Taylor development of symbols $a, b$, and can be found in detail in the appendix.

\begin{prop} \label{a sharp b}
Let $a\in S_{\delta_1, \beta_1}(M_1)$, $b\in S_{\delta_2, \beta_2}(M_2)$, $\delta_1, \delta_2 \in [0,\frac{1}{2}]$, $\delta_1 + \delta_2< 1$, $\beta_1, \beta_2 \ge 0$ such that
\begin{equation} \label{beta in symbolic calculus}
\beta_1=\beta_2\ge 0 \qquad \mbox{or} \qquad \big[\beta_1\ne\beta_2  \,\mbox{and} \, \beta_i=0\,,\beta_j>0\,, i\ne j\in\{1,2\} \big]\,.
\end{equation}
Then $a \sharp b \in S_{\delta, \beta}(M_1 M_2)$, where $\delta = \max\{\delta_1, \delta_2 \}$, $\beta=\max\{\beta_1,\beta_2\}$.
Moreover,
\begin{equation} \label{a sharp b asymptotic formula}
a \sharp b = ab + \frac{h}{2i}\{a, b\} 
 + \sum_{\substack{\alpha= (\alpha_1, \alpha_2)\\
                  2\le|\alpha|\le k}}\left(\frac{h}{2i}\right)^{|\alpha|} \frac{(-1)^{\alpha_1}}{\alpha!}\partial^{\alpha_1}_x\partial_{\xi}^{\alpha_2}a \, \partial_x^{\alpha_2}\partial_{\xi}^{\alpha_1}b + r_k \,,
\end{equation}
where $\{a, b\} = \partial_{\xi}a \partial_x b - \partial_{\xi}b \partial_x a$, $r_k \in h^{(k+1)(1-(\delta_1 + \delta_2))}S_{\delta, \beta}(M_1 M_2)$ and
\begin{equation} \label{r_k}
\begin{split}
r_k(x,\xi) = \, \left(\frac{h}{2i}\right)^{k+1} \frac{k+1}{(\pi h)^2}  \sum_{\substack{\alpha= (\alpha_1, \alpha_2)\\
|\alpha|= k+1}} \frac{(-1)^{\alpha_1}}{\alpha!}\int_{\mathbb{R}^4}e^{\frac{2i}{h}(\eta z - y\zeta)} \Big\{& \int_0^1  \partial_x^{\alpha_1}\partial_{\xi}^{\alpha_2}a(x+tz, \xi +t\zeta)(1-t)^k dt  \\
&  \times \partial_y^{\alpha_2}\partial_{\eta}^{\alpha_1}b(x+y, \xi + \eta)\Big\} \, dy d\eta dz d\zeta \,.
\end{split}
\end{equation}
More generally, if $h^{(k+1)\delta_1}\partial^{\alpha}a \in S_{\delta_1,\beta_1}(M^{k+1}_1)$, $h^{(k+1)\delta_2}\partial^{\alpha}b \in S_{\delta_2,\beta_2}(M^{k+1}_2)$, for $|\alpha|=k+1$, for order functions $M_1^{k+1}, M_2^{k+1}$, then $r_k \in h^{(k+1)(1-(\delta_1+\delta_2))}S_{\delta,\beta}(M^{k+1}_1M^{k+1}_2)$.
\end{prop}
\begin{remark}
Observe that
\begin{equation}
a\sharp b = ab + \frac{h}{2i}\{a,b\} + \left(\frac{h}{2i} \right)^2 \left[\frac{1}{2}\partial^2_x a \partial^2_{\xi}b + \frac{1}{2}\partial^2_{\xi}a \partial^2_x b - \partial_x \partial_{\xi}a \partial_x \partial_{\xi}b\right] +r^{a\sharp b}_2 \,,
\end{equation}
so the contribution of order two (and all other contributions of even order) disappears when we calculate the symbol associated to a commutator
\begin{equation} \label{symbol of commutator }
a\sharp b - b\sharp a = \frac{h}{i}\{a,b\} + r_2 \,,
\end{equation}
where $r_2 = r_2^{a\sharp b}- r_2^{b\sharp a}\in h^{3(1-(\delta_1+\delta_2))}S_{\delta,\beta}(M_1 M_2)$.
\end{remark}

\vspace{0.5cm}
\noindent The result of proposition \ref{a sharp b} is still true also when one of order functions, or both, has the form $\langle\frac{x+f(\xi)}{\sqrt{h}}\rangle^{-1}$, for a smooth function $f(\xi)$, $f'(\xi)$ bounded, as stated below and proved as well in the appendix.

\begin{lem} \label{lem of composition with Gamma Tilde}
Let $f(\xi)$ be a smooth function, $f'(\xi)$ bounded.
Consider $a \in S_{\delta_1,\beta_1}(\langle \frac{x+f(\xi)}{\sqrt{h}}\rangle^{-d})$, $d\in\mathbb{N}$, and $b\in S_{\delta_2,\beta_2}(M)$, for $M$ order function or $M(x,\xi)=\langle\frac{x+f(\xi)}{\sqrt{h}}\rangle^{-l}$, $l \in \mathbb{N}$, some $\delta_1 \in [0,\frac{1}{2}]$, $\delta_2 \in [0,\frac{1}{2}[$, $\beta_1, \beta_2\ge 0$ as in \eqref{beta in symbolic calculus}.
Then $a\sharp b \in S_{\delta,\beta}(\langle \frac{x+f(\xi)}{\sqrt{h}}\rangle^{-d}M)$, where $\delta = \max\{\delta_1,\delta_2\}$, $\beta = \max\{\beta_1,\beta_2\}$,  and the asymptotic expansion \eqref{a sharp b asymptotic formula} holds, with $r_k$ given by \eqref{r_k}, $r_k \in h^{(k+1)(1-(\delta_1+\delta_2))}S_{\delta,\beta}(\langle \frac{x+f(\xi)}{\sqrt{h}}\rangle^{-d}M)$. \\
More generally, if $h^{(k+1)\delta_1}\partial^{\alpha}a \in S_{\delta_1,\beta_1}(\langle \frac{x+f(\xi)}{\sqrt{h}}\rangle^{-d'})$ and $h^{(k+1)\delta_2}\partial^{\alpha}b \in S_{\delta_2,\beta_2}(M^{k+1})$, $|\alpha|=k+1$, $M^{k+1}$ order function or $M^{k+1}(x,\xi)=\langle \frac{x+f(\xi)}{\sqrt{h}}\rangle^{-l'}$, for others $d',l' \in \mathbb{N}$, then $r_k \in h^{(k+1)(1-(\delta_1+\delta_2))}S_{\delta,\beta}(\langle \frac{x+f(\xi)}{\sqrt{h}}\rangle^{-d'}M^{k+1})$.
\end{lem}

\subsection{Some Technical Estimates}

\noindent This subsection is mostly devoted to the introduction of some technical results about symbols and operators we will often use in the entire paper, first of all continuity on Sobolev spaces.
We also introduce multi-linear quantizations which will be used in the next section (and which are fully described in \cite{delort:semiclassical}), especially because they make our notations easier and clearer at first.
Moreover, from now on we follow the notation $p(\xi):=\sqrt{1+\xi^2}$.  

\vspace{0.5cm}
\noindent The first statement is about continuity on spaces $H^s_h(\mathbb{R})$, and generalises theorem 7.11 in \cite{dimassi:spectral}.
The second statement concerns instead a result of continuity from $L^2$ to $W_h^{\rho, \infty}$.
In the spirit of \cite{ifrimtataru:global_bounds} for the Schr\"{o}dinger equation, it allows to pass from uniform norms to the $L^2$ norm losing only a power $h^{-\frac{1}{4}-\sigma}$ for a small $\sigma>0$, and not a $h^{-\frac{1}{2}}$ as for the Sobolev injection. 

\begin{prop}[Continuity on $H^s_h$]\label{Continuity from $L^2$ to $L^2$}
Let $s\in \mathbb{R}$.
Let $a \in S_{\delta,\beta}(\langle\xi\rangle^r)$, $r\in\mathbb{R}$, $\delta \in [0, \frac{1}{2}]$, $\beta \ge 0$.
Then $Op_h^w(a)$ is uniformly bounded : $H^s_h(\mathbb{R})\rightarrow H^{s-r}_h(\mathbb{R})$, and there exists a positive constant $C$ independent of $h$ such that 
\begin{equation}
\|Op_h^w(a)\|_{\mathcal{L}(H^s_h;H^{s-r}_h)}\le C\, ,  \qquad \forall h\in ]0,1]\, .
\end{equation}
\end{prop}

\begin{prop} [Continuity from $L^2$ to $W_h^{\rho ,\infty}$]\label{Continuity from $L^2$ to L^inf}
Let $\rho \in \mathbb{N}$.
Let $a \in S_{\delta,\beta}(\langle \frac{x+p'(\xi)}{\sqrt{h}}\rangle^{-1})$, $\delta \in [0, \frac{1}{2}]$, $\beta>0$.
Then $Op_h^w(a)$ is bounded : $L^2(\mathbb{R})\rightarrow W_h^{\rho, \infty}(\mathbb{R})$, and there exists a positive constant $C$ independent of $h$ such that 
\begin{equation}
\|Op_h^w(a)\|_{\mathcal{L}(L^2;W_h^{\rho, \infty})}\le C h^{-\frac{1}{4}-\sigma}\, ,  \qquad \forall h\in ]0,1]\, ,
\end{equation}
where $\sigma>0$ depends linearly on $\beta$. 
\proof
Firstly, remark that thanks to symbolic calculus of lemma \ref{lem of composition with Gamma Tilde}, to estimate the $W_h^{k,\infty}$ norm of an operator whose symbol is rapidly decaying in $|h^{\beta}\xi|$ corresponds actually to estimate the $L^{\infty}$ norm of an operator associated to another symbol (namely, $\tilde{a}(x,\xi)= \langle \xi \rangle^k a(x,\xi)$) which is still in the same class as $a$, up to a small loss on $h$, of order $h^{-k\beta}$.

\vspace{0.5cm}
\noindent From the definition of $Op_h^w(a)v$, and using thereafter integration by part, Cauchy-Schwarz inequality, and Young's inequality for convolutions, we derive what follows :
\begin{equation}
\begin{split}
& |Op_h^w(a)v| = \\
& = \left| \frac{1}{2\pi}\int_{\mathbb{R}}\int_{\mathbb{R}}e^{i(\frac{x}{\sqrt{h}}-y)\xi}a(\frac{x+\sqrt{h} y}{2},\sqrt{h}\xi) v(\sqrt{h}y ) \, dyd\xi \right|\\
& = \left|\frac{1}{(2\pi)^2\sqrt{h}}\int_{\mathbb{R}}\hat{v}(\frac{\eta}{\sqrt{h}})d\eta \int_{\mathbb{R}}\int_{\mathbb{R}} e^{i(\frac{x}{\sqrt{h}}-y)\xi + i\eta y}a(\frac{x+\sqrt{h} y}{2},\sqrt{h}\xi)\, dyd\xi \right| \\
& \le \left|\frac{1}{(2\pi)^2\sqrt{h}}\int_{\mathbb{R}}\hat{v}(\frac{\eta}{\sqrt{h}}) \int_{\mathbb{R}}\int_{\mathbb{R}}\left(\frac{1-i(\frac{x}{\sqrt{h}}-y)\partial_{\xi}}{1+(\frac{x}{\sqrt{h}}-y)^2}\right)^2\left(\frac{1+i(\xi -\eta)\partial_y}{1+(\xi -\eta)^2}\right)^2 \left[e^{i(\frac{x}{\sqrt{h}}-y)\xi + i\eta y}\right] \right. \\
& \left. \hspace{0.5 cm} \times \, a(\frac{x+\sqrt{h} y}{2},\sqrt{h}\xi)\, dyd\xi d\eta  \right| \\
& \le \frac{C}{\sqrt{h}} \int_{\mathbb{R}}\left|\hat{v}(\frac{\eta}{\sqrt{h}})\right| \int_{\mathbb{R}}\int_{\mathbb{R}}\langle \frac{x}{\sqrt{h}}- y\rangle^{-2} \langle \xi - \eta \rangle^{-2}\langle h^{\beta}\sqrt{h}\xi \rangle^{-N} \Big\langle\frac{\frac{x+\sqrt{h} y}{2}+p'(\sqrt{h}\xi)}{\sqrt{h}}\Big\rangle^{-1} dy d\xi d\eta\\
& \le \frac{C}{\sqrt{h}}\left\|\hat{v}(\frac{\eta}{\sqrt{h}})\right\|_{L^2_{\eta}} \|\langle \eta \rangle^{-2}\|_{L^1_{\eta}}\, \left\| \int_{\mathbb{R}} \langle \frac{x}{\sqrt{h}}-y\rangle^{-2}\langle h^{\beta}\sqrt{h}\xi \rangle^{-N} \Big\langle \frac{\frac{x+ \sqrt{h}y}{2}+p'(\sqrt{h}\xi)}{\sqrt{h}}\Big\rangle^{-1} dy \right\|_{L^2_\xi} \\
& \le  Ch^{-\frac{1}{4}}\|v\|_{L^2} \int_{\mathbb{R}} \langle \frac{x}{\sqrt{h}}-y\rangle^{-2} \Big\|\langle h^{\beta}\sqrt{h}\xi \rangle^{-N}\Big\langle \frac{\frac{x+ \sqrt{h}y}{2}+p'(\sqrt{h}\xi)}{\sqrt{h}}\Big\rangle^{-1} \Big\|_{L^2_\xi} dy \, ,
\end{split}
\end{equation}
where $N>0$ is properly chosen later.
We draw attention to the fact that, when we integrated by parts, we used that $a$ belongs to $S_{\delta,\beta}(1)$ with a $\delta \le \frac{1}{2}$, for the loss of $h^{-\delta}$ is offset by the factor $\sqrt{h}$ coming from the derivation of $a(\frac{x+\sqrt{h}y}{2},\sqrt{h}\xi)$ with respect to $y$ and $\xi$.

\noindent To estimate $\|\langle h^{\beta}\sqrt{h}\xi \rangle^{-N}\big\langle\frac{\frac{x+\sqrt{h} y}{2}+p'(\sqrt{h}\xi)}{\sqrt{h}}\big\rangle^{-1}\|_{L^2_\xi}$ we consider a Littlewood-Paley decomposition, i.e.
\begin{equation}
1 = \sum_{k=0}^{+\infty} \varphi_k(\xi)\, , 
\end{equation}
where $\varphi_k(\xi)\in C^{\infty}_0(\mathbb{R})$, $supp\, \varphi_0 \subset B(0,1)$, $\varphi_k(\xi)=\varphi(2^{-k}\xi)$ and $supp\, \varphi \subset \{A^{-1}\le|\xi|\le A \}$, for a constant $A>0$.
Then,
\begin{equation} \label{summation over k}
\begin{split}
\left\|\langle h^{\beta}\sqrt{h}\xi \rangle^{-N}\big\langle\frac{\frac{x+\sqrt{h} y}{2}+p'(\sqrt{h}\xi)}{\sqrt{h}}\big\rangle^{-1}\right\|_{L^2_\xi}^2 & = \frac{1}{\sqrt{h}}\sum_{k\ge 0} \int_{\mathbb{R}} \langle h^{\beta}\xi \rangle^{-2N}\Big\langle\frac{\frac{x+\sqrt{h} y}{2}+p'(\xi)}{\sqrt{h}}\Big\rangle^{-2} \varphi_k(\xi) d\xi \\
&= \frac{1}{\sqrt{h}}\sum_{k\ge 0} I_k \, , 
\end{split}
\end{equation}
where
\begin{equation}
I_0= \int_{\mathbb{R}}\langle h^{\beta}\xi \rangle^{-2N} \Big\langle\frac{\frac{x+\sqrt{h} y}{2}+p'(\xi)}{\sqrt{h}}\Big\rangle^{-2} \varphi_0(\xi) d\xi \, ,
\end{equation}
and
\begin{equation} \label{I_k}
\begin{split}
I_k &= \int_{\mathbb{R}}\langle h^{\beta}\xi \rangle^{-2N} \Big\langle\frac{\frac{x+\sqrt{h} y}{2}+p'(\xi)}{\sqrt{h}}\Big\rangle^{-2} \varphi(2^{-k}\xi) d\xi \\
& = 2^k \int_{\mathbb{R}} \langle h^{\beta}2^k \xi \rangle^{-2N}\Big\langle\frac{\frac{x+\sqrt{h} y}{2}+p'(2^k\xi)}{\sqrt{h}}\Big\rangle^{-2} \varphi(\xi) d\xi \, ,  \\
& \le A^{2N} 2^{(-2N+1)k}h^{-2\beta N} \int_{\mathbb{R}}\Big\langle\frac{\frac{x+\sqrt{h} y}{2}+p'(2^k\xi)}{\sqrt{h}}\Big\rangle^{-2} \varphi(\xi) d\xi \, .
\end{split} \qquad k\ge 1
\end{equation}
For $k\le k_0$, for a fixed $k_0$, $p''(2^k\xi)\ne 0$ on the support of $\varphi$.
As $\xi\rightarrow \pm\infty$ we have the expansion
\begin{equation}
p'(\xi)= \frac{\xi}{\sqrt{1+\xi^2}} = \pm 1 \mp \frac{1}{2\xi^2}+O(|\xi|^{-4}) \,,
\end{equation}
and then
\begin{equation}
p'(2^k\xi)= \pm 1 \mp \frac{2^{-2k}}{2\xi^2} + O(|2^k\xi|^{-4})\, .
\end{equation}
For $k\ge k_0$, the function $\xi \rightarrow g_k(\xi)= 2^{2k}(\frac{x+\sqrt{h} y}{2})+ 2^{2k}p'(2^k\xi)$ is such that $|g_k'(\xi)|= |\xi|^{-3}\sim 1$ on the support of $\varphi$, so for every $k$ we can perform a change of variables $z=g_k(\xi)$ in the last line of \eqref{I_k}.
Hence,
\begin{equation}
\begin{split}
I_k &\le A^{2N} 2^{(-2N+1)k} h^{-2\beta N}\int \langle\frac{z}{2^{2k}\sqrt{h}}\rangle^{-2} \varphi(g_k^{-1}(z)) dz \\
& \le A^{2N} 2^{(-2N+3)k}h^{-2\beta N}\sqrt{h}\int \langle z \rangle^{-2}  dz \\
& \le C 2^{(-2N+3)k}h^{-2\beta N}\sqrt{h} \, ,
\end{split}
\end{equation}
so taking the summation of all $I_k$ for $k\ge 0$ we deduce
\begin{equation}
\left\|\langle h^{\beta}\sqrt{h}\xi\rangle^{-N}\Big\langle\frac{\frac{x+\sqrt{h} y}{2}+p'(\sqrt{h}\xi)}{\sqrt{h}}\Big\rangle^{-1}\right\|_{L^2_\xi}\le C h^{-\beta N}\sum_{k\ge 0}2^{(\frac{-2N+3}{2})k} \le  C' h^{-\beta N} \, ,  
\end{equation}
if we choose $N>0$ such that $\frac{-2N+3}{2}<0$ (e.g. $N=2$).
Finally
\begin{equation}
\|Op_h^w(a)\|_{\mathcal{L}(L^2; W_h^{\rho, \infty})} = O( h^{-\frac{1}{4}- \sigma })\,,
\end{equation}
where $\sigma(\beta) = (N + \rho)\beta$ depends linearly on $\beta$.
\endproof
\end{prop}

\vspace{0.5 cm}
\noindent The following lemma shows that we have nice upper bounds for operators acting on $v$ whose symbols are supported for $|\xi|\ge h^{-\beta}$, $\beta >0$, provided that we have an a priori $H^s_h$ estimate on $v$, with large enough $s$.

\begin{lem} \label{new estimate 1-chi}
Let $s'\ge 0$.
Let $\chi \in C^{\infty}_0(\mathbb{R})$, $\chi \equiv 1$ in a neighbourhood of zero, e.g.
\begin{equation}
\begin{split}
& \chi(\xi)= 1 \, , \qquad \emph{for}\, |\xi|<C_1 \\
& \chi(\xi) = 0 \, , \qquad \emph{for}\, |\xi|> C_2 \,.
\end{split}
\end{equation}
Then
\begin{equation}
\|Op_h((1-\chi)(h^{\beta}\xi))v \|_{H^{s'}_h}\le C h^{\beta(s-s')}\|v\|_{H^s_h} \, , \qquad \qquad \forall s> s'\, .
\end{equation}
\proof
The result is a simple consequence of the fact that $(1-\chi)(h^{\beta}\xi)$ is supported for $|\xi|\ge C_1 h^{-\beta}$, because
\begin{equation}
\begin{split}
\|Op_h((1-\chi)(h^{\beta}\xi))v \|_{H^{s'}_h}^2 & = \int (1+ |h\xi|^2)^{s'}|(1-\chi)(h^{\beta}h\xi)|^2 |\hat{v}(\xi)|^2 d\xi \\
& = \int (1+ |h\xi|^2)^s (1+ |h\xi|^2)^{s'-s}|(1-\chi)(h^{\beta}h\xi)|^2 |\hat{v}(\xi)|^2 d\xi \\
& \le C h^{2\beta(s-s')}\|v\|^2_{H^s_h}\, ,
\end{split}
\end{equation}
where the last inequality follows from an integration on $|h\xi|>C_1 h^{-\beta}$, and from the two following conditions $s'-s< 0$, $(1+|h\xi|^2)^{s'-s}\le C h^{-2\beta(s'-s)}$.
\endproof
\end{lem}
\noindent This result is useful when we want to reduce essentially to symbols rapidly decaying in $|h^{\beta}\xi|$, for example in the intention of using proposition \ref{Continuity from $L^2$ to L^inf} or when we want to pass from a symbol of a certain positive order to another one of order zero, up to small losses of order $O(h^{-\sigma})$, $\sigma>0$ depending linearly on $\beta$. We can always split a symbol using that $1= \chi(h^{\beta}\xi) + (1-\chi)(h^{\beta}\xi)$, and consider as remainders all contributions coming from the latter.

\vspace{0.5cm}
\noindent Define the set $\Lambda := \{(x,\xi)\in \mathbb{R}\times\mathbb{R}\, |\, x+ p'(\xi)= 0 \}$, i.e. the graph of the function $x\in ]-1,1[\rightarrow d\varphi(x)$, $\varphi(x) = \sqrt{1-x^2}$, as drawn in picture \ref{fig: Lambda portrait}.
We will use the following technical lemma, whose proof can be found in lemma 1.2.6 in \cite{delort:semiclassical} :

\begin{lem} \label{lem on e and etilde}
Let $\gamma \in C^{\infty}_0(\mathbb{R})$.
If the support of $\gamma$ is sufficiently small, the two functions defined below
\begin{equation} \label{def of e and etilde}
e_{\pm}(x,\xi) = \frac{x + p'(\pm\xi)}{\xi \mp d\varphi(x)} \gamma\left(\langle\xi\rangle^2 (x+p'(\pm\xi))\right) \quad \emph{and} \quad \tilde{e}_{\pm}(x,\xi) = \frac{\xi \mp d\varphi(x)}{x + p'(\pm\xi)} \gamma\left(\langle\xi\rangle^2(x+p'(\pm\xi))\right)
\end{equation}
verify estimates
\begin{equation}
\begin{split}
|\partial^{\alpha}_x \partial_{\xi}^{\beta}e_{\pm}(x,\xi)|& \le C_{\alpha\beta} \langle\xi\rangle^{-3+2\alpha-\beta}\,,\\
|\partial^{\alpha}_x \partial_{\xi}^{\beta}\tilde{e}_{\pm}(x,\xi)|& \le C_{\alpha\beta} \langle\xi\rangle^{3+2\alpha-\beta} \,.
\end{split}
\end{equation}
Moreover, if $supp\gamma$ is small enough, then on the support of $\gamma(\langle\xi\rangle^2(x+p'(\pm\xi)))$ one has $\langle d\varphi\rangle\sim \langle\xi\rangle$ and there is a constant $A>0$ such that, on that support
\begin{equation} \label{bound for x pm 1}
\begin{split}
A^{-1}\langle\xi\rangle^{-2} & \le \pm x +1 \le A \langle\xi\rangle^{-2} \, , \qquad \xi \rightarrow +\infty \\
A^{-1}\langle\xi\rangle^{-2} & \le \mp x +1 \le A \langle\xi\rangle^{-2} \, , \qquad \xi \rightarrow -\infty 
\end{split}
\end{equation}
Finally, for every $k\in \mathbb{N}$
\begin{equation}
\partial^k(d\varphi(x))=O(\langle d\varphi\rangle^{1+2k}) \,.
\end{equation}
\end{lem}

\vspace{0.5cm}
\begin{lem} \label{Class of Gamma}
Let $\gamma\in C^{\infty}_0(\mathbb{R})$ such that $\gamma \equiv 1$ in a neighbourhood of zero, and define $\Gamma(x,\xi)= \gamma(\frac{x+p'(\xi)}{\sqrt{h}})$.
Then $\Gamma \in S_{\frac{1}{2}, 0}(\langle\frac{x+p'(\xi)}{\sqrt{h}} \rangle^{-N})$, for all $N\ge 0$.
\proof
Let $N\in \mathbb{N}$.
Since $\gamma \in C^{\infty}_0(\mathbb{R})$, $p'' \in S_{0,0}(1)$, we have
\begin{equation}
\begin{split}
| \Gamma(x,\xi)| & \le \|\langle x \rangle ^N \gamma(x) \|_{L^{\infty}} \langle\frac{x+p'(\xi)}{\sqrt{h}}\rangle^{-N} \, , \\
| \partial_x \Gamma (x,\xi) | & = \big|\gamma'(\frac{x+p'(\xi)}{\sqrt{h}})\frac{1}{\sqrt{h}}\big| \le h^{-\frac{1}{2}} \| \langle x\rangle^N \gamma'(x) \|_{L^{\infty}} \langle\frac{x+p'(\xi)}{\sqrt{h}}\rangle^{-N} \, , \\
| \partial_{\xi} \Gamma (x,\xi) | & = \big|\gamma'(\frac{x+p'(\xi)}{\sqrt{h}})\frac{p''(\xi)}{\sqrt{h}}\big| \lesssim h^{-\frac{1}{2}} \|\langle x \rangle ^N \gamma'(x) \|_{L^{\infty}} \langle\frac{x+p'(\xi)}{\sqrt{h}}\rangle^{-N} \, ,
\end{split}
\end{equation}
and going on one can prove that $|\partial_x^{\alpha_1}\partial_{\xi}^{\alpha_2}\Gamma| \le C_{\alpha_1, \alpha_2, N} h^{-\frac{1}{2}(\alpha_1+\alpha_2)} \langle\frac{x+p'(\xi)}{\sqrt{h}}\rangle^{-N}$.
\endproof
\end{lem}

\paragraph{Multi-linear Operators.}
We briefly generalise some definitions given at the beginning of this section in order to introduce multi-linear operators.

\vspace{0.5cm}
\noindent Let $n \in \mathbb{N}^*$ and set $\xi = (\xi_1, \dotsc, \xi_n)$.
An order function on $\mathbb{R}\times \mathbb{R}^n$ will be a smooth function $(x,\xi)\rightarrow M(x,\xi)$ satisfying \eqref{def ineq order function}, where $\langle \xi-\eta\rangle^{N_0}$ is replaced by
 $$ \prod_{i=1}^n \langle \xi_i - \eta_i\rangle^{N_0}\, .$$   Equivalently, we define the class $S_{\delta,\beta}(M,n)$, for some $\delta\ge 0$, $\beta\ge 0$ and $M(x,\xi)$ order function on $\mathbb{R}\times \mathbb{R}^n$, to be the set of smooth functions 
\begin{align*}
(x,\xi_1,\dotsc, \xi_n, h)  & \rightarrow a(x,\xi, h) \\
\mathbb{R}\times\mathbb{R}^n\times ]0,1] & \rightarrow \mathbb{C}
\end{align*}
satisfying the inequality \eqref{symbol in S delta beta M}, $\forall \alpha_1 \in \mathbb{N}, \alpha_2 \in \mathbb{N}^n$, $\forall k, N \in \mathbb{N}$.

\begin{defn}
Let $a$ be a symbol in $S_{\delta,\beta}(M, n)$ for some order function $M$, some $\delta\ge 0$, $\beta \ge 0$.
\begin{enumerate}[(i)]
\item We define the $n$-linear operator $Op(a)$ acting on test functions $v_1, \dotsc, v_n$ by
\begin{equation}
Op(a)(v_1, \dotsc, v_n) = \frac{1}{(2\pi)^n}\int_{\mathbb{R}^n}e^{ix(\xi_1 + \dots + \xi_n)} a(x, \xi_1, \dotsc, \xi_n)\prod_{l=1}^n \hat{v}_j(\xi_l)\, d\xi_1 \dots d\xi_n \, .
\end{equation}
\item We also define the $n$-linear semiclassical operator $Op_h(a)$ acting on test functions $v_1, \dotsc, v_n$ by
\begin{equation}
Op_h(a)(v_1, \dotsc, v_n) = \frac{1}{(2\pi h)^n}\int_{\mathbb{R}^n}e^{\frac{i}{h}x(\xi_1 + \dots + \xi_n)} a(x, \xi_1, \dotsc, \xi_n)\prod_{l=1}^n \hat{v}_j(\xi_l)\, d\xi_1 \dots d\xi_n \, .
\end{equation}
\end{enumerate}
\end{defn}

\noindent For a further need of compactness in our notations, we introduce $I=(i_1,\dots,i_n)$ a $n$-dimensional vector, $i_k\in \{ 1,-1\}$ for every $k=1,\dots,n$. We set $|I|=n$ and define
\begin{equation}
w_I=(w_{i_1},\dots,w_{i_n})\, , \qquad w_1 = w \, , \, w_{-1}= \bar{w}\, ,
\end{equation}
while $m_I(\xi) \in S_{\delta,\beta}(M, n)$ will be always in what follows a symbol of the form
\begin{equation} \label{multilinear symbol}
m_I(\xi)= m_1^I(\xi_1)\cdots m_n^I(\xi_n) \,.
\end{equation}
We warn the reader that in following sections, when we focus on a fixed general symbol $m_I(\xi)$, we will often refer to components $m_k^I(\xi_k)$ as $m_k(\xi_k)$, forgetting the superscript $I$ in order to make notations lighter. 
Sometimes we will also write $m_k(\xi)$ if this makes no confusion.

\section{Semiclassical Reduction to an ODE.} \label{Semiclassical Reduction to an ODE}

In this section we want to reformulate the Cauchy problem \eqref{KG} and to deduce a new equation which can be transformed into an ODE.
It is organised in three subsections.
In the first one, we introduce semiclassical coordinates, rewrite the problem in this new framework and state the main theorem.
The second and third sections are devoted to the proof of the main theorem.
In particular, in the second one we introduce some technical lemmas we often refer to and we estimate $v$ when it is away from $\Lambda$. 
In the third one, we first cut symbols in the cubic nonlinearity near $\Lambda$ and away from points $x=\pm 1$, and develop them at $\xi = d\varphi(x)$, transforming multi-linear pseudo-differential operators in smooth functions of $x$; then, we repeat the development argument for $Op_h^w(x\xi + p(\xi))$.

\subsection{Semiclassical Coordinates and Statement of the Main Result} 
Let $u$ be a solution of \eqref{KG} and set
\begin{equation} \label{def of w}
\begin{cases}
& w = (D_t + \sqrt{1+D_x^2}) u \\
& \bar{w} = -(D_t - \sqrt{1+D_x^2})u
\end{cases} \, , \qquad \qquad
\begin{cases}
&u  = \langle D_x \rangle^{-1}(\frac{w + \bar{w}}{2}) \\
& D_t u = \frac{w-\bar{w}}{2}
\end{cases}\,.
\end{equation}
With notations introduced in \eqref{polynomial decomposition}, the function $w$ satisfies the following equation
\begin{equation} \label{half KG in w with polynomials}
\begin{split}
(D_t - \sqrt{1+D_x^2}) w & = \sum_{k=0}^3 i^k P'_k\left(\langle D_x\rangle^{-1}(\frac{w+\bar{w}}{2}); \frac{w-\bar{w}}{2}, D_x \langle D_x\rangle^{-1}(\frac{w+\bar{w}}{2})\right)\\
& + \sum_{k=0}^2 i^k P''_k \left( \langle D_x\rangle^{-1}(\frac{w+\bar{w}}{2}), D_x (\frac{w-\bar{w}}{2}), D_x^2 \langle D_x\rangle^{-1}(\frac{w+\bar{w}}{2})\right. ; \\
&\left. \hspace{2.2cm} \, \frac{w-\bar{w}}{2}, D_x \langle D_x\rangle^{-1}(\frac{w+\bar{w}}{2}) \right) \, .
\end{split}
\end{equation}
Observe that operators which take the place of second derivatives have symbols of order one, while all other symbols are of order zero or negative ($-1$). 
Let us simplify the notation for the rest of the section by rewriting the nonlinearity in term of multi-linear pseudo-differential operators introduced in the previous section, namely as
\begin{equation}
\sum_{|I|=3} Op(m_I)(w_I) + \sum_{|I|=3}Op(\widetilde{m}_I)(w_I) \,,
\end{equation}
where symbols $m_I$, $\widetilde{m}_I$ are of the form \eqref{multilinear symbol}.
Moreover, $m_I$ will denote symbols of order equal or less than zero, in the sense that all occurring symbols $m_k^I$ are of order equal or less than zero, while in $\widetilde{m}_I$ there will be exactly one symbol of order one, thanks to the quasi-linear nature of the starting equation.
Therefore \eqref{half KG in w with polynomials} is rewritten as
\begin{equation} \label{half KG in w with multi-operators}
(D_t - \sqrt{1+D_x^2}) w =\sum_{|I|=3} Op(m_I)(w_I) + \sum_{|I|=3}Op(\widetilde{m}_I)(w_I) \, ,
\end{equation}
and passing to the semiclassical framework by
\begin{equation} \label{def of v}
w(t,x)= \frac{1}{\sqrt{t}} v (t, \frac{x}{t})\, , \qquad h := \frac{1}{t} \, ,
\end{equation}
we obtain
\begin{equation} \label{half KG in v with multi-operators}
\big(D_t - Op_h^w(x \xi + p(\xi)\big)v = h \sum_{|I|=3} Op_h(m_I)(v_I) + h \sum_{|I|=3}Op_h(\widetilde{m}_I)(v_I) \, , 
\end{equation}
where $p(\xi)=\sqrt{1+\xi^2}$ and where we used the equality $Op_h(x\xi + p(\xi) + \frac{h}{2i})= Op_h^w(x \xi + p(\xi))$ following from
\begin{equation*}
\begin{split}
Op_h^w(x\xi) &= \frac{h}{2}D_x x + \frac{h}{2}x D_x \\
& = \frac{h}{2i} + x\, hD_x = \frac{h}{2i} + Op_h(x\xi)\, .
\end{split}
\end{equation*}
Furthermore, we write explicitly the nonlinearity of the equation, which will be useful hereinafter 
\begin{equation} \label{half KG for v with polynomials}
\begin{split}
\big(D_t - Op_h^w(x \xi + p(\xi)\big)v &= h \sum_{k=0}^3 i^k P'_k\left(\langle hD\rangle^{-1}(\frac{v+\bar{v}}{2}); \frac{v-\bar{v}}{2}, (hD) \langle hD\rangle^{-1}(\frac{v+\bar{v}}{2})\right)\\
& + h \sum_{k=0}^2 i^k P''_k \left( \langle hD\rangle^{-1}(\frac{v+\bar{v}}{2}), (hD) (\frac{v-\bar{v}}{2}), (hD)^2 \langle hD\rangle^{-1}(\frac{v+\bar{v}}{2});\right. \\
& \left. \hspace{2.3 cm} \, \frac{v-\bar{v}}{2}, (hD) \langle hD\rangle^{-1}(\frac{v+\bar{v}}{2}) \right) \, .
\end{split}
\end{equation}
Let us also define the operator $\mathcal{L}$ to be 
\begin{equation}
\mathcal{L}:= \frac{1}{h}Op_h^w(x+p'(\xi))\, .
\end{equation}

\vspace{0.5cm}
\noindent The equation \eqref{half KG in v with multi-operators} represents for us the starting point to deduce an ODE satisfied by $v$, from which it will be easier to derive an estimate on the $L^{\infty}$ norm of $v$. 
In reality, we will need more than an uniform estimate for $v$, namely we have to involve also a certain number of its derivatives, and then to control its $W^{\rho, \infty}_h$ norm for a fixed $\rho > 0$. 
With this in mind, we set $\Gamma(x,\xi) = \gamma (\frac{x + p'(\xi)}{\sqrt{h}})$, for a function $\gamma \in C^{\infty}_0(\mathbb{R})$, $\gamma \equiv 1$ in a neighbourhood of zero, with a small support.
From lemma \ref{Class of Gamma}, $\Gamma \in S_{\frac{1}{2},0}(\langle \frac{x+ p'(\xi)}{\sqrt{h}}\rangle^{-N})$ for every $N \in \mathbb{N}^*$, and case by case we will choose the right power we need. 
We consider also $\Sigma(\xi)= \langle \xi \rangle^{\rho}$ (in practice, at times we consider $\rho-1 \in \mathbb{N}$, with $\rho$ introduced for $u$ in theorem \ref{bootstrap theorem},  when we prove the bootstrap, or $\rho =-1$ when we develop asymptotics), and define
\begin{equation} \label{def of v Sigma}
v^{\Sigma}  := Op_h(\Sigma)v \, , 
\end{equation}
together with 
\begin{equation}
\begin{split}
& v^{\Sigma}_{\Lambda} := Op_h^w(\Gamma) v^{\Sigma} \, , \\
& v^{\Sigma}_{\Lambda^c}  := Op_h^w(1-\Gamma)v^{\Sigma} \, ,
\end{split}
\end{equation}
and symbols
\begin{equation} \label{def of m_I Sigma}
\begin{split}
m_I^{\Sigma}(\xi)& = \prod_{k=1}^3m_k^{I,\Sigma}(\xi_k) : = \prod_{k=1}^3m_k^I(\xi_k)\Sigma(\xi_k)^{-1}\,, \\
\widetilde{m}_I^{\Sigma}(\xi)& = \prod_{k=1}^3\widetilde{m}_k^{I,\Sigma}(\xi_k):=  \prod_{k=1}^3\widetilde{m}_k^I(\xi_k)\Sigma(\xi_k)^{-1}\, .
\end{split}
\end{equation}

\vspace{0.5cm}
\noindent The main result we want to prove in this section is the following:
\begin{thm}[Reformulation of the PDE] \label{thm : From the PDE problem to an ODE}
Suppose that we are given constants $A', B'>0$, some $T>1$ and a solution $v \in L^{\infty}([1,T]; H^s_h)\cap L^{\infty}([1,T]; W^{\rho,\infty}_h)$ of the equation \eqref{half KG in v with multi-operators} (or, equivalently, of \eqref{half KG for v with polynomials}), satisfying the following a priori bounds, for any $\varepsilon \in ]0,1]$, $t \in [1,T]$,
\begin{align} 
\|v(t,\cdot)\|_{W^{\rho,\infty}_h}& \le A'\varepsilon \,, \label{a priori estimates on v} \\
\|\mathcal{L}v(t,\cdot)\|_{L^2}+ \|v(t,\cdot)\|_{H^s_h} &\le B'h^{-\sigma}\varepsilon \label{a priori estimate on Lv and v}\,,
\end{align}
for some $\sigma>0$ small enough.
Then, with preceding notations, $v^{\Sigma}_{\Lambda}$ is solution of 
\begin{equation} \label{ODE in section 4}
\begin{split}
D_t v^{\Sigma}_{\Lambda} & = \varphi(x)\theta_h(x) v^{\Sigma}_{\Lambda} + h \Phi_1^{\Sigma}(x)\theta_h(x) |v^{\Sigma}_{\Lambda}|^2v^{\Sigma}_{\Lambda} \\
& + h Op_h^w(\Gamma)\left[\Phi_3^{\Sigma}(x)\theta_h(x)(v^{\Sigma}_{\Lambda})^3 + \Phi_{-1}^{\Sigma}(x)\theta_h(x) |v^{\Sigma}_{\Lambda}|^2\overline{v^{\Sigma}_{\Lambda}} + \Phi_{-3}^{\Sigma}(x)\theta_h(x) (\overline{v^{\Sigma}_{\Lambda}})^3 \right] + hR(v) \, ,
\end{split}
\end{equation}
with $(\theta_h(x))_h$ a family of smooth functions compactly supported in $]-1,1[$, some smooth coefficients $\Phi_j^{\Sigma}(x)$, $|\Phi_j^{\Sigma}(x)|=O(h^{-\sigma'})$ on the support of $\theta_h$, for $j \in \{3,1,-1,-3\}$ and a small $\sigma'>0$.
Moreover, $R(v)$ is a remainder verifying the following estimates
\begin{align}
\|R(v)\|_{L^2}  & \le C h^{\frac{1}{2}-\sigma}\, (\|\mathcal{L}v\|_{L^2}+ \|v\|_{H^s_h})\, ,  \label{estimates L^2 to L^2 of R(v)} \\ 
\|R(v)\|_{L^{\infty}}& \le C h^{\frac{1}{4}-\sigma} \, (\|\mathcal{L}v\|_{L^2}+ \|v\|_{H^s_h})\, , \label{estimates L^2 to L^inf of R(v)}
\end{align}
for a new small $\sigma\ge 0$.
\end{thm}

\noindent Smooth coefficients $\Phi^{\Sigma}_j(x)$ in \eqref{ODE in section 4} may be explicitly calculated starting from the nonlinearity in \eqref{half KG for v with polynomials}, and in particular this will be done for $\Phi^{\Sigma}_1(x)$ at the beginning of section \ref{Study of the ODE and End of the Proof}.
Afterwards, we will use the notation $R_1(v)$ to refer to a remainder satisfying the following estimates:
\begin{align} \label{L^2-L^2 estimate of R}
\|R_1(v)\|_{H^{\rho}_h} &\le C h^{\frac{1}{2}-\sigma} (\|\mathcal{L}v\|_{L^2}+ \|v\|_{H^s_h})\,, \\
\|R_1(v)\|_{L^ { \infty}} & \le C h^{\frac{1}{4}-\sigma} (\|\mathcal{L}v\|_{L^2}+ \|v\|_{H^s_h})\,,  \label{L^2-L^inf estimate of R}
\end{align}
for a small $\sigma\ge 0$.

\subsection{Technical Results}
We estimate $v^{\Sigma}_{\Lambda^c}$ as follows :
\begin{lem}\label{lemma GammaTilde}
Let $\widetilde{\Gamma}(\xi)$ a smooth function such that $|\partial^{\alpha}\widetilde{\Gamma}|\lesssim \langle \xi \rangle^{-\alpha}$, $\chi$ as in lemma \ref{new estimate 1-chi}, $\beta > 0$. Then
\begin{equation}
Op_h^w(\widetilde{\Gamma}(\frac{x+p'(\xi)}{\sqrt{h}}))v^{\Sigma}= Op_h^w\left(\Sigma(\xi)\chi(h^{\beta}\xi)\widetilde{\Gamma}(\frac{x+p'(\xi)}{\sqrt{h}})\right)v + R_1(v) \,,
\end{equation}
where $R_1(v)$ is a remainder satisfying \eqref{L^2-L^2 estimate of R}, \eqref{L^2-L^inf estimate of R}.
\proof
We consider a function $\chi$ as in lemma \ref{new estimate 1-chi}, and we write
\begin{equation}
\begin{split}
Op_h^w(\widetilde{\Gamma}(\frac{x+p'(\xi)}{\sqrt{h}}))v^{\Sigma} & = Op_h^w(\widetilde{\Gamma}(\frac{x+p'(\xi)}{\sqrt{h}}))Op_h^w(\Sigma(\xi)\chi(h^{\beta}\xi))v  \\
& + Op_h^w(\widetilde{\Gamma}(\frac{x+p'(\xi)}{\sqrt{h}}))Op_h^w(\Sigma(\xi)(1-\chi)(h^{\beta}\xi))v \,,
\end{split}
\end{equation} 
for $\beta>0$.
The second term in the right hand side represents a remainder $R_1(v)$ satisfying the two inequalities of the statement just because $\widetilde{\Gamma}(\frac{x+p'(\xi)}{\sqrt{h}}) \in S_{\frac{1}{2},0}(1)$ (so, for instance, $\|Op_h^w(\widetilde{\Gamma}(\frac{x+p'(\xi)}{\sqrt{h}}))\|_{\mathcal{L}(H^{\rho+1}_h ; W_h^{\rho, \infty})}=O(h^{-\frac{1}{2}})$ by Sobolev inequality \eqref{semiclassical Sobolev injection} and proposition \ref{Continuity from $L^2$ to $L^2$}) and $(1-\chi)(h^{\beta}\xi)$ is supported for $|\xi|\ge h^{-\beta}$, so that we can use essentially lemma \ref{new estimate 1-chi}.

\vspace{0.5cm}
\noindent On the other hand, since $|\partial^{\alpha}\widetilde{\Gamma}|\le \langle \xi \rangle^{-\alpha}$ and $\Sigma(\xi)\chi(h^{\beta}\xi)\in h^{-\sigma}S_{0,\beta}(1)$, with 
\begin{equation} \label{cases for sigma}
\sigma = 
\begin{cases}
\rho \beta  \quad &\emph{if}\,\, \rho \in \mathbb{N} \\
0\quad &\emph{if}\,\, \rho<0
\end{cases}
\end{equation}
we use the composition formula of lemma \ref{lem of composition with Gamma Tilde} for the first term in the right hand side, i.e.
\begin{equation} 
Op_h^w(\widetilde{\Gamma}(\frac{x+p'(\xi)}{\sqrt{h}}))Op_h^w(\Sigma(\xi)\chi(h^{\beta}\xi))v = Op_h^w\left(\Sigma(\xi)\chi(h^{\beta}\xi)\widetilde{\Gamma}(\frac{x+p'(\xi)}{\sqrt{h}})\right) v + Op_h(r_0)v \, ,
\end{equation}
where $r_0 \in h^{\frac{1}{2}-\sigma}S_{\frac{1}{2},\beta}(\langle \frac{x+p'(\xi)}{\sqrt{h}}\rangle^{-1})$.
So $Op_h(r_0)v$ satisfies inequalities \eqref{L^2-L^2 estimate of R}, \eqref{L^2-L^inf estimate of R} respectively by propositions \ref{Continuity from $L^2$ to $L^2$} and \ref{Continuity from $L^2$ to L^inf}.
\endproof
\end{lem}
 
\begin{lem} \label{composition gamma-1 and its argument}
Let $\widetilde{\Gamma}(\xi)$ be a smooth function such that $|\partial^{\alpha}\widetilde{\Gamma}| \lesssim \langle \xi\rangle^{-\alpha}$, $c(x,\xi) \in S_{\delta,\beta}(1)$, $c'(x,\xi) \in S_{\delta',0}(1)$, with $\delta, \delta'\in [0,\frac{1}{2}[$, $\beta>0$. Then
\begin{equation} \label{sharp gammatilde and its argument}
c(x,\xi)\widetilde{\Gamma}(\frac{x+p'(\xi)}{\sqrt{h}}) \sharp (x+p'(\xi)) =c(x,\xi)\widetilde{\Gamma}(\frac{x+p'(\xi)}{\sqrt{h}}) (x+p'(\xi)) +h^{1-2\delta}r \, 
\end{equation} 
with $r \in S_{\frac{1}{2},\beta}(1)$, and
\begin{align} 
\big\|Op_h^w\big(c(x,\xi)\widetilde{\Gamma}(\frac{x+p'(\xi)}{\sqrt{h}})(x+p'(\xi))\big)Op_h^w(c')v \big\|_{L^2} &\le h^{1-\sigma}(\|\mathcal{L}v\|_{L^2} + \|v\|_{H^s_h}) \, , \label{est L2 in lemma 4.3 for GammaTilde}\\
\big\|Op_h^w\big(c(x,\xi)\widetilde{\Gamma}(\frac{x+p'(\xi)}{\sqrt{h}})(x+p'(\xi))\big)Op_h^w(c')v \big\|_{L^{\infty}} & \le h^{\frac{1}{2}-\sigma} (\|\mathcal{L}v\|_{L^2} + \|v\|_{H^s_h}) \label{est Linf in lemma 4.3 for GammaTilde}\, ,
\end{align}
with $\sigma=\sigma(\delta,\delta',\beta)\rightarrow 0$ as $\delta,\delta',\beta\rightarrow 0$. \\
Moreover, if $\widetilde{\Gamma}= \Gamma_{-1}$, with $|\partial^{\alpha}\Gamma_{-1}|\lesssim \langle\xi\rangle^{-1-\alpha}$, then in \eqref{sharp gammatilde and its argument} $r \in S_{\frac{1}{2},\beta}(\langle \frac{x+p'(\xi)}{\sqrt{h}} \rangle^{-1})$, and the $L^{\infty}$ estimate can be improved
\begin{equation} \label{est Gamma-1}
\big\|Op_h^w\big(c(x,\xi)\Gamma_{-1}(\frac{x+p'(\xi)}{\sqrt{h}})(x+p'(\xi))\big)Op_h^w(c')v \big\|_{L^{\infty}}  \le h^{\frac{3}{4}-\sigma} (\|\mathcal{L}v\|_{L^2} + \|v\|_{H^s_h}) \, . 
\end{equation}

\proof
The result is immediate if we use the development of proposition \ref{a sharp b} at order one,
\begin{equation} \label{sharp chi Gamma-1}
\begin{split}
c(x,\xi)\widetilde{\Gamma}(\frac{x+p'(\xi)}{\sqrt{h}}) \sharp (x+p'(\xi)) & = c(x,\xi)\widetilde{\Gamma}(\frac{x+p'(\xi)}{\sqrt{h}}) (x+p'(\xi))\\
& + \frac{h}{2i}\left\{c(x,\xi)\widetilde{\Gamma}(\frac{x+p'(\xi)}{\sqrt{h}}), (x+p'(\xi))\right\} + h r_1\, , 
\end{split}
\end{equation}
where $r_1 \in h^{-2\delta}S_{\frac{1}{2},\beta}(1)$, while by direct calculation the Poisson bracket is equal to:
\begin{equation*}
\left\{c(x,\xi)\widetilde{\Gamma}(\frac{x+p'(\xi)}{\sqrt{h}}), (x+p'(\xi))\right\}  =  \widetilde{\Gamma}(\frac{x+p'(\xi)}{\sqrt{h}})(\partial_{\xi}c - p''\partial_x c )\, ,
\end{equation*}
$\widetilde{\Gamma}(\frac{x+p'(\xi)}{\sqrt{h}})(\partial_{\xi}c - p''\partial_x c ) \in h^{-\delta}S_{\frac{1}{2},\beta}(1)$.
Therefore
\begin{equation} \label{dev in lemma Gammatilde}
\begin{split}
& Op_h^w\big(c(x,\xi)\widetilde{\Gamma}(\frac{x+p'(\xi)}{\sqrt{h}})(x+p'(\xi))\big)Op_h^w(c')v = \\
&=  h Op_h^w\big(c(x,\xi)\widetilde{\Gamma}(\frac{x+p'(\xi)}{\sqrt{h}})\big)\mathcal{L}Op_h^w(c')v \\
& + h Op_h^w\big(\widetilde{\Gamma}(\frac{x+p'(\xi)}{\sqrt{h}})(\partial_{\xi}c-p''\partial_xc ) +r_1 \big)Op_h^w(c')v \, ,
\end{split}
\end{equation}
and the application of proposition \ref{Continuity from $L^2$ to $L^2$}, along with Sobolev injection \eqref{semiclassical Sobolev injection}, immediately implies that the second term in the right hand side satisfies estimates \eqref{est L2 in lemma 4.3 for GammaTilde}, \eqref{est Linf in lemma 4.3 for GammaTilde}. 
Moreover, $[\mathcal{L},Op_h^w(c')]= i(\partial_{\xi}c'- p''\partial_x c') + h^{1-2\delta'}r_1$, $r_1$ being a symbol in $S_{\delta',0}(1)$, hence it belongs to $h^{-\delta'}S_{\delta',0}(1)$, and its quantization is a bounded operator from $L^2$ to $L^2$ by proposition \ref{Continuity from $L^2$ to $L^2$} up to a small loss in $h^{-\delta'}$. 
This remark, together with $c(x,\xi)\widetilde{\Gamma}(\frac{x+p'(\xi)}{\sqrt{h}}) \in S_{\frac{1}{2},\beta}(1)$, $c' \in S_{\delta',0}(1)$, proposition \ref{Continuity from $L^2$ to $L^2$}, and Sobolev injection imply that also the first term in the right hand side verifies estimates in \eqref{est L2 in lemma 4.3 for GammaTilde}, \eqref{est Linf in lemma 4.3 for GammaTilde}.
The same reasoning as above can be applied when $\widetilde{\Gamma}=\Gamma_{-1}$ with $|\partial^{\alpha}\Gamma_{-1}|\lesssim \langle \xi \rangle^{-1-\alpha}$.
In this case, the development in \eqref{sharp chi Gamma-1} is justified by lemma \ref{lem of composition with Gamma Tilde}.
Moreover, symbols involving $c(x,\xi)\Gamma_{-1}(\frac{x+p'(\xi)}{\sqrt{h}})$ stay in $S_{\frac{1}{2},\beta}(\langle \frac{x+p'(\xi)}{\sqrt{h}}\rangle^{-1})$, so we can apply proposition \ref{Continuity from $L^2$ to L^inf}, instead of Sobolev injection, to control the $L^{\infty}$ norm, losing only a power $h^{-\frac{1}{4}-\sigma}$, for a small $\sigma>0$ (and not $h^{-\frac{1}{2}}$ due to Sobolev estimate) and so deriving the improved estimate \eqref{est Gamma-1}. 
\endproof
\end{lem}

\begin{prop} [Estimates on $v^{\Sigma}_{\Lambda^c}$] \label{Estimates on v Lambda complementary}
There exist $s\in \mathbb{N}$ and a constant $C>0$ independent of $h$ such that $v^{\Sigma}_{\Lambda^c}$ can be considered as a remainder $R(v)$ satisfying \eqref{estimates L^2 to L^2 of R(v)}, \eqref{estimates L^2 to L^inf of R(v)}.
\proof
Firstly, we want to reduce to the study of the action of $Op_h^w(1-\Gamma)$ on $v$ and not on $v^{\Sigma}$, so we can use lemma \ref{lemma GammaTilde} with $\widetilde{\Gamma}= 1-\gamma$, obtaining
\begin{equation}
Op_h^w\Big((1-\gamma)(\frac{x+p'(\xi)}{\sqrt{h}})\Big)v^{\Sigma}= Op_h^w\Big(\Sigma(\xi)\chi(h^{\beta}\xi)(1-\gamma)(\frac{x+p'(\xi)}{\sqrt{h}})\Big)v + R(v) \,,
\end{equation}
where $R(v)$ satisfies \eqref{estimates L^2 to L^2 of R(v)}, \eqref{estimates L^2 to L^inf of R(v)}.
Then it remains to analyse $$Op_h^w\left(\Sigma(\xi)\chi(h^{\beta}\xi)(1-\gamma)(\frac{x+p'(\xi)}{\sqrt{h}})\right)v \,.$$ 
We write the symbol of the operator as $\Sigma(\xi)\chi(h^{\beta}\xi) \Gamma_{-1}(\frac{x+p'(\xi)}{\sqrt{h}})(\frac{x+p'(\xi)}{\sqrt{h}})$, with $\Gamma_{-1}(\xi) = \frac{(1-\gamma)(\xi)}{\xi}$, and we can apply the previous lemma with $c(x,\xi)= \Sigma(\xi)\chi(h^{\beta}\xi)\in h^{-\sigma}S_{0,\beta}(1)$, $\sigma$ as in \eqref{cases for sigma}, $c'(x,\xi)\equiv 1$, to derive that it is a remainder $R(v)$ satisfying \eqref{estimates L^2 to L^2 of R(v)}, \eqref{estimates L^2 to L^inf of R(v)}.
\endproof
\end{prop}

\vspace{0.5cm}
\noindent We want to apply first $Op_h^w(\Sigma(\xi))$ to \eqref{half KG in v with multi-operators}.
As $Op_h^w(\Sigma(\xi))$ commutes with $D_t - Op_h^w(x\xi + p(\xi))$ (because $\Sigma(D)$ commutes with $D_t - p(D)$), we obtain the equation:
\begin{equation} \label{half KG with Sigma}
(D_t - Op_h^w(x\xi + p(\xi)))v^{\Sigma} = h Op_h^w(\Sigma) \Big[\sum_{|I|=3}Op_h(m_I)(v_I) + \sum_{|I|=3}Op_h(\widetilde{m}_I)(v_I)\Big] \, .
\end{equation}
Then, we apply also $Op_h^w(\Gamma)$ to \eqref{half KG with Sigma}, so we have to calculate its commutator with the linear part of the equation, as done in the following:

\begin{lem} \label{commutation of Op(gamma) with linear part}
\begin{equation}
\big[ D_t - Op_h^w(x\xi + p(\xi)), Op_h^w(\Gamma(x,\xi))\big] = Op_h^w(b) \, , 
\end{equation}
where
\begin{equation}
b(x,\xi) = h \Gamma_{-1}(\frac{x+p'(\xi)}{\sqrt{h}})(\frac{x+p'(\xi)}{\sqrt{h}}) + h^{\frac{3}{2}}r \, ,
\end{equation}
$r \in S_{\frac{1}{2},0}(\langle\frac{x+p'(\xi)}{\sqrt{h}}\rangle^{-1})$, and $\Gamma_{-1}$ satisfies $|\partial^{\alpha}\Gamma_{-1}(\xi)|\lesssim \langle \xi \rangle^{-1-\alpha}$.
\proof
First we start by calculating $[D_t, Op_h^w(\Gamma)]= D_t Op_h^w(\Gamma) - Op_h^w(\Gamma)D_t$ : 
\begin{equation} \label{Dt vLambda}
\begin{split}
D_t Op_h^w(\Gamma) v & = \frac{1}{i}\partial_t \left[\frac{1}{2\pi} \int_{\mathbb{R}} \int_{\mathbb{R}} e^{i(x-y)\xi} \gamma(\frac{\frac{x+y}{2} + p'(h \xi)}{\sqrt{h}})v(t,y) \, dy d\xi  \right] \\
& = \frac{-h^2}{2\pi i} \partial_h \left[\frac{1}{2\pi} \int_{\mathbb{R}} \int_{\mathbb{R}} e^{i(x-y)\xi} \gamma(\frac{\frac{x+y}{2} + p'(h \xi)}{\sqrt{h}})v(t,y) \, dy d\xi  \right] \\
& = -\frac{h}{2\pi i} \int_{\mathbb{R}} \int_{\mathbb{R}} e^{i(x-y)\xi}\,\gamma'(\frac{\frac{x+y}{2} + p'(h \xi)}{\sqrt{h}}) \frac{p''(h\xi)h\xi}{\sqrt{h}} v(t,y) \, dy d\xi \\ 
& + \frac{h}{4 \pi i}\int_{\mathbb{R}} \int_{\mathbb{R}} e^{i(x-y)\xi}\,\gamma'(\frac{\frac{x+y}{2} + p'(h \xi)}{\sqrt{h}})(\frac{\frac{x+y}{2} + p'(h \xi)}{\sqrt{h}}) v(t,y) \, dy d\xi \\
& + \frac{1}{2\pi}\int_{\mathbb{R}} \int_{\mathbb{R}} e^{i(x-y)\xi} \gamma(\frac{\frac{x+y}{2} + p'(h \xi)}{\sqrt{h}})D_tv(t,y) \, dy d\xi  \\
& = ih\, Op_h^w\left(\gamma'(\frac{x+p'(\xi)}{\sqrt{h}})(\frac{p''(\xi)\xi}{\sqrt{h}})\right)v - \frac{ih}{2}\, Op_h^w \left(\gamma'(\frac{x+p'(\xi)}{\sqrt{h}})(\frac{x+p'(\xi)}{\sqrt{h}}) \right)v \\
& + Op_h^w(\Gamma)D_tv \, .
\end{split}
\end{equation}
Then, using \eqref{r_k} and \eqref{symbol of commutator } we write
\begin{equation} \label{brackets}
[Op_h^w(\Gamma(x,\xi)), Op_h^w(x\xi + p(\xi))] = \frac{h}{i}Op_h^w\Big(\big\{\gamma(\frac{x+p'(\xi)}{\sqrt{h}}), x\xi + p(\xi)\big\}\Big) + r_2 \, ,
\end{equation}
with $r_2 \in h^{\frac{3}{2}}S_{\frac{1}{2},0}(\langle\frac{x+p'(\xi)}{\sqrt{h}} \rangle^{-1})$ from lemma \ref{lem of composition with Gamma Tilde}, since $\partial^\alpha\Gamma \in h^{-\frac{|\alpha|}{2}}S_{\frac{1}{2},0}(\langle \frac{x+p'(\xi)}{\sqrt{h}}\rangle^{-1})$, $\partial^\alpha(x\xi+p'(\xi))\in S_{0,0}(1)$ for $|\alpha|=3$.
On the other hand, developing the braces in \eqref{brackets} we find
\begin{equation*}
\begin{split}
\frac{h}{i} Op_h^w\Big(\big\{ \gamma(\frac{x+p'(\xi)}{\sqrt{h}}), x\xi +p(\xi)\big\}\Big) & = -ih \, Op_h^w\left(\gamma'(\frac{x+p'(\xi)}{\sqrt{h}})\frac{p''(\xi)\xi}{\sqrt{h}}\right) \\
& + i h \, Op^w_h\left(\gamma'(\frac{x+p'(\xi)}{\sqrt{h}})(\frac{x+p'(\xi)}{\sqrt{h}})\right)\,,
\end{split}
\end{equation*}
so when we add it to $[D_t, Op_h^w(\Gamma)]$ calculated before, we obtain the result just choosing $\Gamma_{-1}(\xi) = \frac{1}{2} \gamma'(\xi)$.
\endproof
\end{lem}

\vspace{0.5 cm}
\noindent We apply $Op_h^w(\Gamma)$ to equation \eqref{half KG with Sigma}. Using lemma \ref{commutation of Op(gamma) with linear part}, we obtain
\begin{equation}
\begin{split}
(D_t - Op_h^w(x\xi +p(\xi))v^{\Sigma}_{\Lambda} & = h\, Op_h^w(\Gamma) Op_h^w(\Sigma) \Big[\sum_{|I|=3}Op_h(m_I)(v_I) + \sum_{|I|=3}Op_h(\widetilde{m}_I)(v_I)\Big] \\
& + h Op_h^w\left(\Gamma_{-1}(\frac{x+p'(\xi)}{\sqrt{h}})(\frac{x+p'(\xi)}{\sqrt{h}})\right)v^{\Sigma} + h^{\frac{3}{2}}Op_h^w(r)v^{\Sigma} \, , 
\end{split}
\end{equation}
$r \in S_{\frac{1}{2},0}(\langle \frac{x+p'(\xi)}{\sqrt{h}}\rangle^{-1})$, where the second and third term in the right hand side are of the form $hR(v)$, $R(v)$ satisfying the estimates \eqref{estimates L^2 to L^2 of R(v)},\eqref{estimates L^2 to L^inf of R(v)}.
In fact, using lemma \ref{lemma GammaTilde} with $\widetilde{\Gamma}(\xi)=\Gamma_{-1}(\xi)\xi$, and lemma \ref{composition gamma-1 and its argument}, we have
\begin{equation}
\begin{split}
Op_h^w\left(\Gamma_{-1}(\frac{x+p'(\xi)}{\sqrt{h}})(\frac{x+p'(\xi)}{\sqrt{h}})\right)v^{\Sigma} & = Op_h^w\left(\Sigma(\xi)\chi(h^{\beta}\xi)\Gamma_{-1}(\frac{x+p'(\xi)}{\sqrt{h}})(\frac{x+p'(\xi)}{\sqrt{h}})\right)v + R(v) \\
& = R(v) \,,
\end{split}
\end{equation}
while $r$ can be split via a function $\chi$ as in lemma \ref{new estimate 1-chi}, with $\beta>0$, obtaining $r(x,\xi)\chi(h^{\beta}\xi) \in S_{\frac{1}{2},\beta}(\langle \frac{x+p'(\xi)}{\sqrt{h}}\rangle^{-1})$ to which we can apply proposition \ref{Continuity from $L^2$ to L^inf}, and $r(x,\xi)(1-\chi)(h^{\beta}\xi)$ to which can be instead applied lemma \ref{new estimate 1-chi}.
Then also $h^{\frac{3}{2}}Op_h^w(r)v^{\Sigma}=hR(v)$.
Therefore, the equation we want to deal with becomes
\begin{equation} \label{KG for v-lambda}
(D_t - Op_h^w(x\xi +p(\xi))v^{\Sigma}_{\Lambda}  = h\, Op_h^w(\Gamma)Op_h^w(\Sigma) \Big[\sum_{|I|=3}Op_h(m_I)(v_I) + \sum_{|I|=3}Op_h(\widetilde{m}_I)(v_I)\Big] + h R(v) \, ,
\end{equation}
with a remainder $R(v)$ which satisfies \eqref{estimates L^2 to L^2 of R(v)}, \eqref{estimates L^2 to L^inf of R(v)}.

\vspace{0.5cm}
\subsection{Development at $\xi = d\varphi(x)$} 
The next step now is to develop the symbol of the \emph{characteristic} term in the nonlinearity, i.e. the one corresponding to $I=(1,1,-1)$, at $\xi = d\varphi(x)$. 
This will allow us to write it from $|v^{\Sigma}_{\Lambda}|^2v^{\Sigma}_{\Lambda}$ up to some remainders, transforming the action of pseudo-differential operators acting on it into a product of smooth functions of $x$.
Such development is not essential on \emph{non characteristic} terms, which will be eliminated through a normal form argument in the next section.
However, some results, such as proposition \ref{prop Op(a) develop } and lemma \ref{dev of Sigma}, apply suitably also to \emph{non characteristic} terms, so we will freely make use of them to get some simplifications.

\noindent We want to prove the following result:

\begin{prop} \label{development of nonlinearity}
Suppose that $v$ satisfies the a priori estimates in \eqref{a priori estimates on v}, \eqref{a priori estimate on Lv and v}. 
Then there exists a family of functions $\theta_h(x)\in C^{\infty}_0(]-1,1[)$, real valued, equal to one on an interval $[-1+ch^{2\beta},1-ch^{2\beta}]$, $\|\partial^{\alpha}\theta_h\|_{L^{\infty}}=O(h^{-2\beta\alpha})$, such that  
the nonlinearity in \eqref{KG for v-lambda} can be written as
\begin{equation}
\begin{split}
& h \Phi_1^{\Sigma}(x)\theta_h(x)|v^{\Sigma}_{\Lambda}|^2v^{\Sigma}_{\Lambda} + h Op_h^w(\Gamma)\left[\Phi_3^{\Sigma}(x)\theta_h(x)(v^{\Sigma}_{\Lambda})^3 + \Phi_{-1}^{\Sigma}(x)\theta_h(x) |v^{\Sigma}_{\Lambda}|^2\overline{v^{\Sigma}_{\Lambda}} + \Phi_{-3}^{\Sigma}(x)\theta_h(x) (\overline{v^{\Sigma}_{\Lambda}})^3\right] \\
& + hR(v) \,,
\end{split}
\end{equation}
where $\Phi_j^{\Sigma}(x)$ are smooth functions of $x$, $|\Phi_j^{\Sigma}(x)|=O(h^{-\sigma})$ on the support of $\theta_h$, for $j \in \{3,1,-1,-3\}$, $\sigma>0$ small, and where the remainder $R(v)$ satisfies estimates \eqref{estimates L^2 to L^2 of R(v)}, \eqref{estimates L^2 to L^inf of R(v)}.
\end{prop}

\noindent Before proving this proposition, we need the following general result.

\begin{prop} \label{prop Op(a) develop }
Let $a(x,\xi)$ be a real symbol in $S_{\delta,0}(\langle\xi\rangle^q)$, $q \in\mathbb{R}$, for some $\delta>0$ small.
There exists a family $(\theta_h(x))_h$ of $C^{\infty}$ functions, real valued, supported in some interval $[-1+ch^{2\beta},1-ch^{2\beta}]$, for a small $\beta>0$, with $(h\partial_h)^k\theta_h$ bounded for every $k$, such that
\begin{equation} \label{Op(a)v development}
Op_h^w(a)v = \theta_h(x)a(x,d\varphi(x))v + R_1(v) \,,
\end{equation}
where $R_1(v)$ is a remainder satisfying estimates \eqref{L^2-L^2 estimate of R}, \eqref{L^2-L^inf estimate of R}, with $\sigma=\sigma(\beta,\delta)>0$, $\sigma\rightarrow 0$ as $\beta, \delta \rightarrow 0$.
The same equality holds replacing $v$ by $v^{\Sigma}$.
\proof
In order to develop the symbol $a(x,\xi)$ at $\xi = d\varphi(x)$ we need to stay away from points $x=\pm 1$, so the idea is to truncate near $\Lambda$ and to estimate different terms coming out.

\vspace{0.5cm}
\noindent First, let us consider a function $\chi \in C^{\infty}_0(\mathbb{R})$ as in lemma \ref{new estimate 1-chi}, $\beta>0$.
We decompose $a(x,\xi)$ as follows
\begin{equation}
a(x,\xi)= a(x,\xi) \chi(h^{\beta}\xi) + a(x,\xi)(1-\chi)(h^{\beta}\xi) \,.
\end{equation}
It turns out from symbolic calculus, proposition \ref{Continuity from $L^2$ to $L^2$}, lemma \ref{new estimate 1-chi} and Sobolev injection \eqref{semiclassical Sobolev injection}, that $Op_h^w(a(x,\xi)(1-\chi)(h^{\beta}\xi))v$ is of the form $R_1(v)$, if we choose $s\gg 1$ sufficiently large, so we have just to deal with $a(x,\xi)\chi(h^{\beta}\xi)$.
Since this symbol is rapidly decaying in $|h^{\beta}\xi|$, we notice that, to prove that the estimate \eqref{L^2-L^2 estimate of R} holds for terms of interest, we can turn the $H^{\rho}_h$ norm into the $L^2$ norm up to a small loss in $h$, and then simply estimate the $L^2$ norm of these terms.
This is obvious when $\rho <0$, for $H^{\rho}_h$ injects in $L^2$, while for $\rho \in \mathbb{N}$ this follows using the definition \ref{def of Sobolev spaces} (\ref{def of H^s_h}), symbolic calculus, and the fact that $\langle\xi \rangle^\rho \chi(h^{\beta}\xi) \le h^{-\rho\beta}$.
Therefore, it is sufficient for our aim to prove that terms coming out are remainders $R(v)$, in the sense of inequalities \eqref{estimates L^2 to L^2 of R(v)}, \eqref{estimates L^2 to L^inf of R(v)}. \\
Secondly, we consider a smooth compactly supported function $\gamma \in C^{\infty}_0(\mathbb{R})$, $\gamma\equiv 1$ in a neighbourhood of zero, with a sufficiently small support, and we split $a(x,\xi)\chi(h^{\beta}\xi)$ as
\begin{equation}
a(x,\xi)\chi(h^{\beta}\xi) = a(x,\xi)\chi(h^{\beta}\xi)\gamma(\frac{x+p'(\xi)}{\sqrt{h}}) + a(x,\xi)\chi(h^{\beta}\xi)(1-\gamma)(\frac{x+p'(\xi)}{\sqrt{h}})\, .
\end{equation}
Again, the second symbol in the right hand side gives us a remainder.
In fact, since $(1-\gamma)(\xi)$ is supported for $|\xi|>\alpha_1$, we can write
\begin{equation}
a(x,\xi)\chi(h^{\beta}\xi)(1-\gamma)(\frac{x+p'(\xi)}{\sqrt{h}})= a(x,\xi)\chi(h^{\beta}\xi) \Gamma_{-1}(\frac{x+p'(\xi)}{\sqrt{h}}) (\frac{x+p'(\xi)}{\sqrt{h}}) \, , 
\end{equation}
where $\Gamma_{-1}(\xi)= \frac{(1-\gamma)(\xi)}{\xi}$, $|\partial^{\alpha}\Gamma_{-1}(\xi)|\lesssim \langle \xi\rangle^{-1-\alpha}$.
Lemma \ref{composition gamma-1 and its argument} with $c(x,\xi)=a(x,\xi)\chi(h^{\beta}\xi) \in h^{-\sigma}S_{\delta,\beta}(1)$, $\sigma \ge 0$ small (either equal to $q\beta$ for $q\ge 0$, or to 0 for $q< 0$), $c'(x,\xi) \equiv 1$, implies that 
$Op_h^w\left( a(x,\xi)\chi(h^{\beta}\xi) \Gamma_{-1}(\frac{x+p'(\xi)}{\sqrt{h}})(\frac{x+p'(\xi)}{\sqrt{h}})\right)v$ satisfies \eqref{estimates L^2 to L^2 of R(v)}, \eqref{estimates L^2 to L^inf of R(v)}.

\vspace{0.5cm}
\noindent The last remaining symbol is $a(x,\xi)\chi(h^{\beta}\xi)\gamma(\frac{x+p'(\xi)}{\sqrt{h}})$, which is supported in $\{(x,\xi)\in \mathbb{R}\times \mathbb{R}\, \big|\, |\xi|<C_2h^{-\beta}\, , |\frac{x+p'(\xi)}{\sqrt{h}}|<\alpha_2 \}$, so $x$ is bounded in a compact subset of $]-1,1[$ of the form $[-1 + ch^{2\beta}, 1-ch^{2\beta}]$, for a suitable positive constant $c$.
We may find a smooth function $\theta_h(x) \in C^{\infty}_0(]-1,1[)$, $\theta_h \equiv 1$ on $[-1 + ch^{2\beta}, 1-ch^{2\beta}]$, satisfying $\|\partial^{\alpha}\theta_h\|_{L^{\infty}}=O(h^{-2\beta\alpha})$, and write
\begin{equation}
a(x,\xi)\chi(h^{\beta}\xi) \gamma(\frac{x+p'(\xi)}{\sqrt{h}}) = a(x,\xi) \theta_h(x)\chi(h^{\beta}\xi) \gamma(\frac{x+p'(\xi)}{\sqrt{h}}) \, .
\end{equation}
Since on the support of $\theta_h$ we are away from $x=\pm 1$, we may now develop $a(x, \xi)$ at $\xi = d\varphi(x)$, 
\begin{equation}
\begin{split}
a(x,\xi) & = a(x, d\varphi(x)) + \int_0^1 \partial_\xi a(x, t\xi + (1-t)d\varphi(x))  dt \, (\xi -d\varphi(x)) \\
& = a(x, d\varphi(x)) + b(x,\xi) (x+p'(\xi)) \, ,
\end{split}
\end{equation}
where 
\begin{equation}
b(x,\xi) = \int_0^1 \partial_\xi a(x, t\xi + (1-t)d\varphi(x)) \, dt \, \frac{\xi -d\varphi(x)}{x+p'(\xi)} \, .
\end{equation}
Then,
\begin{equation} \label{develop at dfi}
\begin{split}
a(x,\xi) \theta_h(x)\chi(h^{\beta}\xi) \gamma(\frac{x+p'(\xi)}{\sqrt{h}}) & = a(x, d\varphi(x)) \theta_h(x) + a(x, d\varphi(x))\theta_h(x) \big[\chi(h^{\beta}\xi)\gamma(\frac{x+p'(\xi)}{\sqrt{h}})-1\big] \\
& + b(x,\xi) \chi(h^{\beta}\xi)\gamma(\frac{x+p'(\xi)}{\sqrt{h}})(x+p'(\xi))\, .
\end{split}
\end{equation}
Let us call $I_1$ and $I_2$ the Weyl quantizations respectively of the second and third term in the right hand side of \eqref{develop at dfi}. We want to show that they satisfy \eqref{estimates L^2 to L^2 of R(v)}, \eqref{estimates L^2 to L^inf of R(v)}. \\
\noindent First we analyse the third term in the right hand side of \eqref{develop at dfi}.
We may find another smooth function $\widetilde{\gamma}$, with a small support, such that
\begin{equation} \label{we can consider gamma tilde}
\chi(h^{\beta}\xi)\gamma(\frac{x+p'(\xi)}{\sqrt{h}}) = \chi(h^{\beta}\xi)\gamma(\frac{x+p'(\xi)}{\sqrt{h}}) \widetilde{\gamma}(\langle \xi \rangle^2 (x+p'(\xi))) \,.
\end{equation}
From $a \in S_{\delta,0}(\langle\xi\rangle^q)$ and lemma \ref{lem on e and etilde}, $B(x,\xi):= b(x,\xi)\chi(h^{\beta}\xi)\widetilde{\gamma}(\langle\xi\rangle^2(x+p'(\xi)))$ is an element of 
$$h^{-\delta} S_{2\beta,\beta}(\langle\xi\rangle^{3+q})\subset h^{-\sigma}S_{\delta',\beta}(1) \,, $$
for $\delta'=\max\{\delta, 2\beta\}$, $\sigma>0$ small depending on $\beta$ and $\delta$.
Moreover, $|\partial^{\alpha}\gamma(\xi)|\le \langle\xi\rangle^{-1-\alpha}$, so lemma \ref{composition gamma-1 and its argument} implies that $Op_h^w\big(B(x,\xi)\gamma(\frac{x+p'(\xi)}{\sqrt{h}})(x+p'(\xi))\big)$ is a remainder $h^{\frac{1}{2}}R(v)$.
\\
On the other hand, $I_1$ has a symbol whose support is included in the union $\{|\xi|> C_1h^{-\beta}\}\cup \{|\frac{x+p'(\xi)}{\sqrt{h}}|>\alpha_1 \}$.
Take $\widetilde{\chi}\in C^{\infty}_0(\mathbb{R})$, $\widetilde{\chi}\equiv 1$ in a neighbourhood of zero, $supp\widetilde{\chi}\subset \{|\xi|<C_1 h^{-\beta}\}$, so that $\chi \widetilde{\chi}\equiv \widetilde{\chi}$.
We make a further decomposition,
\begin{equation} \label{further decomposition with chi tilde}
\begin{split}
& \chi(h^{\beta}\xi) \gamma(\frac{x+p'(\xi)}{\sqrt{h}}) - 1  =\\
& = \left(\chi(h^{\beta}\xi) \gamma(\frac{x+p'(\xi)}{\sqrt{h}}) - 1 \right) \widetilde{\chi}(h^{\beta}\xi) + \left( \chi(h^{\beta}\xi) \gamma(\frac{x+p'(\xi)}{\sqrt{h}}) - 1\right)(1-\widetilde{\chi})(h^{\beta}\xi) \\
& = \left( \gamma(\frac{x+p'(\xi)}{\sqrt{h}})-1\right)\widetilde{\chi}(h^{\beta}\xi)+ \left( \chi(h^{\beta}\xi) \gamma(\frac{x+p'(\xi)}{\sqrt{h}}) - 1\right)(1-\widetilde{\chi})(h^{\beta}\xi) \, .
\end{split}
\end{equation}
The first term in the right hand side is supported for $|\frac{x+p'(\xi)}{\sqrt{h}}|>\alpha_1$, so it can be written as
\begin{equation*}
\widetilde{\chi}(h^{\beta}\xi) \Gamma_{-1}(\frac{x+p'(\xi)}{\sqrt{h}})(\frac{x+p'(\xi)}{\sqrt{h}}) \, ,
\end{equation*}
and it is a remainder by lemma \ref{composition gamma-1 and its argument}.
Besides, the second term in the right hand side is supported for $|\xi|>C_1 h^{-\beta}$, so it is essentially an application of lemma \ref{new estimate 1-chi} to show that it is a remainder $R(v)$.
This shows finally that the development in \eqref{Op(a)v development} holds.
For the last statement of the proposition, one can show that the same proof we did for $v$ can be applied for $v^{\Sigma}$, just changing $a(x,\xi)$ into $a(x,\xi)\Sigma(\xi)$ trough lemma \ref{lemma GammaTilde} when needed, and for a new $\sigma>0$ depending also on $\rho$.  

\endproof
\end{prop}

\proof[Proof of Proposition \ref{development of nonlinearity}]
The idea of the proof is to develop all symbols $m_I(\xi), \widetilde{m}_I(\xi)$ occurring in the cubic nonlinearity at $\xi = d\varphi(x)$ using the previous proposition.
Remark that, when $i_k=-1$, $v_{i_k}=\bar{v}$ and we have 
\begin{equation}
Op_h(m_k(\xi))\bar{v} = \overline{Op_h(m_k(-\xi))v} = \overline{Op_h(m_k(i_k \xi))v} \, ,
\end{equation}
so the development at $\xi=d\varphi(x)$ will give us a term $m_k(i_k d\varphi(x))v_{i_k}$, since $m_k(\xi), d\varphi(x)$ are real valued.

\vspace{0.5cm}
\noindent We first write $Op_h^w(m_k(\xi))v_{i_k}=Op_h^w(m_k(\xi)\Sigma(\xi)^{-1})v_{i_k}^{\Sigma}=Op_h^w(m_k^{\Sigma}(\xi))v^{\Sigma}$ (following the notation introduced in \eqref{def of m_I Sigma}) and then we apply proposition \ref{prop Op(a) develop }. 
From bounds \eqref{a priori estimates on v}, \eqref{a priori estimate on Lv and v}, we have $\|m^{\Sigma}_{k}(i_kd\varphi(x))\theta_h(x) v_{i_k}^{\Sigma}\|_{L^{\infty}}=O(h^{-\sigma})$, $\|m^{\Sigma}_{k}(i_kd\varphi(x))\theta_h(x) v_{i_k}^{\Sigma}\|_{H^{\rho}_h}=O(h^{-\sigma})$, for a $\sigma>0$ depending on $\beta$, so we immediately obtain that $Op_h(m_I)(v_I)= \displaystyle\prod_{k=1}^3m^{\Sigma}_{k}(i_kd\varphi(x))\theta_h(x) v_{i_k}^{\Sigma} + R_1(v)$, $R_1(v)$ satisfying estimates \eqref{L^2-L^2 estimate of R}, \eqref{L^2-L^inf estimate of R}, and we can perform the passage from the term
\begin{equation}
\sum_{|I|=3}Op_h(m_I)(v_I) + \sum_{|I|=3}Op_h(\widetilde{m}_I)(v_I)
\end{equation}
to its development
\begin{equation}
\sum_{|I|=3} m_I^{\Sigma}(d\varphi_I(x))\theta_h(x) v^{\Sigma}_I +  \sum_{|I|=3} \widetilde{m}^{\Sigma}_I(d\varphi_I(x))\theta_h(x) v^{\Sigma}_I + R_1(v)\,.
\end{equation}
The nonlinearity in \eqref{KG for v-lambda} becomes
\begin{equation}
\begin{split}
& h\, Op_h^w(\Gamma)  Op_h^w(\Sigma(\xi))\left[ \sum_{|I|=3} m_I^{\Sigma}(d\varphi_I(x))\theta_h(x) v^{\Sigma}_I + \sum_{|I|=3} \widetilde{m}_I^{\Sigma}(d\varphi_I(x))\theta_h(x) v^{\Sigma}_I \right] \\
&+ hOp_h^w(\Gamma)Op_h^w(\Sigma(\xi))R_1(v) \,,
\end{split}
\end{equation} 
where $R_1(v)$ satisfies \eqref{L^2-L^2 estimate of R}, so that $Op_h^w(\Gamma)Op_h^w(\Sigma(\xi))R_1(v)$ is a remainder of the form $R(v)$, satisfying the estimates \eqref{estimates L^2 to L^2 of R(v)}, \eqref{estimates L^2 to L^inf of R(v)}, by propositions \ref{Continuity from $L^2$ to $L^2$} and \ref{Continuity from $L^2$ to L^inf}. \\
The following three lemmas allow us to conclude the proof.
In particular, we underline that in lemma \ref{dev of Sigma} we only need an $L^2$ estimate on what we denote $R(v)$, because we will apply to it the operator $Op_h^w(\Gamma)$, which is continuous from $L^2$ to $L^{\infty}$ with norm $\|Op_h^w(\Gamma)\|_{\mathcal{L}(L^2;L^{\infty})}=O(h^{-\frac{1}{4}-\sigma})$ by proposition \ref{Continuity from $L^2$ to L^inf}.
\endproof

\begin{lem} \label{dev of Sigma}
Let $I=(i_1,i_2,i_3)$, $i_k \in \{1,-1 \}$ for $k=1,2,3$, be a fixed vector.
Denote by $A(\xi)$ the function $\Sigma(\xi)\chi(h^{\beta}\xi)$, with $\chi$ as in lemma \ref{new estimate 1-chi}, $\beta>0$.
Then
\begin{equation}
\begin{split}
Op_h^w(\Sigma(\xi)) \left(m^{\Sigma}_I(d\varphi_I(x))\theta_h(x)v^{\Sigma}_I\right) & = A\Big(\sum_{l=1}^3i_ld\varphi(x)\Big)m^{\Sigma}_I(d\varphi_I(x))\theta_h(x)v^{\Sigma}_I + h^{\frac{1}{2}}R(v) \,, \\
Op_h^w(\Sigma(\xi)) \left(\widetilde{m}^{\Sigma}_I(d\varphi_I(x))\theta_h(x)v^{\Sigma}_I \right)& = A\Big(\sum_{l=1}^3i_ld\varphi(x)\Big)\widetilde{m}^{\Sigma}_I(d\varphi_I(x))\theta_h(x)v^{\Sigma}_I + h^{\frac{1}{2}}R(v) \,,
\end{split}
\end{equation}
where $R(v)$ satisfies the estimate \eqref{estimates L^2 to L^2 of R(v)}.
\proof
Before proving the result, let us make some observations: first, in all the proof we will use generically the letter $\sigma$ to denote a small non-negative constant depending on $\beta$, that goes to zero when $\beta$ goes to zero; 
the symbol $\Sigma(\xi)$ can be reduced to $\Sigma(\xi)\chi(h^{\beta}\xi) \in h^{-\sigma}S_{0,\beta}(1)$, $\sigma$ as in \eqref{cases for sigma}, up to remainders (essentially using lemma \ref{new estimate 1-chi}); from the \emph{a priori} estimates \eqref{a priori estimates on v}, \eqref{a priori estimate on Lv and v} made on $v$, we have $\|m^{\Sigma}_I(d\varphi_I(x))\theta_h(x)v^{\Sigma}_I\|_{L^2}=O(h^{-\sigma})$.

\vspace{0.5cm}
\noindent Let us consider a smooth function $\tilde{\theta}_h(x) \in C^{\infty}_0(]-1,1[)$, such that $\tilde{\theta}_h \theta_h \equiv \theta_h$, and let us write
$$m^{\Sigma}_I(d\varphi_I(x))\theta_h(x)v^{\Sigma}_I = \tilde{\theta}_h(x) m^{\Sigma}_I(d\varphi_I(x))\theta_h(x)v^{\Sigma}_I \,.$$ 
We enter $\tilde{\theta}_h(x)$ in $Op_h^w(\Sigma(\xi)\chi(h^{\beta}\xi))$ applying symbolic calculus of proposition \ref{a sharp b}, to be able to develop the symbol at $\xi = \displaystyle\sum_{l=1}^3i_ld\varphi(x)$.
We have
\begin{equation}
\Sigma(\xi)\chi(h^{\beta}\xi)\sharp \tilde{\theta}_h(x)= \Sigma(\xi)\chi(h^{\beta}\xi)\tilde{\theta}_h(x) + r_0 \,,
\end{equation} with $r_0 \in h^{1-\sigma}S_{\delta,\beta}(1)$, $\delta>0$  small, so proposition \ref{Continuity from $L^2$ to $L^2$} implies that its quantization gives a remainder as in the statement, when applied to $m^{\Sigma}_I(d\varphi_I(x))\theta_h(x)v^{\Sigma}_I$.
Now, since we are away from $x=\pm 1$, we can develop $A(\xi)=\Sigma(\xi)\chi(h^{\beta}\xi)$ at $\xi = \displaystyle\sum_{l=1}^3 i_l d\varphi(x)$ by Taylor's formula, i.e.
\begin{equation}
A(\xi)= A\Big(\sum_{l=1}^3 i_l d\varphi(x)\Big) + A'(x,\xi)(\xi - \sum_{l=1}^3 i_l d\varphi(x))\, ,
\end{equation}
with
\begin{equation}
A'(x,\xi)= \int_0^1 A'\Big(t\xi + (1-t)\sum_{l=1}^3 i_l d\varphi(x))\Big) \, dt \,,
\end{equation}
$A'(x,\xi)\tilde{\theta}_h(x)$ belonging to $h^{-\sigma}S_{\delta,0}(1)$.
To conclude the proof, we need to show that also
$Op_h^w\Big(A'(x,\xi)\tilde{\theta}_h(x)(\xi - \displaystyle\sum_{l=1}^3 i_l d\varphi(x))\Big)\left(m^{\Sigma}_I(d\varphi_I(x))\theta_h(x)v^{\Sigma}_I\right)=h^{\frac{1}{2}}R(v)$.
So let us consider a new function $\tilde{\tilde{\theta}}_h(x) \in C^{\infty}_0(]-1,1[)$, such that $\tilde{\tilde{\theta}}_h\tilde{\theta}_h \equiv \tilde{\theta}_h$.
Since $\tilde{\theta}_h(\xi - \displaystyle\sum_{l=1}^3 i_l d\varphi(x)) \in h^{-\sigma}S_{\delta,0}(\langle\xi\rangle)$, and using symbolic calculus of proposition \ref{a sharp b}, we write
\begin{equation} 
A'(x,\xi)\tilde{\tilde{\theta}}_h(x) \sharp \Big( \tilde{\theta}_h(\xi - \sum_{l=1}^3 i_l d\varphi(x))\Big) = A'(x,\xi) \tilde{\theta}_h(x)(\xi - \sum_{l=1}^3 i_l d\varphi(x)) + r'_0 \, ,
\end{equation}
where $r'_0 \in h^{1-\sigma}S_{\delta,0}(1)$. Again proposition \ref{Continuity from $L^2$ to $L^2$} and the uniform bound on $v$ imply that $Op_h^w(r'_0) ( m^{\Sigma}_I(d\varphi_I(x))\theta_h(x) v_I^{\Sigma})$ is a remainder $h^{\frac{1}{2}}R(v)$.
We can focus on the term
\begin{equation} \label{two op con BSigma}
Op_h^w\big(A'(x,\xi)\tilde{\tilde{\theta}}_h(x) \big) Op_h^w\Big(\tilde{\theta}_h(x)(\xi -\sum_{l=1}^3 i_l d\varphi(x)) \Big) ( m^{\Sigma}_I(d\varphi_I(x))\theta_h(x) v_I^{\Sigma})\, ,
\end{equation}
and we can go further, limiting ourselves to consider the action of these operators when $v^{\Sigma}_I$ is replaced by
\begin{equation} \label{v^0_I def}
v^0_I:= \displaystyle\prod_{k=1}^3Op_h^w(\Sigma(\xi)\chi(h^{\beta}\xi))v_{i_k}\,,
\end{equation}
again up to terms with symbols supported for $|\xi|\ge h^{-\beta}$, which are remainders from lemma \ref{new estimate 1-chi}. 
The operator $Op_h^w\Big(\tilde{\theta}_h(x) (\xi -\displaystyle\sum_{l=1}^3 i_l d\varphi(x))\Big)$ has a symbol linear in $\xi$, so 
\begin{equation} \label{xi - dvarphi is linear in xi}
\begin{split}
Op_h^w\Big(\tilde{\theta}_h(x) (\xi -\sum_{l=1}^3 i_l d\varphi(x))\Big) & = \frac{1}{2} hD \tilde{\theta}_h(x) + \frac{1}{2}\tilde{\theta}_h(x) hD - \tilde{\theta}_h(x)\sum_{l=1}^3 i_l d\varphi(x) \\
& = h\frac{\tilde{\theta}'_h(x)}{2i} + \tilde{\theta}_h(x) (hD -\sum_{l=1}^3 i_l d\varphi(x))\,,
\end{split}
\end{equation}
and if we choose $\tilde{\theta}_h$ such that $\tilde{\theta}_h \theta_h \equiv \theta_h$, we have that $\tilde{\theta}'_h \equiv 0$ on the support of $\theta_h$, giving no contributions when $h\frac{\tilde{\theta}'_h(x)}{2}$ is multiplied by $m^{\Sigma}_I(d\varphi_I(x))\theta_h(x)v^0_I$.
Here $(hD-\displaystyle\sum_{l=1}^3 i_l d\varphi(x))$ acts like a derivation on $v^0_I$, so the Leibniz's rule holds and
\begin{equation} \label{Leibnitz action}
\begin{split}
& Op_h^w\Big(\tilde{\theta}_h(x) (\xi -\sum_{l=1}^3 i_l d\varphi(x))\Big)\left(m^{\Sigma}_I(d\varphi_I(x))\theta_h(x)v^0_I \right)  = \\
& = \tilde{\theta}_h(x) (hD-\sum_{l=1}^3 i_l d\varphi(x))\left(m^{\Sigma}_I(d\varphi_I(x))\theta_h(x)v^0_I\right) \\
& = h D(m^{\Sigma}_I(d\varphi_I(x))\theta_h(x))v^0_I  + m^{\Sigma}_I(d\varphi_I(x))\theta_h(x) \tilde{\theta}_h(x)(hD-\sum_{l=1}^3 i_l d\varphi(x))(v^0_I) \, .
\end{split}
\end{equation}
Then, if for instance $v^0_I=(v^0)^3$ (i.e. $I=(1,1,1)$, and the same idea can be applied with $|v^0|^2v^0$, $|v^0|^2\overline{v^0}$ and $(\overline{v^0})^3$), we derive
\begin{equation} \label{leibnitz rule on v^3}
\begin{split}
\tilde{\theta}_h(x)(hD- 3 d\varphi(x))(v^0)^3 & = 3(v^0)^2 \tilde{\theta}_h(x)(hD-d\varphi(x))v^0 \\
& = 3(v^0)^2 Op_h^w(\tilde{\theta}_h(x)(\xi-d\varphi(x)))v^0  - \frac{3}{2i}h \tilde{\theta}'_h(x)(v^0)^3\, ,
\end{split}
\end{equation}  
using the relation \eqref{xi - dvarphi is linear in xi} in the last equality (however, the second term in the right hand side disappears when we inject \eqref{leibnitz rule on v^3} in \eqref{Leibnitz action}).
Now we can re-express the first term in the right hand side from $h\mathcal{L}v^0$.
In fact, up to further remainders, $Op_h^w(\tilde{\theta}_h(x)(\xi-d\varphi(x)))v^0$ can be reduced to $Op_h^w(\tilde{\theta}_h(x)\chi(h^{\beta}\xi)(\xi-d\varphi(x)))v^0$, and this term can be split with the help of a $\gamma \in C^{\infty}_0(\mathbb{R})$, $\gamma \equiv 1$ in zero, namely
\begin{equation} \label{xi - dfi split con gamma}
\begin{split}
Op_h^w\left(\tilde{\theta}_h(x)\chi(h^{\beta}\xi)(\xi-d\varphi(x))\right)v^0 & = Op_h^w\left(\tilde{\theta}_h(x)\chi(h^{\beta}\xi)\gamma(\frac{x+p'(\xi)}{\sqrt{h}})(\xi-d\varphi(x))\right)v^0 \\
& + Op_h^w\left(\tilde{\theta}_h(x)\chi(h^{\beta}\xi)(1-\gamma)(\frac{x+p'(\xi)}{\sqrt{h}})(\xi-d\varphi(x))\right)v^0 \, .
\end{split}
\end{equation}
Both terms have an $L^2$ norm controlled from above by
$$C h^{1-\sigma}(\|\mathcal{L}v\|_{L^2} + \|v\|_{H^s_h})\, . $$
In fact, on one hand, we can take up the observation made in \eqref{we can consider gamma tilde}, and rewrite the first term in the right hand side as
\begin{equation}\label{split con etilde}
Op_h^w\left(\tilde{\theta}_h(x)\chi(h^{\beta}\xi)\gamma(\frac{x+p'(\xi)}{\sqrt{h}})\tilde{e}_+(x+p'(\xi))\right)v^0
\end{equation}
where $\tilde{e}_+$ is defined in \eqref{def of e and etilde}. 
On the other hand, also the symbol associated to the second operator in the right hand side 
can be rewritten highlighting the factor $(x+p'(\xi))$, as follows
$$ \tilde{\theta}_h(x)\chi(h^{\beta}\xi)\left(\frac{\xi - d\varphi(x)}{x+p'(\xi)}\right)(1-\gamma)(\frac{x+p'(\xi)}{\sqrt{h}})  (x+p'(\xi)) \,,$$
with $\tilde{\theta}_h(x)\chi(h^{\beta}\xi)\left(\frac{\xi - d\varphi(x)}{x+p'(\xi)}\right)(1-\gamma)(\frac{x+p'(\xi)}{\sqrt{h}})  \in h^{-\sigma}S_{\frac{1}{2},\beta}(1)$.
Then, to both operators we can apply lemma \ref{composition gamma-1 and its argument}, for $c(x,\xi)$ respectively equal to $\tilde{\theta}_h(x)\chi(h^{\beta}\xi)\tilde{e}_+$ and $\tilde{\theta}_h(x)\chi(h^{\beta}\xi)\left(\frac{\xi - d\varphi(x)}{x+p'(\xi)}\right)$, $c'(x,\xi)=\Sigma(\xi)\chi(h^{\beta}\xi)$, obtaining that they satisfy \eqref{est L2 in lemma 4.3 for GammaTilde}. 
Summing all up, together with \eqref{two op con BSigma}, \eqref{Leibnitz action},  \eqref{leibnitz rule on v^3}, the fact that $A'(x,\xi)\tilde{\tilde{\theta}}_h(x) \in h^{-\sigma}S_{\delta,0}(1)$, and propositions \ref{Continuity from $L^2$ to $L^2$}, we obtain the result of the lemma.
\endproof
\end{lem}

\noindent
From now on, we will denote by $\Phi^{\Sigma}_3(x), \Phi^{\Sigma}_1(x), \Phi^{\Sigma}_{-1}(x), \Phi^{\Sigma}_{-3}(x)$ respectively the coefficients of $(v^{\Sigma})^3, |v^{\Sigma}|^2v^{\Sigma}, |v^{\Sigma}|^2\overline{v^{\Sigma}}, (\overline{v^{\Sigma}})^3$. Since they come from $m^{\Sigma}_I(d\varphi_I(x))\theta_h(x), \widetilde{m}^{\Sigma}_I(d\varphi_I(x))\theta_h(x)$ which are $O(h^{-\sigma})$, for a small $\sigma>0$, they are also $O(h^{-\sigma})$.

\begin{lem}
With same notations as before,
\begin{equation}
Op_h^w(\Gamma)(\Phi_1^{\Sigma}(x)\theta_h(x)|v^{\Sigma}|^2v^{\Sigma}) = \Phi_1^{\Sigma}(x)\theta_h(x)|v^{\Sigma}|^2v^{\Sigma} + R(v) \,,
\end{equation}
where $R(v)$ satisfies estimates \eqref{estimates L^2 to L^2 of R(v)}, \eqref{estimates L^2 to L^inf of R(v)}.
\proof
Let us write $Op_h^w(\Gamma)=1- Op_h^w(1-\Gamma)$.
We want to show that the action of $Op_h^w(1-\Gamma)$ on $\Phi_1^{\Sigma}(x)\theta_h(x)|v^{\Sigma}|^2v^{\Sigma}$ gives us a remainder $R(v)$.
First, we can reduce the symbol $1-\Gamma$ to $(1-\Gamma)\chi(h^{\beta}\xi)$, with $\chi$ cut-off function, $\beta>0$, up to remainders from lemma \ref{new estimate 1-chi}.
Moreover, we can consider a smooth function $\tilde{\tilde{{\theta}}}_h(x) \in C^{\infty}_0(]-1,1[)$ such that $\tilde{\tilde{\theta}}_h\theta_h \equiv \theta_h$, and use symbolic calculus to enter $\tilde{\tilde{\theta}}_h(x)$ in $Op_h^w((1-\Gamma)\chi(h^{\beta}\xi))$ (again up to a remainder $R(v)$).
Then, we can write
\begin{equation}
(1-\Gamma)\chi(h^{\beta}\xi)\tilde{\tilde{\theta}}_h(x) = \frac{1}{\sqrt{h}} b(x,\xi) \Gamma_{-1}(\frac{x+p'(\xi)}{\sqrt{h}})\tilde{\theta}_h(x)(\xi - d\varphi(x))\,,
\end{equation}
where $b(x,\xi)= \chi(h^{\beta}\xi)\tilde{\tilde{\theta}}_h(x)(\frac{x+p'(\xi)}{\xi -d\varphi(x)}) \in h^{-\sigma}S_{\delta,\beta}(1)$, $\Gamma_{-1}(\xi)=\frac{(1-\gamma)(\xi)}{\xi}$, $\sigma, \delta>0$ small depending on $\beta$, and $\tilde{\theta}_h(x)$ a new smooth function in $C^{\infty}_0(]-1,1[)$, identically equal to 1 on the support of $\tilde{\tilde{\theta}}_h(x)$.
Applying symbolic calculus of lemma \ref{lem of composition with Gamma Tilde}, we derive
\begin{equation} \label{other calc}
\begin{split}
\frac{1}{\sqrt{h}} b(x,\xi) \Gamma_{-1}(\frac{x+p'(\xi)}{\sqrt{h}}) \sharp \tilde{\theta}_h(x)(\xi - d\varphi(x)) & = \frac{1}{\sqrt{h}} b(x,\xi) \Gamma_{-1}(\frac{x+p'(\xi)}{\sqrt{h}})\tilde{\theta}_h(x)(\xi - d\varphi(x))\\
& + \frac{\sqrt{h}}{2i}\left\{b(x,\xi) \Gamma_{-1}(\frac{x+p'(\xi)}{\sqrt{h}}), \tilde{\theta}_h(x)(\xi - d\varphi(x)) \right\} \\
& + r_1 \,,
\end{split}
\end{equation}
with $r_1\in h^{\frac{1}{2}-\sigma}S_{\frac{1}{2},\beta}(\langle\frac{x+p'(\xi)}{\sqrt{h}}\rangle^{-1})$, for a new small $\sigma>0$.
An explicit calculation, and the observation that $\tilde{\theta}_h' \equiv 0$ on the support of $\tilde{\tilde{\theta}}_h$, show that the Poisson bracket is equal to
\begin{equation} \label{calulation on brackets}
\begin{split}
& \Gamma_{-1}(\frac{x+p'(\xi)}{\sqrt{h}})\left[\tilde{\theta}_h(x)( -\partial_{\xi}b(x,\xi)d^2\varphi(x) - \partial_xb(x,\xi))\right] + \\
& + \Gamma'_{-1}(\frac{x+p'(\xi)}{\sqrt{h}})(\frac{x+p'(\xi)}{\sqrt{h}}) \chi(h^{\beta}\xi)\tilde{\theta}_h(x) \left[ \frac{-d^2\varphi(x)p''(\xi) - 1}{\xi -d\varphi(x)}\right] \, ,
\end{split}
\end{equation}
and since $x+p'(d\varphi)=0$, we have $-d^2\varphi(x) = \frac{1}{p''(d\varphi)}$ and 
\begin{equation}
\chi(h^{\beta}\xi)\tilde{\theta}_h(x)\left[\frac{-d^2\varphi(x)p''(\xi) - 1}{\xi -d\varphi(x)}\right] = \frac{\chi(h^{\beta}\xi)\tilde{\theta}_h(x)}{p''(d\varphi(x))}\int_0^1 p'''(t\xi + (1-t)d\varphi(x)) dt  \in h^{-\sigma}S_{\delta,\beta}(1)\,.
\end{equation}
Therefore, from $\Gamma_{-1}(\frac{x+p'(\xi)}{\sqrt{h}}), \,\Gamma'_{-1}(\frac{x+p'(\xi)}{\sqrt{h}})(\frac{x+p'(\xi)}{\sqrt{h}}) \in S_{\frac{1}{2},0}(\langle\frac{x+p'(\xi)}{\sqrt{h}}\rangle^{-1})$, other appearing symbols in \eqref{calulation on brackets} belonging to $h^{-\sigma}S_{\delta,\beta}(1)$,
we can rewrite the second and third term in the right hand side of \eqref{other calc} as $h^{\frac{1}{2}-\sigma}r$, with $r \in S_{\frac{1}{2},\beta}(\langle\frac{x+p'(\xi)}{\sqrt{h}}\rangle^{-1})$, so their action on $\Phi^{\Sigma}_1(x)\theta_h(x)|v^{\Sigma}|^2v^{\Sigma}$ gives us a remainder $R(v)$ by propositions \ref{Continuity from $L^2$ to $L^2$}, \ref{Continuity from $L^2$ to L^inf}.
In this way, we are reduce to estimate
\begin{equation} \label{op b etc}
\frac{1}{\sqrt{h}}Op_h^w\left(b(x,\xi)\Gamma_{-1}(\frac{x+p'(\xi)}{\sqrt{h}})\right) Op_h^w\big(\tilde{\theta}_h(x)(\xi -d\varphi(x)) \big) (\Phi^{\Sigma}_1(x)\theta_h(x)|v^{\Sigma}|^2v^{\Sigma})\, .
\end{equation}
Taking up \eqref{v^0_I def}, \eqref{xi - dvarphi is linear in xi}, \eqref{Leibnitz action} for $I=(1,1,-1)$, we obtain that $Op_h^w\big(\tilde{\theta}_h(x)(\xi -d\varphi(x))\big)$ acts like a derivation on its argument and 
\begin{equation}
\|Op_h^w\big(\tilde{\theta}_h(x)(\xi -d\varphi(x))\big) \Phi^{\Sigma}_1(x)\theta_h(x)|v^{\Sigma}|^2v^{\Sigma}\|_{L^2}\le C h^{1-\sigma}(\|\mathcal{L}v\|_{L^2} + \|v\|_{H^s_h}) \,,
\end{equation}
for a new small $\sigma>0$ still depending on $\beta$, so the fact that $b(x,\xi)\Gamma_{-1}(\frac{x+p'(\xi)}{\sqrt{h}})$ belongs to $S_{\frac{1}{2},\beta}(\langle \frac{x+p'(\xi)}{\sqrt{h}}\rangle^{-1})$, along with propositions \ref{Continuity from $L^2$ to $L^2$}, \ref{Continuity from $L^2$ to L^inf}, imply that the term in \eqref{op b etc} is a remainder $R(v)$ satisfying \eqref{estimates L^2 to L^2 of R(v)}, \eqref{estimates L^2 to L^inf of R(v)}.
This concludes the proof.
\endproof
\end{lem}

\vspace{0.5cm}
\noindent Proposition \ref{development of nonlinearity} allows us to arrive at the following equation 
\begin{equation}
\begin{split}
(D_t - Op_h^w(x\xi +p(\xi))v^{\Sigma}_{\Lambda} =  &\, h \Phi_1^{\Sigma}(x)\theta_h(x)|v^{\Sigma}|^2v^{\Sigma} +  h Op_h^w(\Gamma)\left[\Phi_3^{\Sigma}(x)\theta_h(x)(v^{\Sigma})^3 \right.\\
& \left. +\, \Phi_{-1}^{\Sigma}(x)\theta_h(x) |v^{\Sigma}|^2\overline{v^{\Sigma}} + \Phi_{-3}^{\Sigma}(x)\theta_h(x) (\overline{v^{\Sigma}})^3\right] + hR(v) \, ,
\end{split}
\end{equation}
which is not entirely in $v^{\Sigma}_{\Lambda}$, so to transform to the right equation we use the following lemma, whose proof comes directly from proposition \ref{Estimates on v Lambda complementary}, and this is the reason why we omit the details.

\begin{lem}
Under the same hypothesis as in theorem \ref{thm : From the PDE problem to an ODE}, there exists $s>0$ sufficiently large, and a constant $C>0$ independent of $h$, such that
\begin{align}
\| v^{\Sigma}_I - (v^{\Sigma}_{\Lambda})_I\|_{L^2} &\le C h^{\frac{1}{2}-\sigma}\left( \|\mathcal{L}v\|_{L^2} + \|v\|_{H^s_h}\right) \, ,\\
\| v^{\Sigma}_I - (v^{\Sigma}_{\Lambda})_I\|_{L^{\infty}} &\le C h^{\frac{1}{4}-\sigma}\left( \|\mathcal{L}v\|_{L^2} + \|v\|_{H^s_h}\right) \,, 
\end{align}
for a small $\sigma>0$ depending on $\beta$.
\end{lem}
\noindent Therefore $v^{\Sigma}_{\Lambda}$ is solution of the following equation :
\begin{equation} \label{ODE with the last pseudo op}
\begin{split}
(D_t - Op_h^w(x\xi +p(\xi))v^{\Sigma}_{\Lambda} = & \,h \Phi_1^{\Sigma}(x)\theta_h(x)|v^{\Sigma}_{\Lambda}|^2v^{\Sigma}_{\Lambda} + h Op_h^w(\Gamma)\left[\Phi_3^{\Sigma}(x)\theta_h(x)(v^{\Sigma}_{\Lambda})^3  \right.\\
& \left. +\,\Phi_{-1}^{\Sigma}(x)\theta_h(x) |v^{\Sigma}_{\Lambda}|^2\overline{v^{\Sigma}_{\Lambda}} + \Phi_{-3}^{\Sigma}(x)\theta_h(x) (\overline{v^{\Sigma}_{\Lambda}})^3\right] + hR(v)\,,
\end{split}
\end{equation}
$R(v)$ being a remainder satisfying estimates \eqref{estimates L^2 to L^2 of R(v)}, \eqref{estimates L^2 to L^inf of R(v)}.

\vspace{0.5cm}
\noindent To conclude this section, we develop $Op_h^w(x\xi + p(\xi))v^{\Sigma}_{\Lambda}$ at $\xi= d\varphi(x)$.

\begin{prop}
Under the same hypothesis as in theorem \ref{thm : From the PDE problem to an ODE}, 
\begin{equation}
Op_h^w(x\xi + p(\xi))v^{\Sigma}_{\Lambda}= \varphi(x)\theta_h(x)v^{\Sigma}_{\Lambda} + hR(v) \, ,
\end{equation}
where $R(v)$ satisfies the estimates in \eqref{estimates L^2 to L^2 of R(v)}, \eqref{estimates L^2 to L^inf of R(v)}.
Then, $v^{\Sigma}_{\Lambda}$ is solution of the following equation:
\begin{equation} \label{ODE for vLambda}
\begin{split}
D_t v^{\Sigma}_{\Lambda} & =\varphi(x) \theta_h(x) v^{\Sigma}_{\Lambda} + h \Phi_1^{\Sigma}(x)\theta_h(x)|v^{\Sigma}_{\Lambda}|^2v^{\Sigma}_{\Lambda} \\
& + h Op_h^w(\Gamma)\left[\Phi_3^{\Sigma}(x)\theta_h(x)(v^{\Sigma}_{\Lambda})^3 + \Phi_{-1}^{\Sigma}(x)\theta_h(x) |v^{\Sigma}_{\Lambda}|^2\overline{v^{\Sigma}_{\Lambda}} + \Phi_{-3}^{\Sigma}(x)\theta_h(x) (\overline{v^{\Sigma}_{\Lambda}})^3\right] \\
& + hR(v)\,,
\end{split}
\end{equation}
\proof
Consider a cut-off function $\chi$ as in lemma \ref{new estimate 1-chi}, and $\beta>0$.
One can split as follows
\begin{equation} \label{split con chi su vLambda}
v^{\Sigma}_{\Lambda}= Op_h^w(\chi(h^{\beta}\xi)\Gamma(x,\xi))v^{\Sigma} + Op_h^w((1-\chi)(h^{\beta}\xi)\Gamma(x,\xi))v^{\Sigma} \, , 
\end{equation}
and easily show that $Op_h^w(x\xi+p(\xi))Op_h^w((1-\chi)(h^{\beta}\xi)\Gamma(x,\xi))v^{\Sigma}$ is a remainder of the form $hR(v)$, $R(v)$ satisfying estimates \eqref{estimates L^2 to L^2 of R(v)}, \eqref{estimates L^2 to L^inf of R(v)}, just using symbolic calculus and lemma \ref{new estimate 1-chi}.

\vspace{0.5cm}
\noindent Therefore, we have to deal with $Op_h^w(x\xi +p(\xi))Op_h^w(\chi(h^{\beta}\xi)\Gamma(x,\xi))v^{\Sigma}$.
We have already observed that for $(x,\xi)$ in the support of $\chi(h^{\beta}\xi)\gamma(\frac{x+p'(\xi)}{\sqrt{h}})$, $x$ is bounded on a compact set $[-1+ch^{2\beta}, 1-ch^{2\beta}]$, which allows us to consider a smooth function $\theta_h(x) \in C^{\infty}_0(]-1,1[)$, identically equal to one on this interval, and then on the support of $\chi(h^\beta\xi)\gamma(\frac{x+p'(\xi)}{\sqrt{h}})$, so that:
\begin{equation}
\chi(h^{\beta}\xi)\gamma(\frac{x+p'(\xi)}{\sqrt{h}})= \theta_h(x)\chi(h^{\beta}\xi)\gamma(\frac{x+p'(\xi)}{\sqrt{h}}) \,.
\end{equation}
Moreover, the derivatives of $\theta_h$ are zero on the support of $\partial^{\alpha}(\chi(h^{\beta}\xi)\gamma(\frac{x+p'(\xi)}{\sqrt{h}}))$, for every multi-index $\alpha$.
This implies that products of $\theta'_h(x)$ with $\chi(h^{\beta}\xi)\gamma(\frac{x+p'(\xi)}{\sqrt{h}})$ and all its derivatives are always zero so, by lemma \ref{lem of composition with Gamma Tilde},
\begin{equation}
\theta_h(x) \sharp \chi(h^{\beta}\xi)\gamma(\frac{x+p'(\xi)}{\sqrt{h}}) = \theta_h(x)\chi(h^{\beta}\xi)\gamma(\frac{x+p'(\xi)}{\sqrt{h}}) + r_{\infty} \,,
\end{equation}
where $r_{\infty}\in h^NS_{\frac{1}{2},\beta}(\langle x \rangle^{-\infty})$, for $N$ as large as we want.
In this way we can factor out $\theta_h(x)$, write the equality
\begin{equation}
\begin{split}
& Op_h^w(x\xi + p(\xi))Op_h^w\left(\theta_h(x)\chi(h^{\beta}\xi)\gamma(\frac{x+p'(\xi)}{\sqrt{h}})\right)v^{\Sigma} = \\
& = Op_h^w(x\xi +p(\xi))\theta_h(x) Op_h^w\left(\chi(h^{\beta}\xi)\gamma(\frac{x+p'(\xi)}{\sqrt{h}})\right)v^{\Sigma} + hR(v) \, ,
\end{split}
\end{equation}
and return to $v^{\Sigma}_{\Lambda}$ by 
\begin{equation}
Op_h^w\Big(\chi(h^{\beta}\xi)\gamma(\frac{x+p'(\xi)}{\sqrt{h}})\Big)v^{\Sigma} = v^{\Sigma}_{\Lambda} -  Op_h^w\Big((1-\chi(h^{\beta}\xi))\gamma(\frac{x+p'(\xi)}{\sqrt{h}})\Big)v^{\Sigma}\,.
\end{equation}
Then,
\begin{equation}
\begin{split}
& Op_h^w(x\xi +p(\xi))\theta_h(x) Op_h^w\Big(\chi(h^{\beta}\xi)\gamma(\frac{x+p'(\xi)}{\sqrt{h}})\Big)v^{\Sigma} = \\
& = Op_h^w(x\xi +p(\xi))\theta_h(x) v_{\Lambda}^{\Sigma} - Op_h^w(x\xi +p(\xi))\theta_h(x) Op_h^w\Big((1-\chi(h^{\beta}\xi))\gamma(\frac{x+p'(\xi)}{\sqrt{h}})\Big)v^{\Sigma}\,,
\end{split}
\end{equation}
and one can show that the second term in the right hand side is a remainder $hR(v)$ essentially using symbolic calculus, lemma \ref{new estimate 1-chi}, and Sobolev injection.
Symbolic calculus enables us also to put $\theta_h(x)$ in $Op_h^w(x\xi + p(\xi))$, as the following deduction shows,
\begin{equation}
\begin{split}
Op_h^w(x\xi + p(\xi))\theta_h(x) v^{\Sigma}_{\Lambda}& = Op_h^w\big((x\xi + p(\xi))\theta_h(x) \big) v^{\Sigma}_{\Lambda} + \frac{h}{2i}Op_h^w\big(\theta_h'(x)(x + p'(\xi))\big)v^{\Sigma}_{\Lambda} + hR(v) \\
& = Op_h^w\big((x\xi + p(\xi))\theta_h(x) \big) v^{\Sigma}_{\Lambda} + hR(v)\, ,
\end{split}
\end{equation}
with $R(v)$ satisfying \eqref{estimates L^2 to L^2 of R(v)}, \eqref{estimates L^2 to L^inf of R(v)}, using proposition \ref{Continuity from $L^2$ to $L^2$} and Sobolev injection.
In the last equality, $\frac{h}{2i}Op_h^w\big(\theta_h'(x)(x + p'(\xi))\big)v^{\Sigma}_{\Lambda}$ enters in the remainder, for $\gamma(\frac{x+p'(\xi)}{\sqrt{h}}) \in S_{\frac{1}{2},0}(\langle \frac{x+p'(\xi)}{\sqrt{h}}\rangle^{-2})$ by lemma \ref{Class of Gamma}, $\theta_h'(x)(x + p'(\xi)) \in h^{-\delta}S_{\delta,0}(1)$ for a small $\delta>0$, and using symbolic calculus. 
Actually, we first write
\begin{equation}
\frac{h}{2i}Op_h^w\big(\theta_h'(x)(x + p'(\xi))\big)v^{\Sigma}_{\Lambda} = \frac{h^{\frac{3}{2}}}{2i}Op_h^w\Big(\theta'_h(x)\gamma(\frac{x+p'(\xi)}{\sqrt{h}})(\frac{x+p'(\xi)}{\sqrt{h}})\Big)v^{\Sigma} + h^{\frac{3}{2}} Op_h^w(r_0)v^{\Sigma} \,,
\end{equation}
where $r_0\in h^{-2\delta}S_{\frac{1}{2},0}(\langle \frac{x+p'(\xi)}{\sqrt{h}}\rangle^{-1})$, and then we use lemma \ref{lemma GammaTilde} with $\tilde{\Gamma}(\xi)=\gamma(\xi)\xi$, and lemma \ref{composition gamma-1 and its argument} to deduce that it is a remainder $hR(v)$.

\vspace{0.5cm}
\noindent We can now analyse $Op_h^w((x\xi + p(\xi))\theta_h(x))v^{\Sigma}_{\Lambda}$.
As we are away from points $x= \pm 1$, we can develop the symbol $x\xi + p(\xi)$ at $\xi = d\varphi(x)$, and since $\partial_{\xi}(x\xi +p(\xi))|_{\xi = d\varphi(x)}=0$ we have
\begin{equation} 
\begin{split}
x\xi + p(\xi) & = x d\varphi(x) + p(d\varphi(x)) + \int_0^1p''(t\xi +(1-t)d\varphi(x))(1-t) dt \, (\xi -d\varphi(x))^2 \\
& =  x d\varphi(x) + p(d\varphi(x)) + b(x,\xi) (x+p'(\xi))^2 \, ,
\end{split}
\end{equation}
where 
\begin{equation*}
b(x,\xi)= \int_0^1p''(t\xi +(1-t)d\varphi(x))(1-t) dt\,  \left(\frac{\xi -d\varphi(x)}{x+p'(\xi)}\right)^2  \, . 
\end{equation*}
Observe that  $x d\varphi(x) + p(d\varphi(x))= \varphi(x)$.
To conclude the proof, we need to show that $$Op_h^w\left(b(x,\xi)\theta_h(x)(x+p'(\xi))^2\right)v^{\Sigma}_{\Lambda}$$ gives rise to a remainder, too.
First, we may consider a function $\chi$ as in lemma \ref{new estimate 1-chi}, $\beta>0$, and split $b(x,\xi)\theta_h(x)(x+p'(\xi))^2$ as the sum of $b(x,\xi)\theta_h(x)(x+p'(\xi))^2(1-\chi(h^{\beta}\xi))\in h^{-\sigma}S_{\delta,0}(\langle\xi\rangle^2)$, for small $\delta,\sigma>0$, whose quantization acts on $v^{\Sigma}_{\Lambda}$ as a remainder because of lemma \ref{new estimate 1-chi}, and $b(x,\xi)\theta_h(x)(x+p'(\xi))^2\chi(h^{\beta}\xi)$.
For $b(x,\xi)\theta_h(x)\chi(h^{\beta}\xi)(x+p'(\xi))^2$, we can perform a further splitting through a function $\tilde{\gamma}\in C^{\infty}_0(\mathbb{R})$, such that $\tilde{\gamma}\left(\langle\xi\rangle^2(x+p'(\xi)) \right)\equiv 1$ on the support of $\chi(h^{\beta}\xi)\gamma(\frac{x+p'(\xi)}{\sqrt{h}})$, i.e.
\begin{equation}
\begin{split}
& b(x,\xi)\theta_h(x)\chi(h^{\beta}\xi)(x+p'(\xi))^2 \tilde{\gamma}\left(\langle\xi\rangle^2(x+p'(\xi)) \right) \\
& + b(x,\xi)\theta_h(x)\chi(h^{\beta}\xi)(x+p'(\xi))^2 (1-\tilde{\gamma})\left(\langle\xi\rangle^2(x+p'(\xi)) \right)\,.
\end{split}
\end{equation}
As discussed before, this implies that $(1-\tilde{\gamma})(\langle\xi\rangle^2(x+p'(\xi)))$ and all its derivatives are equal to zero on that support.
Since $b(x,\xi)\theta_h(x)\chi(h^{\beta}\xi)(x+p'(\xi))^2(1-\tilde{\gamma})\left(\langle\xi\rangle^2(x+p'(\xi))\right) \in h^{-\sigma}S_{\delta,\beta}(1)$ for $\sigma, \delta>0$ small depending on $\beta$, one can apply symbolic calculus (up to a large enough order) to obtain
\begin{equation}
b(x,\xi)\theta_h(x)\chi(h^{\beta}\xi)(x+p'(\xi))^2(1-\tilde{\gamma})\left(\langle\xi\rangle^2(x+p'(\xi))\right) \sharp \gamma(\frac{x+p'(\xi)}{\sqrt{h}}) = r'_{\infty} \,,
\end{equation}
with $r'_{\infty}=h^NS_{\frac{1}{2},\beta}(1)$, $N$ sufficiently large, to have $$Op_h^w\left( b(x,\xi)\theta_h(x)\chi(h^{\beta}\xi)(x+p'(\xi))^2(1-\tilde{\gamma})\left(\langle\xi\rangle^2(x+p'(\xi))\right)\right)v^{\Sigma}_{\Lambda} = hR(v) \,.$$
On the other hand, $B(x,\xi):= b(x,\xi)\theta_h(x)\chi(h^{\beta}\xi)\tilde{\gamma}(\langle \xi \rangle^2(x+p'(\xi)))$ belongs to $h^{-\sigma}S_{\delta,\beta}(1)$, for $\delta\ge 2\beta$, by lemma \ref{lem on e and etilde}.
Using twice lemma \ref{lem of composition with Gamma Tilde}, together with the fact that $\gamma(\frac{x+p'(\xi)}{\sqrt{h}})\in S_{\frac{1}{2},0}(\langle \frac{x+p'(\xi)}{\sqrt{h}} \rangle^{-3})$ and $B(x,\xi)(x+p'(\xi))^2 \in h^{1-\sigma}S_{\delta,\beta}(\langle\frac{x+p'(\xi)}{\sqrt{h}}\rangle^2)$, we derive
\begin{equation}
\left(B(x,\xi)(x+p'(\xi))^2 \right)\sharp \gamma(\frac{x+p'(\xi)}{\sqrt{h}}) = B(x,\xi) \gamma(\frac{x+p'(\xi)}{\sqrt{h}})(x+p'(\xi))^2 + hr_0 \,, 
\end{equation}
and
\begin{equation}
\left(B(x,\xi)\gamma(\frac{x+p'(\xi)}{\sqrt{h}})(x+p'(\xi)) \right)\sharp (x+p'(\xi)) = B(x,\xi)\gamma(\frac{x+p'(\xi)}{\sqrt{h}})(x+p'(\xi))^2 + hr'_0 \, , 
\end{equation}
where $r_0, r'_0 \in h^{\frac{1}{2}-\sigma} S_{\frac{1}{2},\beta}(\langle \frac{x+p'(\xi)}{\sqrt{h}}\rangle^{-1})$.
Therefore
\begin{equation}
\left(B(x,\xi) (x+p'(\xi))^2\right) \sharp \gamma(\frac{x+p'(\xi)}{\sqrt{h}}) = \left(B(x,\xi)\gamma(\frac{x+p'(\xi)}{\sqrt{h}})(x+p'(\xi)) \right)\sharp (x+p'(\xi)) + h(r_0-r'_0) \,,
\end{equation}
and 
\begin{equation}
Op_h^w(B(x,\xi)(x+p'(\xi))^2)v_{\Lambda}^{\Sigma} = h Op_h^w\Big(B(x,\xi)\gamma(\frac{x+p'(\xi)}{\sqrt{h}})(x+p'(\xi))\Big)\mathcal{L}v^{\Sigma} + hOp_h^w(r_0-r'_0)v^{\Sigma}\,,
\end{equation}
so one can show that the right hand side is a remainder $hR(v)$, commutating $\mathcal{L}$ with $Op_h^w(\Sigma(\xi))$, using that $B(x,\xi)\gamma(\frac{x+p'(\xi)}{\sqrt{h}})(x+p'(\xi)), r_0-r'_0 \in h^{\frac{1}{2}-\sigma}S_{\frac{1}{2},\beta}(\langle \frac{x+p'(\xi)}{\sqrt{h}}\rangle^{-1})$ , and propositions \ref{Continuity from $L^2$ to $L^2$}, \ref{Continuity from $L^2$ to L^inf}.
We finally obtain 
\begin{equation}
Op_h^w(x\xi + p(\xi))v^{\Sigma}_{\Lambda}= \varphi(x) \theta_h(x) v^{\Sigma}_{\Lambda} + h R(v) \,, 
\end{equation}
and according to \eqref{ODE with the last pseudo op}, $v^{\Sigma}_{\Lambda}$ is solution of 
\begin{equation}
\begin{split}
D_t v^{\Sigma}_{\Lambda} & =\varphi(x) \theta_h(x) v^{\Sigma}_{\Lambda} + h \Phi_1^{\Sigma}(x)\theta_h(x)|v^{\Sigma}_{\Lambda}|^2v^{\Sigma}_{\Lambda} \\
& + h Op_h^w(\Gamma)\left[\Phi_3^{\Sigma}(x)\theta_h(x)(v^{\Sigma}_{\Lambda})^3 + \Phi_{-1}^{\Sigma}(x)\theta_h(x) |v^{\Sigma}_{\Lambda}|^2\overline{v^{\Sigma}_{\Lambda}} + \Phi_{-3}^{\Sigma}(x)\theta_h(x) (\overline{v^{\Sigma}_{\Lambda}})^3\right] \\
& + hR(v)\,,
\end{split}
\end{equation}
where $R(v)$ is a remainder satisfying estimates \eqref{estimates L^2 to L^2 of R(v)}, \eqref{estimates L^2 to L^inf of R(v)}.
\endproof
\end{prop}

\section{Study of the ODE and End of the Proof} \label{Study of the ODE and End of the Proof}

\subsection{The $L^{\infty}$ estimate}
The goal of this subsection is to the derive from the equation \eqref{ODE for vLambda} an ODE for a new function $f^{\Sigma}_{\Lambda}$ obtained from $v^{\Sigma}_{\Lambda}$, from which we can deduce uniform bounds for $v^{\Sigma}_{\Lambda}$, and for the starting function $v$, with a certain number $\rho\in \mathbb{N}$ of its derivatives.
The idea is to get rid of contributions of \emph{non characteristic} terms (i.e. of  cubic terms different from $|v^{\Sigma}_{\Lambda}|^2v^{\Sigma}_{\Lambda}$) by a reasoning of normal forms.
This will allow us to eliminate all terms still containing pseudo-differential operators, to finally write an ODE, and to prove the required $L^{\infty}$ estimate, if the \emph{null condition} is satisfied.

\vspace{0.5cm}
\noindent In the previous section, we denoted by $\Phi^{\Sigma}_3(x)$, $\Phi^{\Sigma}_1(x)$, $\Phi^{\Sigma}_{-1}(x)$ and $\Phi^{\Sigma}_{-3}(x)$
respectively the coefficients of $(v^{\Sigma}_{\Lambda})^3$, $|v^{\Sigma}_{\Lambda}|^2v^{\Sigma}_{\Lambda}$, $|v^{\Sigma}_{\Lambda}|^2\overline{v^{\Sigma}_{\Lambda}}$, $(\overline{v^{\Sigma}_{\Lambda}})^3$ in the right hand side of \eqref{ODE for vLambda}.
One can calculate them explicitly, using both the expression of the nonlinearity obtained in proposition \ref{development of nonlinearity} and its polynomial representation as in equation \eqref{half KG for v with polynomials}.
In the latter, after the development at $\xi = d\varphi(x)$, we essentially replaced $hD$ by $d\varphi(x)$ when it is applied to $v^{\Sigma}_{\Lambda}$, and by $-d\varphi(x)$ when it is applied to $\overline{v^{\Sigma}_{\Lambda}}$, modulus some new smooth coefficients $a_I(x):=A(\displaystyle\sum_{l=1}^3i_ld\varphi(x))\Sigma(d\varphi(x))^{-3}$, for every $I=(i_1, i_2, i_3)$ (the factor $\Sigma(d\varphi(x))^{-3}$ coming out from $m^{\Sigma}_I(d\varphi_I(x))= m_I(d\varphi(x))\Sigma(d\varphi(x))^{-3}$, according to the notation introduced in \eqref{def of m_I Sigma}, $A(\xi)= \Sigma(\xi)\chi(h^{\beta}\xi)$).

\vspace{0.5cm}
\noindent We are interested in particular in $\Phi^{\Sigma}_1(x)$ or, to be more precise, to its real part. 
In fact, the \emph{null condition} introduced in definition \ref{definition of null condition} at the very beginning is the same as requiring for the coefficient of $|v^{\Sigma}_{\Lambda}|^2v^{\Sigma}_{\Lambda}$ to be real, i.e. its imaginary part must be equal to zero. Since polynomials $P'_k$, $P''_k$ are real as well as $d\varphi(x), \langle d\varphi(x)\rangle$, the only contribution to the imaginary part comes from $P'_k$, $P''_k$ for $k=1, 3$ (which have a factor $i^k$) and produces a multiple of  the function $\Phi(x)$ defined in \eqref{Phi null condition}.
Therefore, if we suppose that the nonlinearity satisfies this \emph{null condition} (as demanded in theorem \ref{main theorem}) then we find for $\Phi^{\Sigma}_1(x)$ that
\begin{equation}
\begin{split}
\Phi^{\Sigma}_1(x) = \frac{1}{8}a_{(1,1,-1)}(x)\langle d\varphi \rangle^{-3} & \left[3 P_0(1, d\varphi \langle d\varphi\rangle, (d\varphi)^2; \, \langle d\varphi\rangle, d\varphi) \right.\\
& \left. + P_2(1, d\varphi \langle d\varphi\rangle, (d\varphi)^2; \, \langle d\varphi\rangle, d\varphi) \right].
\end{split}
\end{equation}

\begin{prop} \label{prop ODE for fLambda}
Suppose we are given two constants $A'', B''>0$, some $T>1$ and a $\sigma>0$ small.
Let $v^{\Sigma}_{\Lambda}$ be a solution of the equation \eqref{ODE for vLambda} on the interval $[1,T]$, $v^{\Sigma}_{\Lambda}$ satisfying the a priori estimates
\begin{align} 
\|v^{\Sigma}_{\Lambda}(t,\cdot)\|_{L^{\infty}(\mathbb{R})} & \le A'' \varepsilon  \, ,  \label{a priori estimate on vLambda} \\
\|v^{\Sigma}_{\Lambda}(t,\cdot)\|_{L^{2}(\mathbb{R})} & \le B'' \varepsilon h^{-\sigma} \, , \label{a priori L2 estimate on vLambda} 
\end{align}
for all $t \in [1,T]$.
Let $\tilde{\theta}_h(x) \in C^{\infty}_0(]-1,1[)$, such that $\tilde{\theta}_h \theta_h \equiv \theta_h$, and define
\begin{equation} \label{definition of fLambda}
f^{\Sigma}_{\Lambda}:= v^{\Sigma}_{\Lambda} + Op_h^w(\Gamma)\left[- \frac{h}{2}\frac{\tilde{\theta}_h(x)}{\varphi(x)}\Phi^{\Sigma}_3(x) (v^{\Sigma}_{\Lambda})^3 + \frac{h}{2}\frac{\tilde{\theta}_h(x)}{\varphi(x)}\Phi^{\Sigma}_{-1}(x) |v^{\Sigma}_{\Lambda}|^2\overline{v^{\Sigma}_{\Lambda}} + \frac{h}{4}\frac{\tilde{\theta}_h(x)}{\varphi(x)}\Phi^{\Sigma}_{-3}(x) (\overline{v^{\Sigma}_{\Lambda}})^3 \right]\, .
\end{equation}
Then $f^{\Sigma}_{\Lambda}$ is well defined and it is solution of the ODE:
\begin{equation} \label{ODE for fLambda}
D_t f^{\Sigma}_{\Lambda}= \varphi(x)\theta_h(x) f^{\Sigma}_{\Lambda} + h\theta_h(x) \Phi^{\Sigma}_1(x) |f^{\Sigma}_{\Lambda}|^2 f^{\Sigma}_{\Lambda} + hR(v) \, ,
\end{equation}
where $R(v)$ is a remainder satisfying estimates \eqref{estimates L^2 to L^2 of R(v)}, \eqref{estimates L^2 to L^inf of R(v)}.
\proof
Firstly, we would like to underline that, if we suppose bounds in \eqref{a priori estimates on v} and \eqref{a priori estimate on Lv and v} on $v$, then hypothesis \eqref{a priori estimate on vLambda} and \eqref{a priori L2 estimate on vLambda} follow immediately, because of the definition of $v^{\Sigma}_{\Lambda}$ as $Op_h^w(\Gamma)v^{\Sigma}$.
In fact, estimate \eqref{a priori L2 estimate on vLambda} follows from proposition \ref{Continuity from $L^2$ to $L^2$} and the \emph{a priori} estimate \eqref{a priori estimate on Lv and v}, with $B''=B'$. Regarding the estimate \eqref{a priori estimate on vLambda}, we can write
\begin{equation}
v^{\Sigma}_{\Lambda} = v^{\Sigma} - v^{\Sigma}_{\Lambda^c} \, ,
\end{equation}
and since $\|v^{\Sigma}(t,\cdot)\|_{L^{\infty}}= \|v(t,\cdot)\|_{W^{\rho,\infty}_h}$,
\begin{equation}
\begin{split}
\|v^{\Sigma}_{\Lambda}(t,\cdot)\|_{L^{ \infty}} & \le \|v^{\Sigma}(t,\cdot)\|_{L^{\infty}} + \|v^{\Sigma}_{\Lambda^c}(t,\cdot)\|_{L^{\infty}} \\
& = \|v(t,\cdot)\|_{W_h^{\rho, \infty}} + \|v^{\Sigma}_{\Lambda^c}(t,\cdot)\|_{L^{\infty}} \, ,
\end{split}
\end{equation}
where we estimated $\|v^{\Sigma}_{\Lambda^c}(t,\cdot)\|_{L^{\infty}}$ in proposition \ref{Estimates on v Lambda complementary}. 
Therefore, using that for $\sigma>0$ sufficiently small $h^{\frac{1}{4}-\sigma}\le h^{\frac{1}{8}}$, we have
\begin{equation}
\begin{split}
\|v^{\Sigma}_{\Lambda}(t,\cdot)\|_{L^{\infty}} & \le \|v(t,\cdot)\|_{W_h^{\rho, \infty}} + C h^\frac{1}{8}( \|\mathcal{L}v(t,\cdot)\|_{L^2} + \|v(t,\cdot) \|_{H^s_h}) \\
& \le  A' \varepsilon + CB' \varepsilon \, h^{\frac{1}{8}-\sigma} \\
& \le A''\varepsilon \, ,
\end{split}
\end{equation}
if we choose $A''>0$ sufficiently large to have $A',CB' \le \frac{A''}{2}$.

\vspace{0.5cm}
\noindent Secondly, $\varphi(x)\ne 0$ for all $x$ in the support of $\tilde{\theta}_h$.
In fact, we consider $\tilde{\theta}_h$ such that $\tilde{\theta}_h\theta_h \equiv \theta_h$, so we can suppose that its support is of the form $[-1+C'h^{2\beta},1-C'h^{2\beta}]$, for a suitable small positive constant $C'$. On this interval $x^2 \le (1-C'h^{2\beta})^2= 1 + C'^2h^{4\beta} -2C'h^{2\beta}$, so
\begin{equation}
\varphi(x) = \sqrt{1-x^2} \ge \sqrt{C'h^{2\beta}(2-C'h^{2\beta})}\gtrsim h^{\beta} \,,
\end{equation}
which implies that the quotient $\frac{\tilde{\theta}_h(x)}{\varphi(x)}$ is well defined and $|\frac{\tilde{\theta}_h(x)}{\varphi(x)}|\le h^{-\beta}$.
Then, set 
\begin{equation}  \label{early def of fLambda}
f^{\Sigma}_{\Lambda}:= v^{\Sigma}_{\Lambda} + Op_h^w(\Gamma) \left[h \frac{\tilde{\theta}_h(x)}{\varphi(x)} \left(k_1 \Phi^{\Sigma}_3(x) (v^{\Sigma}_{\Lambda})^3 + k_2 \Phi^{\Sigma}_{-1}(x)|v^{\Sigma}_{\Lambda}|^2\overline{v^{\Sigma}_{\Lambda}} + k_3 \Phi^{\Sigma}_{-3}(x)(\overline{v^{\Sigma}_{\Lambda}})^3 \right)\right] \, ,
\end{equation}
with $k_1, k_2, k_3 \in \mathbb{R}$ to be properly chosen, and apply $D_t$ to this expression.
We have already calculated $D_tOp_h^w(\Gamma)$ in \eqref{Dt vLambda}, obtaining that the commutator is
\begin{equation}
[D_t,Op_h^w(\Gamma)]= i h^{\frac{1}{2}}Op_h^w\left(\gamma'(\frac{x+p'(\xi)}{\sqrt{h}})p''(\xi)\xi\right) - \frac{ih}{2}Op_h^w\left(\gamma'(\frac{x+p'(\xi)}{\sqrt{h}})(\frac{x+p'(\xi)}{\sqrt{h}})\right)\,,
\end{equation}
where both appearing symbols belong to $S_{\frac{1}{2},0}(\langle \frac{x+p'(\xi)}{\sqrt{h}}\rangle^{-1})$. 
The truncation of these symbols through a function $\chi(h^{\beta}\xi)$ as in lemma \ref{new estimate 1-chi}, and propositions \ref{Continuity from $L^2$ to $L^2$}, \ref{Continuity from $L^2$ to L^inf}, together with estimates \eqref{a priori estimate on vLambda}, \eqref{a priori L2 estimate on vLambda} on $v^{\Sigma}_{\Lambda}$, show that the action of the commutator on brackets in \eqref{early def of fLambda} gives rise to a remainder $hR(v)$. \\
Denoting by $O(5)$ all terms of order 5 in $(v^{\Sigma}_{\Lambda}, \overline{v^{\Sigma}_{\Lambda}})$, and using \eqref{ODE for vLambda}, we can compute
\begin{equation}
\begin{split}
D_t f^{\Sigma}_{\Lambda} = D_t v^{\Sigma}_{\Lambda} + Op_h^w(\Gamma) & \left[ k_1 h \frac{\tilde{\theta}_h(x)}{\varphi(x)} \Phi^{\Sigma}_3(x) [ 3\varphi(x) \theta_h(x) (v^{\Sigma}_{\Lambda})^3 + h^2 O(5)] \right.\\
& \left. + k_2 h \frac{\tilde{\theta}_h(x)}{\varphi(x)} \Phi^{\Sigma}_{-1}(x) [-\varphi(x) \theta_h(x) |v^{\Sigma}_{\Lambda}|^2\overline{v^{\Sigma}_{\Lambda}} + h^2O(5)]\right. \\
& \left. + k_3 h \frac{\tilde{\theta}_h(x)}{\varphi(x)} \Phi^{\Sigma}_{-3}(x) [-3\varphi(x) \theta_h(x) (\overline{v^{\Sigma}_{\Lambda}})^3 + h^2 O(5)]\right] + hR(v) \, ,
\end{split}
\end{equation}
where $hR(v)$ includes also terms coming out from $D_t(h \tilde{\theta}_h(x))$, and
\begin{equation} \label{step intermediario in ODE}
\begin{split}
D_tf^{\Sigma}_{\Lambda}& = \varphi(x) \theta_h(x) v^{\Sigma}_{\Lambda} +  h\theta_h(x) \Phi^{\Sigma}_1(x) |v^{\Sigma}_{\Lambda}|^2v^{\Sigma}_{\Lambda}\\
& + Op_h^w(\Gamma)\left[h \theta_h(x)\left( (3k_1 + 1) \Phi^{\Sigma}_3(x) (v^{\Sigma}_{\Lambda})^3
+ (-k_2+1) \Phi^{\Sigma}_{-1}(x) |v^{\Sigma}_{\Lambda}|^2\overline{v^{\Sigma}_{\Lambda}} \right.\right.\\
& \left.\left.\hspace{1.2cm}+ (-3k_3 + 1) \Phi^{\Sigma}_{-3}(x)(\overline{v^{\Sigma}_{\Lambda}})^3\right) \right] + hR(v) \, ,
\end{split}
\end{equation}
where $h^2 O(5)$ entered in $hR(v)$ from propositions \ref{Continuity from $L^2$ to $L^2$}, \ref{Continuity from $L^2$ to L^inf}, estimates \eqref{a priori estimate on vLambda}, \eqref{a priori L2 estimate on vLambda}, and the fact that involved coefficients are $O(h^{-\sigma})$, for a small $\sigma>0$.
We use again the definition of $f^{\Sigma}_{\Lambda}$ to replace $v^{\Sigma}_{\Lambda}$ in the linear and in the \emph{characteristic} part.
We have $h\theta_h(x)\Phi_1^{\Sigma}(x)|v^{\Sigma}_{\Lambda}|^2v^{\Sigma}_{\Lambda}=h\theta_h(x)\Phi_1^{\Sigma}(x)|f^{\Sigma}_{\Lambda}|^2f^{\Sigma}_{\Lambda} + h^2O(5)$ and
\begin{equation}
\begin{split}
\varphi(x)\theta_h(x)v^{\Sigma}_{\Lambda} &= \varphi(x)\theta_h(x)f^{\Sigma}_{\Lambda} - \varphi(x)\theta_h(x)Op_h^w(\Gamma)\left[h \frac{\tilde{\theta}_h(x)}{\varphi(x)} \left(k_1 \Phi^{\Sigma}_3(x) (v^{\Sigma}_{\Lambda})^3 + k_2 \Phi^{\Sigma}_{-1}(x)|v^{\Sigma}_{\Lambda}|^2\overline{v^{\Sigma}_{\Lambda}}  \right.\right. \\
& \left. \left. + k_3 \Phi^{\Sigma}_{-3}(x)(\overline{v^{\Sigma}_{\Lambda}})^3 \right)\right] \\
& = \varphi(x)\theta_h(x)f^{\Sigma}_{\Lambda}  - Op_h^w(\Gamma)\left[h \theta_h(x)\left(k_1 \Phi^{\Sigma}_3(x) (v^{\Sigma}_{\Lambda})^3+ k_2 \Phi^{\Sigma}_{-1}(x)|v^{\Sigma}_{\Lambda}|^2\overline{v^{\Sigma}_{\Lambda}} \right.\right. \\
& \left. \left. + k_3 \Phi^{\Sigma}_{-3}(x)(\overline{v^{\Sigma}_{\Lambda}})^3 \right)\right] + hR(v) \,,
\end{split}
\end{equation}
where the last equality is consequence of the fact that, by lemma \ref{lem of composition with Gamma Tilde}, $[\varphi(x)\theta_h(x), Op_h^w(\Gamma)]= h^{\frac{1}{2}-\sigma}Op_h^w(r_0)$, $r_0 \in S_{\frac{1}{2},0}(\langle \frac{x+p'(\xi)}{\sqrt{h}}\rangle^{-1})$, $\sigma>0$ small. 
Again a truncation through $\chi(h^{\beta}\xi)$, and the application of propositions \ref{Continuity from $L^2$ to $L^2$}, \ref{Continuity from $L^2$ to L^inf}, together with estimates on $v^{\Sigma}_{\Lambda}$, ensure that the contribution coming from the action of the commutator on its argument enters in the remainder.
We finally obtain
\begin{equation}
\begin{split}
D_t f^{\Sigma}_{\Lambda}& = \varphi(x) \theta_h(x) f^{\Sigma}_{\Lambda} + h \theta_h(x) \Phi^{\Sigma}_1(x) |f^{\Sigma}_{\Lambda}|^2f^{\Sigma}_{\Lambda} \\
& +Op_h^w(\Gamma)\left[ h \theta_h(x)\left( (2k_1 + 1) \Phi^{\Sigma}_3(x) (v^{\Sigma}_{\Lambda})^3 + (-2k_2 +1) \Phi^{\Sigma}_{-1}(x) |v^{\Sigma}_{\Lambda}|^2\overline{v^{\Sigma}_{\Lambda}} \right.\right.\\
&  \left.\left. \hspace{1.7cm} +(-4k_3 +1) \Phi^{\Sigma}_{-3}(x) (\overline{v^{\Sigma}_{\Lambda}})^3 \right)\right] + hR(v) \,,
\end{split}
\end{equation}
and we get rid of \emph{non-characteristic} terms by requiring
\begin{equation*}
\begin{cases}
2k_1 + 1 &=0\\
-2k_2 +1 &=0 \\
-4k_3 + 1& =0 
\end{cases}\, \qquad \Rightarrow \qquad
\begin{cases}
& k_1 = -\frac{1}{2}\\
& k_2 = \frac{1}{2}\\
& k_3 = \frac{1}{4}\, ,
\end{cases}\, 
\end{equation*}
from which the statement.
\endproof
\end{prop}

\begin{prop}\label{Linfty estimate of fLambda}
Let $f^{\Sigma}_{\Lambda}$ be the function defined in \eqref{definition of fLambda}, solution of the ODE \eqref{ODE for fLambda} under the a priori estimates \eqref{a priori estimate on vLambda}, \eqref{a priori L2 estimate on vLambda}.
Then the following inequality holds :
\begin{equation}
\|f^{\Sigma}_{\Lambda}(t,\cdot)\|_{L^{\infty}} \le \|f^{\Sigma}_{\Lambda}(1,\cdot)\|_{L^{\infty}} + C \int_1^t \tau^{-\frac{5}{4}+\sigma} (\|\mathcal{L}v(\tau,\cdot)\|_{L^2} + \|v(\tau,\cdot)\|_{H^s_h})\, d\tau \, ,
\end{equation}
for $\sigma>0$ small, and a positive constant $C>0$.
\proof
Using the equation \eqref{ODE for fLambda}, we can compute
\begin{equation} \label{control on derivative of flambda2}
\begin{split}
\frac{\partial}{\partial t}|f^{\Sigma}_{\Lambda}(t,x)|^2 & = 2 \Im(f^{\Sigma}_{\Lambda}\overline{D_tf^{\Sigma}_{\Lambda}})(t,x)= 2\Im(\varphi(x)\theta_h(x)|f^{\Sigma}_{\Lambda}|^2 + h\theta_h(x)\Phi^{\Sigma}_1(x) |f^{\Sigma}_{\Lambda}|^4 + h R(v) f^{\Sigma}_{\Lambda})(t,x)\\
& = 2\Im(hR(v) f^{\Sigma}_{\Lambda})(t,x) \le 2h |f^{\Sigma}_{\Lambda}(t,x)| |R(v)| \,,
\end{split}
\end{equation}
from which follows an integral inequality
\begin{equation}
\|f^{\Sigma}_{\Lambda}(t,\cdot)\|_{L^{\infty}} \le \|f^{\Sigma}_{\Lambda}(1,\cdot)\|_{L^{\infty}} + \int_1^t \frac{\|R(v)(\tau,\cdot)\|_{L^{\infty}}}{\tau} \, d\tau \,.
\end{equation}
Using the estimate \eqref{estimates L^2 to L^inf of R(v)} for $R(v)$, we obtain the result
\begin{equation}  \label{integral inequality for fLambda}
\|f^{\Sigma}_{\Lambda}(t,\cdot)\|_{L^{\infty}}
 \le \|f^{\Sigma}_{\Lambda}(1,\cdot)\|_{L^{\infty}} + C \int_1^t \tau^{-\frac{5}{4}+\sigma}(\|\mathcal{L}v (\tau,\cdot)\|_{L^2}) + \|v(\tau,\cdot)\|_{H^s_h}) \, d\tau \,.
\end{equation} 
\endproof
\end{prop}

\noindent Finally, the $L^{\infty}$ estimate we found for $f^{\Sigma}_{\Lambda}$ in the previous proposition enables us to propagate the uniform estimate on $v$, as showed in the following:

\begin{prop}[Propagation of the uniform estimate] \label{Propagation of the uniform estimate}
Let $v$ be a solution of the equation \eqref{half KG for v with polynomials} on some interval $[1,T]$, $T>1$ and $\sigma >0$ small.
Then, for a fixed constant $K>1$, there exist two constants $A', B'>0$ sufficiently large, $\varepsilon_0>0$ sufficiently small, $s, \rho \in \mathbb{N}$ with $s\gg \rho$, such that, if $0<\varepsilon<\varepsilon_0$, and $v$ satisfies
\begin{equation} \label{a priori estimates in propagation of the unif est}
\begin{split}
(A.1) &\, \|v(t,\cdot)\|_{W_h^{\rho, \infty}} \le A' \varepsilon \,, \\
(B.1) & \,\|v(t,\cdot)\|_{H^s_h} \le B'\varepsilon \, h^{-\sigma} \, , \\
(B.2) & \,\|\mathcal{L}v(t,\cdot)\|_{L^2} \le B'\varepsilon \, h^{-\sigma} \, ,
\end{split}
\end{equation}
for every $t \in [1,T]$, then it satisfies also
\begin{equation}
(A.1') \,\|v(t,\cdot)\|_{W_h^{\rho, \infty}} \le \frac{A'}{K}\varepsilon \, , \qquad \forall t\in [1,T]\, .
\end{equation}
\proof
The proof of the proposition comes directly from proposition \ref{Linfty estimate of fLambda} and from the equivalence between $\|v^{\Sigma}_{\Lambda}\|_{L^{\infty}}$ and $\|f^{\Sigma}_{\Lambda}\|_{L^{\infty}}$.
In fact, functions $\Phi^{\Sigma}_j(x)$ are cubic expressions in $d\varphi(x)$ and $\langle d\varphi(x)\rangle$, so they are bounded up to a loss $h^{-\delta}$, $\delta>0$ depending on $\beta$, on the support of $\tilde{\theta}_h(x)$, where also $\varphi(x) \gtrsim  h^{\beta}>0$.
This implies that $|\frac{\tilde{\theta}_h(x)}{\varphi(x)}\Phi^{\Sigma}_j(x)|\le C h^{-\delta}$, $j\in \{3, -1, -3\}$, with  a new $\delta>0$ depending linearly on $\beta$, so that by the definition of $f^{\Sigma}_{\Lambda}$, proposition \ref{Continuity from $L^2$ to L^inf} and estimates \eqref{a priori estimate on vLambda}, \eqref{a priori L2 estimate on vLambda} (which follow from \eqref{a priori estimates in propagation of the unif est}, as already observed in proposition \ref{prop ODE for fLambda}), we find
\begin{equation}
\frac{1}{2} \|v^{\Sigma}_{\Lambda}(t,\cdot)\|_{L^{\infty}}\le \|f^{\Sigma}_{\Lambda}(t,\cdot)\|_{L^{\infty}} \le 2 \|v^{\Sigma}_{\Lambda}(t,\cdot)\|_{L^{\infty}} \, .
\end{equation}
Furthermore, the \emph{a priori} estimate on the $W_h^{\rho, \infty}$ norm of $v$ extends to the $L^{\infty}$ norm of $v^{\Sigma}_{\Lambda}$ just by the decomposition
\begin{equation}
v^{\Sigma}_{\Lambda}= v^{\Sigma} - v^{\Sigma}_{\Lambda^c} \, , 
\end{equation}
and by proposition \ref{Estimates on v Lambda complementary}, so for example at time $t=1$ we have
\begin{equation}
\begin{split}
\|v^{\Sigma}_{\Lambda}(1,\cdot)\|_{L^{\infty}} & \le \|v^{\Sigma}(1,\cdot)\|_{L^{\infty}} + \|v^{\Sigma}_{\Lambda^c}(1,\cdot)\|_{L^{\infty}} \\
& \le \|v(1,\cdot)\|_{W_h^{\rho, \infty}} + C (\|\mathcal{L}v(1,\cdot)\|_{L^2} + \|v(1,\cdot)\|_{H_h^s}) \\
& \le \frac{A'}{32 K}\varepsilon + CB'\varepsilon \\
& \le \frac{A'}{16 K}\varepsilon \, ,
\end{split}
\end{equation}
where we choose $A'>0$ sufficiently large such that $\|v(1,\cdot)\|_{W_h^{\rho, \infty}}\le \frac{A'}{32 K}\varepsilon$ and $CB' < \frac{A'}{32 K}$.
Therefore
\begin{equation} \label{est fSigma Lambda}
\|f^{\Sigma}_{\Lambda}(1,\cdot)\|_{L^{\infty}} \le 2 \|v^{\Sigma}_{\Lambda}(1,\cdot)\|_{L^{\infty}} \le \frac{A'}{8K}\varepsilon \, .
\end{equation}
Using proposition \ref{Linfty estimate of fLambda}, \eqref{est fSigma Lambda}  and the \emph{a priori} estimates (B.1), (B.2), we find that
\begin{equation}
\begin{split}
\|f^{\Sigma}_{\Lambda}(t,\cdot)\|_{L^{\infty}} & \le \frac{A'}{8 K}\varepsilon + CB'\varepsilon\int_1^t \tau^{-\frac{5}{4}+\sigma} \, d\tau\\
& \le \frac{A'}{8K}\varepsilon + C'B'\varepsilon \\
& \le \frac{A'}{4K}\varepsilon \, ,
\end{split}
\end{equation}
where again the last inequality follows from the choice of $A'>0$ large enough to have $C'B'< \frac{A'}{8K}$.
Then we have 
\begin{equation}
\|v^{\Sigma}_{\Lambda}(t,\cdot)\|_{L^{\infty}} \le \frac{A'}{2K}\varepsilon \, ,
\end{equation}
and 
\begin{equation}
\begin{split}
\|v^{\Sigma}(t,\cdot)\|_{L^{\infty}} & \le \|v^{\Sigma}_{\Lambda}(t,\cdot)\|_{L^{\infty}} + \|v^{\Sigma}_{\Lambda^c}(t,\cdot)\|_{L^{\infty}} \\
& \le \frac{A'}{2K}\varepsilon + CB'\varepsilon h^{\frac{1}{4}-\sigma'} \\
& \le \frac{A'}{K}\varepsilon \, .
\end{split}
\end{equation}
\endproof
\end{prop}

\subsection{Asymptotics}

We want now to derive the asymptotic expansion for the function $\langle hD\rangle^{-1}v$, $v$ being the solution of \eqref{half KG for v with polynomials}, when it exists on $[1,+\infty[$.
The reader can refer to the next subsection to find the proof of the global existence of $v$, which implies also the global existence of the solution $u$ of the starting problem \eqref{KG}.

\begin{prop}
Under the same hypothesis as theorem \ref{thm : From the PDE problem to an ODE},  with $T=+\infty$, there exists a family $(\theta_h(x))_h$ of $C^{\infty}$ functions, real valued, supported in some interval $[-1+ch^{2\beta}, 1-ch^{2\beta}]$, $\theta_h\equiv 1$ on an interval of the same form, such that $(h\partial_h)^k\theta_h(x)$ is bounded for any $k$, and a family $(a_{\varepsilon})_{\varepsilon \in ]0,\varepsilon_0]}$ of $\mathbb{C}$-valued functions on $\mathbb{R}$, supported in $[-1,1]$, uniformly bounded, such that
\begin{equation} \label{asymptotics for v}
\langle hD\rangle^{-1}v = \varepsilon a_{\varepsilon}(x) \exp\left[i\varphi(x) \int_1^t\theta_{1/\tau}(x)d\tau + i\varepsilon^2|a_{\varepsilon}(x)|^2\Phi^{\Sigma}_1(x)\int_1^t\theta_{1/\tau}(x)\frac{d\tau}{\tau} \right] + t^{-\frac{1}{4}+\sigma}r(t,x) \,,
\end{equation}
where $h=\frac{1}{t}$, $\sigma>0$ is small and $\displaystyle\sup_{t\ge 1}\|r(t,\cdot)\|_{L^2\cap L^{\infty}}\le C \varepsilon$.
\proof
Let us take $\Sigma(\xi)= \langle \xi \rangle^{-1}$, so that $v^{\Sigma}=\langle hD\rangle^{-1}v$.
Resuming all prevoius results, we have obtained that under the \emph{a priori} estimates \eqref{a priori estimates on v}, \eqref{a priori estimate on Lv and v}, the function $f^{\Sigma}_{\Lambda}$ defined in \eqref{definition of fLambda} satisfies \eqref{ODE for fLambda}, with a remainder $R(v) = O_{L^{\infty}\cap L^2}(\varepsilon t^{-\frac{1}{4}+\sigma})$, for a sufficiently small $\sigma>0$.
Inequality \eqref{control on derivative of flambda2} and the bound \eqref{estimates L^2 to L^inf of R(v)} show that
\begin{equation*}
\|f^{\Sigma}_{\Lambda}(t,\cdot)-f^{\Sigma}_{\Lambda}(t',\cdot)\|_{L^{\infty}} \le C \int_{t'}^t \tau^{-\frac{5}{4}+\sigma}(\|\mathcal{L}v(\tau,\cdot)\|_{L^2}+ \|v(\tau,\cdot)\|_{H^s_h})\, d\tau \,.
\end{equation*}
Combining with the \emph{a priori} estimate \eqref{a priori estimate on Lv and v}, there is a continuous function $x\rightarrow |\tilde{a}(x)|$ such that $\left||f^{\Sigma}_{\Lambda}(t,x)|^2-|\tilde{a}(x)|^2\right|=O(\varepsilon t^{-\frac{1}{2}+\sigma})$, for a new small $\sigma>0$, and replacing this new function in \eqref{ODE for fLambda} we obtain the equation
\begin{equation}
D_t f^{\Sigma}_{\Lambda}= \theta_h(x)\left[\varphi(x) + h \Phi^{\Sigma}_1(x) |\tilde{a}(x)|^2\right] f^{\Sigma}_{\Lambda} + h\, r(t,x) \,,
\end{equation}
for $r= O_{L^{\infty}\cap L^2}(\varepsilon t^{-\frac{1}{4}+\sigma})$, which is a linear non homogeneous ODE for $f^{\Sigma}_{\Lambda}$.
This implies that there is a $O(\varepsilon)$ continuous function $\tilde{a}$ such that
\begin{equation} \label{asymptotics for fLambda}
f^{\Sigma}_{\Lambda}(t,x)=\tilde{a}(x) \exp\left[i\varphi(x) \int_1^t\theta_{1/\tau}(x)d\tau + i|\tilde{a}(x)|^2\Phi^{\Sigma}_1(x)\int_1^t\theta_{1/\tau}(x)\frac{d\tau}{\tau} \right] + t^{-\frac{1}{4}+\sigma}r(t,x)\, ,
\end{equation}
for a new $r$. 
Finally, using the definition of $f^{\Sigma}_{\Lambda}$ and proposition \ref{Estimates on v Lambda complementary}, we have $\|f^{\Sigma}_{\Lambda}- v^{\Sigma}_{\Lambda}\|_{L^2 \cap L^{\infty}}=O(\varepsilon t^{-\frac{3}{4}+\sigma})$ and $\|v^{\Sigma}_{\Lambda}-v^{\Sigma}\|_{L^2\cap L^{\infty}}=O(\varepsilon t^{-\frac{1}{4}+\sigma})$, so we can deduce from \eqref{asymptotics for fLambda} the asymptotic expansion for $v^{\Sigma}=\langle hD\rangle^{-1}v$.
Since \eqref{Op(a)v development} for $a\equiv 1$ shows that $v^{\Sigma}$ vanishes when $x \notin [-1,1]$ and $t\rightarrow +\infty$, we get that $\tilde{a}(x)$ is supported for $x \in [-1,1]$, and we conclude the proof choosing $\tilde{a}(x)=\varepsilon a_{\varepsilon}(x)$ for a bounded $a_{\varepsilon}(x)$ as in the statement.
\endproof
\end{prop}

\subsection{End of the Proof}
\proof[Proof of Theorem \ref{main theorem}]
Let us prove that, for small enough data, the solution of the initial Cauchy problem \eqref{KG} is global.
We show that we can propagate some convenient \emph{a priori} estimates on $u$, as stated in theorem \ref{bootstrap theorem}, namely we want to show that there are some integers $s\gg \rho \gg 1$, some constants $A, B>0$ large enough, $\varepsilon_0 \in ]0,1]$ and $\sigma>0$ small enough such that, if $u \in C^0([1, T[; H^{s+1})\cap C^1([1, T[; H^s)$ is solution of \eqref{KG} for some $T>1$, and satisfies
\begin{align*}
& \|u(t,\cdot)\|_{W^{t, \rho, \infty}}  \le A \varepsilon t^{-\frac{1}{2}}\,, \\
& \|Zu(t,\cdot)\|_{H^1}  \le B\varepsilon t^{\sigma} \, , \qquad  \|\partial_tZu(t,\cdot)\|_{L^2}\le B \varepsilon t^{\sigma} \\
& \|u(t,\cdot)\|_{H^s} \le B\varepsilon t^{\sigma} \, , \qquad \hspace{0.3cm}\|\partial_t u(t,\cdot)\|_{H^{s-1}}\le B \varepsilon t^{\sigma}\, ,
\end{align*}
for every $t\in [1,T]$, then in the same interval it verifies improved estimates,
\begin{align*}
& \|u(t,\cdot)\|_{W^{t, \rho, \infty}}  \le \frac{A}{2} \varepsilon t^{-\frac{1}{2}} \,,\\
& \|Zu(t,\cdot)\|_{H^1}  \le \frac{B}{2}\varepsilon t^{\sigma} \, ,\qquad  \|\partial_tZu(t,\cdot)\|_{L^2}\le \frac{B}{2} \varepsilon t^{\sigma} \\
& \|u(t,\cdot)\|_{H^s} \le 	\frac{B}{2}\varepsilon t^{\sigma}\,, \qquad \hspace{0.3cm}\|\partial_t u(t,\cdot)\|_{H^{s-1}}\le \frac{B}{2} \varepsilon t^{\sigma} \,.
\end{align*}
We can immediately observe that from \eqref{conditions on data}, these bounds are verified at time $t=1$.
In theorem \ref{Propagation of Energy Estimates} in section \ref{section Generalised energy estimates}, we proved that we can improve the energy bounds $\|Zu(t,\cdot)\|_{H^1}$, $\|\partial_tZu(t,\cdot)\|_{L^2}$, $\|u(t,\cdot)\|_{H^s}$ and $\|\partial_t u(t,\cdot)\|_{H^{s-1}}$.
To show that the propagation of the uniform bound $\|u(t,\cdot)\|_{W^{t,\rho,\infty}}$ holds, we passed from equation \eqref{KG} to \eqref{half KG in w with polynomials} at the beginning of section \ref{Semiclassical Reduction to an ODE}, and then we showed that the function $v$ is solution of \eqref{half KG for v with polynomials}. 
The \emph{a priori} assumptions made on $u$ imply the following estimates on $v$,
\begin{equation} \label{end. a priori est on v}
\begin{split}
& \|v(t,\cdot)\|_{W^{\rho-1,\infty}_h}  \le C_1 A\varepsilon \, ,\\
&\|\mathcal{L}v(t,\cdot)\|_{H^1_h}  \le 5 B\varepsilon h^{-\sigma}\,, \qquad \|v(t,\cdot)\|_{H^s_h}  \le B\varepsilon h^{-\sigma} \,,
\end{split}
\end{equation}
for $h^{-1}:=t$ in $[1, T]$.
In fact, from \eqref{def of w}, the definition \eqref{def of v} of $v$ in semiclassical coordinates and the equation \eqref{KG}, $$C_2\|u(t,\cdot)\|_{W^{t,\rho,\infty}}\le t^{-\frac{1}{2}}\|v(t,\cdot)\|_{W_h^{\rho-1,\infty}} \le  C_1 \|u(t,\cdot)\|_{W^{t,\rho,\infty}} \,, $$ $$\|v(t,\cdot)\|_{H^s_h}=\|u(t,\cdot)\|_{H^s} \, ,$$ 
for some positive constants $C_1, C_2$, so the first and third inequality in \eqref{end. a priori est on v} are satisfied.
Moreover, $\mathcal{L}v$ can be expressed in term of $u$, $Zu$, as showed below using equation \eqref{half KG for v with polynomials}, 
\begin{equation} 
\begin{split}
\frac{1}{i}Zu(t,y) & = h^{\frac{1}{2}} \left[ (1-x^2)D_x + tx D_t + i\frac{x}{2} \right]v(t,x) |_{x=\frac{y}{t}}\\
& = \left(h^{\frac{1}{2}} \left[ (1-x^2)D_x + t x \, Op_h^w(x\xi + p(\xi)) + i \frac{x}{2} \right]v + h^{\frac{1}{2}} x \widetilde{P}\right)|_{x=\frac{y}{t}} \\
& = \left(h^{\frac{1}{2}} \left[ D_x + tx \, Op_h^w(p(\xi)) \right]v + h^{\frac{1}{2}} x \widetilde{P}\right)|_{x=\frac{y}{t}} \, ,
\end{split}
\end{equation}
where $\widetilde{P}$ denotes the right hand side of \eqref{half KG for v with polynomials} multiplied by $h^{-1}$.
Using symbolic calculus of proposition \ref{a sharp b},
\begin{equation} \label{relation between Zu and Lv}
\begin{split}
\frac{1}{i}Zu(t,y) & = \left(h^{\frac{1}{2}} \left[ h^{-1}Op_h^w(xp(\xi) + \xi) - \frac{1}{2i}Op_h^w(p'(\xi)) \right] v + h^{\frac{1}{2}} x \widetilde{P} \right)|_{x=\frac{y}{t}}\\
& = \left(h^{\frac{1}{2}} \left[Op_h^w(p(\xi))\mathcal{L}v - \frac{1}{i}Op_h^w(p'(\xi))v + x\widetilde{P} \right] \right)|_{x=\frac{y}{t}}\,, 
\end{split}
\end{equation}
where we used that $p(\xi)= \sqrt{1+\xi^2}$, $p'(\xi)=\xi/p(\xi)$.
Therefore, since $Op_h^w(p(\xi)^{-1}): H^{s-1}_h \rightarrow H^s_h$ is uniformly bounded by proposition \ref{Continuity from $L^2$ to $L^2$}, and from $\|v(t,\cdot)\|_{H^s_h}=\|u(t,\cdot)\|_{H^s}$, we derive
$\|\mathcal{L}v(t,\cdot)\|_{H^1_h}\le \|Zu(t,\cdot)\|_{L^2} + \|u(t,\cdot)\|_{L^2} + \|x\widetilde{P}\|_{L^2}$, where $$\|x \widetilde{P} (t,\cdot)\|_{L^2}\le C \|v(t,\cdot)\|^2_{W^{\rho-1,\infty}_h}(\|\mathcal{L}v(t,\cdot)\|_{L^2} + \|u(t,\cdot)\|_{H^s})\, .$$
Then we can use the uniform estimate $\|v(t,\cdot)\|_{W_h^{\rho-1,\infty}}\le C_1 A\varepsilon$, choose $\varepsilon_0 \ll 1$ small enough such that $CC_1A^2\varepsilon_0^2<\frac{1}{2}$, and use the \emph{a priori} energy bounds on $u$ in \eqref{bootstrap hypothesis}, to have 
$$\|\mathcal{L}v(t,\cdot)\|_{H^1_h}\le 2\|Zu(t,\cdot)\|_{L^2} + 2\|u(t,\cdot)\|_{L^2} + \|u(t,\cdot)\|_{H^s} \le 5B\varepsilon h^{-\sigma}\, .$$
Under these bounds on $v$, in proposition \ref{Propagation of the uniform estimate} we proved that, for $A'=C_1A$ and $B'=5B$, the uniform estimate on $v$ can be propagated, choosing for instance $K=\frac{2C_1}{C_2}$ to obtain $\|v(t,\cdot)\|_{W_h^{\rho-1,\infty}}\le \frac{A C_2}{2}\varepsilon$, and then $\|u(t,\cdot)\|_{W^{t,\rho,\infty}}\le \frac{A}{2}\varepsilon t^{-\frac{1}{2}}$, which concludes the proof of the boostrap and of global existence.

\vspace{0.5cm}
\noindent We prove now the asymptotics.
We consider $\Sigma(\xi)=\langle\xi\rangle^{\rho +1}$ and we write $$\langle hD\rangle^{-1}v = Op_h^w(\langle\xi\rangle^{-1} \langle \xi \rangle^{-\rho -1})v^{\Sigma} \, .$$
Using proposition \ref{prop Op(a) develop }, we develop the symbol $\langle \xi \rangle^{-\rho-2}$ at $\xi = d\varphi(x)$, 
\begin{equation*}
Op_h^w(\langle \xi \rangle^{-\rho-2})v^{\Sigma} = \theta_h(x) \langle d\varphi(x) \rangle^{-\rho-2} v^{\Sigma} + O_{L^{\infty}\cap L^2}(\varepsilon h^{\frac{1}{4}-\sigma}) \, ,
\end{equation*}
and using the expression obtained in \eqref{asymptotics for v}, along with the uniform bound on $v^{\Sigma}$, we derive that in the limit $t \rightarrow +\infty$ the function $\tilde{a}(x)=\varepsilon a_{\varepsilon}(x)$ verifies
\begin{equation}
|\tilde{a}(x)| \le |\theta_h(x) \langle d\varphi(x)\rangle^{-\rho-2}v^{\Sigma}|  + O(\varepsilon t^{-\frac{1}{4}+\sigma}) 
 \overset{t \rightarrow +\infty}{\le} C \varepsilon \langle d\varphi(x) \rangle^{-\rho-2} \, .
\end{equation}
For points $x$ in $]-1,1[$ such that $\langle d\varphi(x) \rangle \ge \alpha h^{-\beta}$, for a small $\alpha>0$, we have $|\tilde{a}(x)| = O(\varepsilon h^{\beta(\rho + 2)})$ and then the corresponding contribution to the right hand side of \eqref{asymptotics for v} is $O(\varepsilon t^{-\min(\beta(\rho +2), \frac{1}{4}-\sigma)})$ in $L^{\infty}\cap L^2$.

\vspace{0.5cm}
\noindent Let us now consider points $x$ in $]-1,1[$ such that $\langle d\varphi(x) \rangle \le \alpha h^{-\beta}$, and remind that the function $\theta_h(x)$ in \eqref{asymptotics for v} is identically equal to one on some interval $[-1+ch^{2\beta}, 1-ch^{2\beta}]$. 
We can write
\begin{equation} \label{decomposition of integral of theta}
\int_1^t \theta_{1/\tau}(x) d\tau = t-1 + \int_1^{\infty} (\theta_{1/\tau}(x) -1) d\tau - \int_t^{\infty} (\theta_{1/\tau}(x) -1) d\tau \, ,
\end{equation}
observing that on the support of $\theta_{1/\tau}(x)-1$, $\tau < \max c^{\frac{1}{2\beta}} (1-x, x+1)^{-\frac{1}{2\beta}}$.
Therefore the last integral is taken on a finite interval and since $|x \pm 1|\sim \langle d\varphi(x)\rangle^{-2}$ as $x\rightarrow \mp 1$ by \eqref{bound for x pm 1}, this implies that at the same time we have $\tau \le c \langle d\varphi(x) \rangle^{\frac{1}{\beta}}$ and $\langle d\varphi(x)\rangle^{\frac{1}{\beta}}\le \alpha t$.
For $t\le \tau$ and $\alpha>0$ small, this leads to a contradiction and to the fact that the last integral in \eqref{decomposition of integral of theta} is equal to zero. 
Then in \eqref{asymptotics for v} we can write
$$
a_{\varepsilon}(x) \exp\left[i\varphi(x) \int_1^t\theta_{1/\tau}(x) d\tau \right]= a_{\varepsilon}(x) \exp[i\varphi(x)t + ig(x)] \,,
$$
with $
g(x)= \varphi(x) \left[\int_1^{\infty}(\theta_{1/\tau}(x)-1) d\tau -1 \right]
$, and similarly, for $x$ satisfying $\langle d\varphi(x) \rangle\le \alpha h^{-\beta}$,
$$|a_{\varepsilon}(x)|^2 \Phi^{\Sigma}_1(x) \int_1^t \theta_{1/\tau}(x) \frac{d\tau}{\tau} = |a_{\varepsilon}(x)|^2 \Phi^{\Sigma}_1(x) \log t + \tilde{g}(x) \,, $$
for $\tilde{g}(x) = |a_{\varepsilon}(x)|^2 \Phi^{\Sigma}_1(x) \left[\int_1^{\infty}(\theta_{1/\tau}(x)-1) \frac{d\tau}{\tau} -1 \right]$.
Moreover, for $\langle d\varphi(x) \rangle \le \alpha h^{-\beta}$ the coefficient $a_{(1,1,-1)}(x)$ appearing in $\Phi^{\Sigma}_1(x)$ is equal to $\langle d\varphi(x)\rangle^{-1}$, since $\chi(h^{\beta}d\varphi(x))\gamma(\frac{x+p'(d\varphi(x))}{\sqrt{h}})\equiv 1$ if $\alpha$ is chosen sufficiently small, which implies that $\Phi_1^{\Sigma}(x)$ is exactly $\Phi_1(x)$ introduced in \eqref{def of Phi1}.
Modifying the function $a_{\varepsilon}(x)$ by a factor of modulus one, we derive from \eqref{asymptotics for v} the asymptotic behaviour for $\langle hD\rangle^{-1}v$:
\begin{equation}
\langle hD\rangle^{-1}v = \varepsilon a_{\varepsilon}(x) \exp \left[i\varphi(x)t + i(\log t)\varepsilon^2|a_{\varepsilon}(x)|^2\Phi_1(x) \right] + t^{-\theta} r(t,x) \, ,
\end{equation} 
for some $\theta>0$ and $\|r(t,\cdot)\|_{L^{\infty}}=O(\varepsilon)$, and reminding the relationship between $v$ and $w$ in \eqref{def of v}, and  between $w$ and $u$ in \eqref{def of w}, we finally obtain the asymptotics for $u$ in \eqref{asymptotics for u}.

\endproof

\section*{Appendix}

This appendix is devoted to the detailed proof of proposition \ref{a sharp b} and lemma \ref{lem of composition with Gamma Tilde}, which are technical.

\proof[Proof of Proposition \ref{a sharp b}] 
Let us expand $a(x+z, \xi + \zeta)$ at $(x,\xi)$ with Taylor's formula :
\begin{equation*}
\begin{split}
a(x+z, \xi + \zeta) = a(x,\xi) & + \sum_{\substack{\alpha= (\alpha_1, \alpha_2)\\
                  1\le|\alpha|\le k}} \frac{1}{\alpha !}\, \partial_x^{\alpha_1}\partial_{\xi}^{\alpha_2}a(x,\xi) z^{\alpha_1}\zeta^{\alpha_2}\\
& + \sum_{\substack{\beta= (\beta_1, \beta_2)\\
                  |\beta|= k+1}} \frac{k+1}{\beta !}\, z^{\beta_1}\zeta^{\beta_2}\int_0^1\partial_x^{\beta_1}\partial_{\xi}^{\beta_2}a(x+tz,\xi+t\zeta) (1-t)^k\, dt \, ,
\end{split}
\end{equation*}                  
and replace this development in \eqref{a sharp b integral formula}, obtaining :
\begin{equation*}
\begin{split}
a\sharp b & = \frac{1}{(\pi h)^2}\int_{\mathbb{R}^4}e^{\frac{2i}{h}(\eta z - y \zeta)} a(x,\xi)b(x+y,\xi + \eta)\, dy d\eta dz d\zeta \\
& + \frac{1}{(\pi h)^2}\int_{\mathbb{R}^4}e^{\frac{2i}{h}(\eta z - y \zeta)}\sum_{\substack{\alpha= (\alpha_1, \alpha_2)\\
                  1\le|\alpha|\le k}}\frac{1}{\alpha !}\, \partial_x^{\alpha_1}\partial_{\xi}^{\alpha_2}a(x,\xi) b(x+y, \xi + \eta)\, z^{\alpha_1}\zeta^{\alpha_2} \, dy d\eta dz d\zeta \\
& + \frac{1}{(\pi h)^2}\int_{\mathbb{R}^4}e^{\frac{2i}{h}(\eta z -y \zeta)} \bigg\{\sum_{\substack{\beta= (\beta_1, \beta_2)\\
                  |\beta|= k+1}} \frac{k+1}{\beta !}\, z^{\beta_1}\zeta^{\beta_2}\int_0^1\partial_x^{\beta_1}\partial_{\xi}^{\beta_2}a(x+tz,\xi+t\zeta) (1-t)^k\, dt\bigg\}\\
& \hspace{4cm}\times b(x+y,\xi + \eta)\, dy d\eta dz d\zeta \\
& := I_1 + I_2 + I_3\,.
\end{split}
\end{equation*}
\renewcommand{\theequation}{A}
From a direct calculation and using that the inverse Fourier transform of the complex exponential is the delta function, i.e.
\begin{equation} \label{delta function}
\frac{1}{\pi h}\int_{\mathbb{R}}e^{\frac{2i}{h}XY} dY = \delta_0(X) \,,
\end{equation}
we derive
\begin{equation*}
\begin{split}
I_1& = \frac{1}{(\pi h)^2}\int_{\mathbb{R}^4}e^{\frac{2i}{h}(\eta z - y \zeta)} a(x,\xi)b(x+y,\xi + \eta)\, dy d\eta dz d\zeta \\
&= a(x,\xi) \int_{\mathbb{R}^2}b(x+y, \xi + \eta) \delta_0(y) \delta_0(\eta) \,dy d\eta = a(x,\xi)b(x,\xi) \,,
\end{split}
\end{equation*}
and 
\begin{equation*}
\begin{split}
& I_2 = \\
& = \frac{1}{(\pi h)^2}\sum_{\substack{\alpha= (\alpha_1, \alpha_2)\\
                  1\le|\alpha|\le k}}\frac{1}{\alpha !} \int_{\mathbb{R}^4}e^{\frac{2i}{h}(\eta z - y \zeta)}\, \partial_x^{\alpha_1}\partial_{\xi}^{\alpha_2}a(x,\xi) b(x+y, \xi + \eta)\, z^{\alpha_1}\zeta^{\alpha_2} \, dy d\eta dz d\zeta \\
& = \frac{1}{(\pi h)^2}\sum_{\substack{\alpha= (\alpha_1, \alpha_2)\\
                  1\le|\alpha|\le k}}\frac{1}{\alpha !}\left(\frac{h}{2i}\right)^{|\alpha|} \int_{\mathbb{R}^4} \partial_{\eta}^{\alpha_1}(-\partial_y^{\alpha_2})e^{\frac{2i}{h}(\eta z - y \zeta)}\, \partial_x^{\alpha_1}\partial_{\xi}^{\alpha_2}a(x,\xi) b(x+y, \xi + \eta)\, dy d\eta dz d\zeta \\
                  & =  \frac{1}{(\pi h)^2}\sum_{\substack{\alpha= (\alpha_1, \alpha_2)\\
                  1\le|\alpha|\le k}}\frac{(-1)^{\alpha_1}}{\alpha !}\left(\frac{h}{2i}\right)^{|\alpha|} \int_{\mathbb{R}^4} e^{\frac{2i}{h}(\eta z - y \zeta)}\, \partial_x^{\alpha_1}\partial_{\xi}^{\alpha_2}a(x,\xi) \, \partial_y^{\alpha_2}\partial_{\eta}^{\alpha_1} b(x+y, \xi + \eta)\, dy d\eta dz d\zeta \\
& = \sum_{\substack{\alpha= (\alpha_1, \alpha_2)\\
                  1\le|\alpha|\le k}}\frac{(-1)^{\alpha_1}}{\alpha !}\left(\frac{h}{2i}\right)^{|\alpha|} \partial_x^{\alpha_1}\partial_{\xi}^{\alpha_2}a(x,\xi)\, \partial_x^{\alpha_2}\partial_{\xi}^{\alpha_1} b(x, \xi)\, .               
\end{split}                  
\end{equation*}
The same calculation shows that $I_3$ is given by
\begin{equation*}
\begin{split}
I_3 = \, \frac{k+1}{(\pi h)^2} \left(\frac{h}{2i}\right)^{k+1} \sum_{\substack{\alpha= (\alpha_1, \alpha_2)\\
|\alpha|= k+1}} \frac{(-1)^{\alpha_1}}{\alpha!}\int_{\mathbb{R}^4}e^{\frac{2i}{h}(\eta z - y\zeta)} \Big\{& \int_0^1  \partial_x^{\alpha_1}\partial_{\xi}^{\alpha_2}a(x+tz, \xi +t\zeta)(1-t)^k dt  \\
&  \times \partial_y^{\alpha_2}\partial_{\eta}^{\alpha_1}b(x+y, \xi + \eta)\Big\} \, dy d\eta dz d\zeta \,,
\end{split}
\end{equation*}
and it belongs to $h^{(k+1)(1-(\delta_1+\delta_2))}S_{\delta,\beta}(M_1M_2)$ since
\begin{equation*}
\begin{split}
& \frac{1}{h^2}\int_{\mathbb{R}^4}e^{\frac{2i}{h}(\eta z - y\zeta)} \Big\{ \int_0^1  \partial_x^{\alpha_1}\partial_{\xi}^{\alpha_2}a(x+tz, \xi +t\zeta)(1-t)^k dt   \partial_y^{\alpha_2}\partial_{\eta}^{\alpha_1}b(x+y, \xi + \eta)\Big\} \, dy d\eta dz d\zeta  = \\
& = \int_{\mathbb{R}^4}e^{2i(\eta z - y\zeta)} \Big\{ \int_0^1  \partial_x^{\alpha_1}\partial_{\xi}^{\alpha_2}a(x+t\sqrt{h} z, \xi +t\sqrt{h} \zeta)(1-t)^k dt\, \partial_y^{\alpha_2}\partial_{\eta}^{\alpha_1}b(x+\sqrt{h} y, \xi +\sqrt{h} \eta)\Big\} \, \\
& \hspace{14cm} dy d\eta dz d\zeta \\
& = \int_{\mathbb{R}^4}\left(\frac{1+2iy\partial_{\zeta}}{1+4y^2}\right)^N\left(\frac{1-2i\eta\partial_z}{1+4\eta^2}\right)^N\left(\frac{1-2iz\partial_{\eta}}{1+4z^2}\right)^N\left(\frac{1+2i\zeta\partial_y}{1+4\zeta^2}\right)^Ne^{2i(\eta z - y\zeta)} \\
& \times \Big\{ \int_0^1  \partial_x^{\alpha_1}\partial_{\xi}^{\alpha_2}a(x+t\sqrt{h} z, \xi +t\sqrt{h} \zeta)(1-t)^k dt\, \partial_y^{\alpha_2}\partial_{\eta}^{\alpha_1}b(x+\sqrt{h} y, \xi +\sqrt{h} \eta)\Big\} \, dy d\eta dz d\zeta \\
& \mbox{so integrating by parts,} \\
& \le C h^{-(\delta_1 + \delta_2)(\alpha_1 + \alpha_2)} \int_{\mathbb{R}^4} \langle y \rangle^{-N} \langle \eta \rangle^{-N} \langle z \rangle^{-N} \langle \zeta \rangle^{-N} \Big\{\int_0^1 M_1(x+t\sqrt{h}z, \xi + t\sqrt{h}\zeta) dt \\
& \hspace{3.5cm}\times M_2(x+\sqrt{h}y, \xi + \sqrt{h}\eta) \Big\}\, dy d\eta dz d\zeta \\
& \le C h^{-(\delta_1 + \delta_2)(k+1)} \int_{\mathbb{R}^4} \langle y \rangle^{-N+N_0} \langle \eta \rangle^{-N+N_0} \langle z \rangle^{-N+N_0} \langle \zeta \rangle^{-N+N_0} \, dy d\eta dz d\zeta\,  M_1(x,\xi) M_2(x,\xi) \\
& \le C h^{-(\delta_1 + \delta_2)(k+1)} M_1(x,\xi)M_2(x,\xi) \, .
\end{split}
\end{equation*}
Equivalently, one can show that $|\partial^{\alpha}I_3|\le C h^{(k+1)(1-(\delta_1+\delta_2)) - \delta |\alpha|}M_1(x,\xi)M_2(x,\xi)$.
The last statement of the proposition follows immediately if we replace in previous inequalities $M_1$ and $M_2$ respectively by $M_1^{k+1}$, $M_2^{k+1}$.
\endproof

\vspace{1cm}
\proof[Proof of Lemma \ref{lem of composition with Gamma Tilde}]
The proof of the lemma is the same as the previous one if, when we calculate to which class the remainder $r_k$ belongs, we remark that
\begin{equation*}
\left\langle \frac{x+t\sqrt{h}z+ f(\xi + t\sqrt{h}\zeta)}{\sqrt{h}}\right\rangle^{-d} = \left\langle \frac{x+f(\xi)}{\sqrt{h}} + tz +tb(\xi,\zeta)\zeta  \right\rangle^{-d} \lesssim \langle tz \rangle^N \langle t\zeta \rangle^N \left\langle \frac{x+f(\xi)}{\sqrt{h}}\right\rangle^{-d}
\end{equation*}
\begin{equation*}
\begin{split}
\left\langle \frac{x+\sqrt{h}y+ f(\xi + \sqrt{h}\eta)}{\sqrt{h}}\right\rangle^{-l}= \left\langle \frac{x+f(\xi)}{\sqrt{h}} + y +b'(\xi,\eta)\eta \right\rangle^{-l} \lesssim \langle y\rangle^N \langle\eta\rangle^N \left\langle \frac{x+f(\xi)}{\sqrt{h}}\right\rangle^{-l}
\end{split}
\end{equation*}
with $b(\xi,\zeta)=\int_0^1f'(\xi + st\sqrt{h}\zeta)ds \lesssim 1$, $b'(\xi,\eta)=\int_0^1 f'(\xi + s\sqrt{h}\eta)ds \lesssim 1$, for a certain $N \in \mathbb{N}$.
\endproof

\bibliographystyle{abbrv}
\bibliography{bibliography}

\end{document}